\documentclass[10pt]{article}
\usepackage{framed,multirow}
\usepackage{amssymb}
\usepackage{latexsym}
\usepackage{url}
\usepackage{xcolor}
\usepackage{setspace}
\definecolor{newcolor}{rgb}{.8,.349,.1}
\makeatletter
\def\ps@pprintTitle{%
	\let\@oddhead\@empty
	\let\@evenhead\@empty
	\def\@oddfoot{}%
	\let\@evenfoot\@oddfoot}
\makeatother
\usepackage[margin=1.5cm]{geometry}
\usepackage[titletoc]{appendix}

\usepackage[normalem]{ulem}

\usepackage{float}
\usepackage{amssymb}
\usepackage{amsmath}
\usepackage{amsthm}
\usepackage{graphics}
\usepackage{psfrag}
\usepackage{graphicx}
\usepackage{color}
\usepackage{fancyhdr}
\usepackage{appendix}
\usepackage{verbatim}
\usepackage{blindtext}
\usepackage{hyperref}
\usepackage{cleveref}
\usepackage{stmaryrd}
\usepackage{soul}
\usepackage{sidecap}
\usepackage{tikz}
\usepackage{rotating}
\usepackage{booktabs}
\usepackage{comment}
\usepackage{todonotes}
\usepackage{multirow}
\usepackage{changepage}
\usepackage{float}
\usepackage{wrapfig}
\usepackage{stmaryrd}
\usepackage{subcaption}  
\usepackage{todonotes}
\usepackage{natbib}
\setcitestyle{numbers}
\setcitestyle{square}
\crefformat{appendix}{#2#1#3}

\definecolor{lightblue}{rgb}{0.32,0.45,0.90}
\definecolor{lightgreen}{rgb}{0.42,0.7,0.40}
\numberwithin{equation}{section}
\numberwithin{figure}{section}

\hypersetup{colorlinks=true,linkcolor = blue,citecolor=blue}

\usetikzlibrary{shapes,arrows}
\usetikzlibrary{arrows.meta}
\usetikzlibrary{decorations.markings}
\numberwithin{figure}{section}

\renewcommand{\div}{\operatorname*{div}}

\def\b{\boldsymbol}

\makeatletter
\newcommand{\vast}{\bBigg@{4}}
\newcommand{\Vast}{\bBigg@{5}}
\makeatother\def\b{\boldsymbol}

\colorlet{cgray}{gray!20!white}

\theoremstyle{definition}
\newtheorem{definition}{Definition}[section]
\newtheorem{assumptions}{Assumptions}[section]

\newtheorem{theorem}{Theorem}[section]
\theoremstyle{remark}
\newtheorem*{remark}{Remark}

\newtheorem{lemma}{Lemma}[section]

\tikzset{every label/.style={font=\footnotesize,inner sep=1pt}}

\usepackage{pgfplots}
\pgfplotsset{compat = newest}
\usepgfplotslibrary{fillbetween,decorations.softclip}

\allowdisplaybreaks

\begin{document}

\title{New Limiter Regions for Multidimensional Flows}

\author{James Woodfield$^{1,2}$, Hilary Weller$^3$, Colin J Cotter$^1$}

\addtocounter{footnote}{1} 
\footnotetext{Department of Mathematics, Imperial College London, South Kensington Campus, London SW7 2AZ, United Kingdom.}
\addtocounter{footnote}{1}
\footnotetext{Department of Mathematics and Statistics, University of Reading, Whiteknights, Reading RG6 6AX, United Kingdom.}
\addtocounter{footnote}{1} 
\footnotetext{Department of Meteorology, University of Reading, University of Reading, Earley Gate, Reading RG6 6BB, United Kingdom.}

\date{\today}

\maketitle

  



\begin{abstract}

Accurate transport algorithms are crucial for computational fluid dynamics and more accurate and efficient schemes are always in development. One dimensional limiting is commonly employed to suppress nonphysical oscillations. However, the application of such limiters can reduce accuracy. It is important to identify the weakest set of sufficient conditions required on the limiter as to allow the development of successful numerical algorithms.

The main goal of this paper is to identify new less restrictive sufficient conditions for flux form in-compressible advection to remain monotonic. We identify additional necessary conditions for incompressible flux form advection to be monotonic, demonstrating that the Spekreijse limiter region is not sufficient for incompressible flux form advection to remain monotonic. Then a convex combination argument is used to derive new sufficient conditions that are less restrictive than the Sweby region for a discrete maximum principle. This allows the introduction of two new more general limiter regions suitable for flux form incompressible advection. 
\end{abstract}

\section{Introduction}
\subsection{Historical context and motivation}

Higher order accuracy and monotonicity have long been at odds. Godunov \cite{godunov:hal-01620642} showed that the order of a linear monotonicity preserving numerical scheme is necessarily first order. Harten, Hyman, Lax and Keyfitz (HHLK)\cite{harten1976finite} generalised the Godunov order barrier theorem, and proved that nonlinear schemes with the stronger HHLK-monotone property are necessarily first order. Harten \cite{harten1984class} developed nonlinear total variation diminishing (TVD) schemes, satisfying a monotonic property (TVD) weaker than the HHLK-monotone property but stronger than the monotonicity preserving property of Godunov \cite{godunov:hal-01620642}. These TVD schemes bypass the Godunov order barrier, can be second order and have a convenient framework and limiter region introduced by Sweby in 1984 \cite{S_1984}.

However, the total variation diminishing framework did not generalise easily for multi-dimensional flows, in part due to Goodman and Leveque showing that a particular two-dimensional definition of TVD would also run into the HHLK order barrier \cite{GL_1985}. Nevertheless, in 1987 Spekreijse \cite{S_1987} developed a monotonicity framework suitable for multidimensional flow, it is based on positivity of coefficients and can be used to enforce a discrete local maximum principle (strictly stronger than the total variation diminishing property). 
Furthermore Spekreijse’s \cite{S_1987} positive coefficient framework was applied to a wide class of flux splitting schemes and generalised to unstructured meshes. The corresponding spatial discretisations are referred to as positive coefficient schemes or Local Extrema Diminishing schemes \cite{jameson1995}. Spekreijse also showed that less stringent requirements were demanded of the limiter functions, these new weaker set of sufficient conditions required on the limiter function allowed a variety of new flux and slope limiters to be introduced. Of particular note is the development of differentiable limiters contained in the second order region.

Spekreijse’s work has since been used to motivate monotonic atmospheric advection:
Hundsdorfer et al. \cite{hundsdorfer1995positive} proposed using one-dimensional flux limiters in a multidimensional flux form finite difference spatial discretisation, for positive mass preserving transport on the sphere. Hundsdorfer et al. \cite{hundsdorfer1995positive} observed preferable accuracy and computational speed using one-dimensional limiters in a method of lines framework as compared with the third order MPDATA scheme in \cite{smolarkiewicz1991monotone}. Wesseling and Zijlema \cite{zijlema1998higher} proposes a very similar scheme under a finite volume interpretation, and suggest some limiters in the Spekreijse region could be used to improve accuracy. We will analyse a generalisation of both the numerical method introduced in Wesseling and Zijlema \cite{zijlema1998higher} and Hundsdorfer et al. \cite{hundsdorfer1995positive}, to find the admissible flux limiter region for the incompressible advection equation in flux form. It is known that using one-dimensional limiters in the Spekreijse region can create a one dimensional TVD scheme for uniform flow. It is also believed that the monotonicity properties carry to the multidimensional flux form advection equation case using the positive coefficient framework of Spekreijse \cite{S_1987}, some examples of this belief can be found in \cite{koren1993robust,zijlema1998higher,hundsdorfer1995positive}. However, it has never been proven and to what extent that this is true has never been resolved.

Whether Spekreijse's limiter region is applicable to multidimensional flux form transport schemes is important as limiters such as OSPRE, van Albada, Albada family, Hemker-Koren, local double-logarithmic reconstruction (LDLR), Čada–Torrilhon, TCDF, MPL2-$\kappa$ and ENO2 may fail to have a provable discrete maximum principle and consequently could lose positivity which could potentially lead to instability when transporting dense scalars such as liquid water. The region in which limiters are provably discrete maximum principle satisfying is crucial to the successful design of monotone advection algorithms and is the main focus of this paper.

In \cref{sec:Background material and summary: Different existing limiter regions} we define and plot both the Sweby Region, and the Spekreijse limiter Region, in the context of one dimensional constant coefficient advection as to put the rest of the paper in context with the existing literature. In \cref{sec:Governing equation and numerical scheme}, we introduce the advection equation and the flux form numerical method of interest, a generalisation of the finite volume method in Wesseling and Zijlema \cite{zijlema1998higher} and the finite difference method in Hundsdorfer et al. \cite{hundsdorfer1995positive}. In \cref{sec:Applicability of the extended Spekreijse region} we explain why Spekreijse's framework does not necessarily apply to flux form advection schemes, and why the extended admissible limiter region of Spekreijse can fail. 
\Cref{sec:theory} contains the main contribution of this paper, the derivation of two new limiter regions sufficient for the incompressible flux form scheme defined in \cref{sec:semi discrete scheme} to remain discrete maximum principle satisfying. 

\Cref{Other theoretical properties of the scheme} contains a wider discussion on the theoretical properties of the numerical scheme introduced in \cref{sec:semi discrete scheme}, in particular how the method of lines interacts with monotonicity properties, a review of the spatial accuracy a short remark on linear invariance and the introduction of several new limiter functions. In \cref{sec:Numerical Demonstrations} we include several numerical tests to see to what extent the theoretical results in \cref{sec:Applicability of the extended Spekreijse region} and \cref{sec:theory} are observed in practice. In \cref{sec:Numerical Demonstrations} we numerically investigate the limiters introduced in \cref{Other theoretical properties of the scheme}.

\subsection{Background: Different existing limiter regions} \label{sec:Background material and summary: Different existing limiter regions}
A semi discrete flux limited numerical scheme for the constant coefficient $v \in \mathbb{R}^{>0}$ advection equation $u_{t}+v u_{x}=0$ can be written as $\frac{du_{i}}{dt} + \Delta x ^{-1}(F_{i+1/2} -F_{i-1/2})=0$, where the flux is given by $F_{i+1/2} = v[u_i + 1/2\psi(R_{i})(u_{i}-u_{i-1})]$, here $\psi$ is the limiter function, and the ratio of successive gradients is given by $R_{i} =(u_{i+1}-u_{i})/(u_{i}-u_{i-1})$. A similar but different construction is given using the inverse ratio  $r_{i} = 1/R_{i}$ with flux given by $F_{i+1/2} = v[u_i +1/2\psi(r_{i})(u_{i+1}-u_{i})]$. We distinguish these different frameworks by the parameter $\theta = 1, 0$, respectively. The total variation of $u$ is defined as $TV(u)= \sum_{\forall i} |u_{i} - u_{i-1}|$, and serves as a measure of how oscillatory $u$ is. A desirable property of a numerical method is that the total variation does not increase with time. This property is referred to as total variation diminishing (TVD).\newline

In 1984, Sweby \citep{S_1984} introduced sufficient conditions for a flux limited scheme in the $\theta = 0$ formulation to be total variation diminishing,
\begin{align}
	\psi(r) &\leq \min(2r,2), \quad \forall r>0, \label{eq:superbee bound condition}\\
	\psi(r) &= 0, \quad \forall r<0. \label{eq:flux limiter zero at extrema}
\end{align}
These conditions serve as bounds on acceptable limiters in the Sweby diagram, and we shall refer to the set of points $(r,\psi(r))$ for which conditions \cref{eq:flux limiter zero at extrema,eq:superbee bound condition} are satisfied as the Sweby region. However, it is sometimes assumed that these bounds are necessary for schemes to be total variation diminishing \citep{versteeg2007introduction}.

Spekreijse showed that this is not the case in the $\theta=1$ formulation by finding a more general (TVD) admissible limiter region,
\begin{align}
	\psi(R) \in[\alpha,M], \quad \psi(S)/S \in [-M,2+\alpha],\quad M\in (0,\infty), \quad  \alpha \in [-\infty,0], \label{eq: spekreijse limiter region}
\end{align} 
for $R,S \in \mathbb{R}$, similarly we shall refer to the set of points $(R,\psi(R))$ for which the conditions in \cref{eq: spekreijse limiter region} are satisfied as the Spekreijse region.
Spekreijse showed that flux splitting schemes with limiters in this region satisfy a local discrete maximum principle, which is strictly stronger than the total variation diminishing property in one dimension and is more convenient in multiple dimensions. The use of limiters in the extended region of Spekreijse has been successful in both flux difference splitting \cite{deconinck1993multidimensional} and flux vector splitting frameworks \cite{S_1987}. \Cref{fig: Regions} shows both the Sweby region $\mathcal{D}_{1}$ in \cref{fig:sweby}, and the Spekreijse region $\mathcal{D}_{2}$ in \cref{fig:spek}. The extended limiter region has allowed a huge variety of new limiters to be introduced into the literature, which can improve accuracy \cite{hua1992possible}, and are not necessarily reduced to the first order upwind scheme at extrema and can be globally smooth \cite{S_1987}. It is worth noting, the $\theta = 0$ framework also possesses a two-parameter generalisation \cite{jacques1989uniformly,hua1992possible}, which we shall also refer to as the Spekreijse region, or more generally an extended limiter region. In the flux splitting framework, there also exists other extended limiter regions, for a selective literature review see \cref{sec:Relationship to other limiters}. As remarked earlier, the extent to which these extended limiter regions are suitable for flux form incompressible advection has not been resolved.

\begin{figure}
\centering
	\begin{subfigure}{0.4\textwidth}
		\begin{tikzpicture}[scale=0.9]
			\pgfdeclarelayer{pre main}
			\pgfsetlayers{pre main,main}
			\begin{axis}[
				axis lines = middle,
				axis equal,
				domain  = -1:3,
				xlabel  = {$r$},
				ylabel  = {$y$},
				xmin    = -0.5,
				xmax    = 2.9,
				ymin    = -0.5,
				ymax    = 2.9,
				samples = 100,
				mark    = none,
				]
				\pgfonlayer{pre main}
				\endpgfonlayer
				\addplot [name path=twor, dashed]  {(2*x)};
				\addplot [name path=zero, thick]  {(0)};
				\addplot [name path=ninf, thin]  {(-2)};
				\addplot [name path=inf, thin]  {(4)};
				\addplot [name path=two, dashed]  {(2)};
				\addplot[color=black!30] fill between[of=zero and twor,
				soft clip={domain=0:1}];
				\addplot[color=black!30] fill between[of=zero and two,
				soft clip={domain=1:3.9}];
				\node at (0.7,0.5) {$\mathcal{D}_{1}$};
				\node at (1.8,2.5) {$y = 2r$};
			\end{axis}
		\end{tikzpicture}
		\caption{Sweby}
		\label{fig:sweby}
	\end{subfigure}\qquad 
	\begin{subfigure}{0.4\textwidth}
		\begin{tikzpicture}[scale=0.9]
			\pgfdeclarelayer{pre main}
			\pgfsetlayers{pre main,main}
			\begin{axis}[
				axis lines = middle,
				axis equal,
				domain  = -2:3,
				xlabel  = {$R$},
				ylabel  = {$y$},
				xmin    = -2,
				xmax    = 2.9,
				ymin    = -1,
				ymax    = 2.9,
				samples = 100,
				mark    = none,
				yticklabels={,,}
				xticklabel=empty,
				ticks=none
				]
				\pgfonlayer{pre main}
				\endpgfonlayer
				\addplot [name path=mMr, dashed]  {(-1.5*x)};
				\addplot [name path=twoalphar, dashed]  {((2-0.5)*x)};
				\addplot [name path=zero, thin]  {(0)};
				\addplot [name path=ninf, thin]  {(-4)};
				\addplot [name path=inf, thin]  {(4)};
				\addplot [name path=alpha, thick]  {(-0.5)};
				\addplot [name path=M, thick]  {(1.5)};
				\addplot[color=black!30] fill between[of=mMr and twoalphar,
				soft clip={domain=0:0.333}];
				\addplot[color=black!30] fill between[of=alpha and twoalphar,
				soft clip={domain=0.333:1}];
				\addplot[color=black!30] fill between[of=alpha and M,
				soft clip={domain=1:3}];
				\addplot[color=black!30] fill between[of=mMr and twoalphar,
				soft clip={domain=-0.333:0}];
				\addplot[color=black!30] fill between[of=mMr and alpha,
				soft clip={domain=-1:-0.333}];
				\addplot[color=black!30] fill between[of=M and alpha,
				soft clip={domain=-2:-1}];
				\node at (0.8,0.5) {$\mathcal{D}_{2}$};
				\node at (1.9,2.5) {$y = (2+\alpha)R$};
				\node at (2.4,1.7) {$y= M$};
				\node at (-1.1,2.5) {$y= -M R$};
				\node at (2.4,-0.7) {$y = \alpha$};
			\end{axis}
		\end{tikzpicture}
		\caption{Spekreijse}
		\label{fig:spek}
	\end{subfigure}
	\caption{ \cref{fig:sweby} is the plot of the sufficient admissible limiter region in the $(r,y)$ plane as defined by Sweby \cite{S_1984}. \Cref{fig:spek} is the plot of the admissible limiter region in the $(R,y)$ plane as defined by Spekreijse \cite{S_1987}. 
	The Sweby region $\mathcal{D}_{1}$ is sufficient for a one-dimensional scheme with flux limiters to be TVD, the Spekreijse region $\mathcal{D}_{2}$, which has two free parameters $\alpha \in (-\infty,0]$, $M \in (0,\infty)$, is also sufficient for the scheme to be TVD. }
	\label{fig: Regions}
\end{figure}

%

%


\section{Governing equation and numerical scheme} \label{sec:Governing equation and numerical scheme}

\subsection{Governing equation}

We consider the evolution of a bounded initial tracer $u_{0}(\b x) \in L^{\infty}(\Omega)$ in a arbitrary spatial domain, $\Omega \subseteq \mathbb{R}^{d}, d\geq 1$, over the time interval $[0,T]$ by a bounded $\b v(\b x, t) \in L^{\infty}(\Omega \times [0,T])$ divergence free $\div (\b v(\b x, t)) = 0$ vector field. This is described by 
the advection equation,
\begin{align}
	\frac{\partial u}{\partial t} + \div (\b v u) &= 0, \quad \forall (\b x,t)\in \Omega\times[0,T], \label{eqn:Transport equation} \\
	u(\b x,0) &= u_{0}(\b x), \quad \forall \b x \in \Omega,
\end{align}
which is a model for the evolution of densities, and for incompressible flows, scalars.

\subsection{Semi discrete scheme}\label{sec:semi discrete scheme}
We introduce the notation and describe a generalisation of the flux form finite volume scheme of Wesseling and Zijlema \cite{zijlema1998higher} and the finite difference scheme of Hundsdorfer et al \cite{hundsdorfer1995positive}.


\begin{enumerate}
	\item Let $u = u_{i,j}$ denote either the cell mean or pointwise value of a tracer within a cell, where if one or other subscript is missing it is assumed to be at position $i$ or $j$ as appropriate.
	\item Reconstruct: A reconstruction operator constructs the values attained to the right, left, up and down of cell $(i,j)$ as follows
	\begin{align}
	&u^{R}_{i} = u_{i} +  \frac{\theta}{2}\psi(R_{i})(u_{i}-u_{i-1}) +  \frac{(1-\theta) }{2}\psi(\frac{1}{R_{i}})(u_{i+1}-u_{i}), \label{eq:UR}\\
	&u_i^L = u_i  +  \frac{\theta}{2} \psi(\frac{1}{R_{i}})(u_i - u_{i+1}) -   \frac{(1-\theta) }{2}\psi(R_{i})(u_{i} - u_{i-1}), \label{eq:UL} \\
	&u_j^U = u_j + \frac{\theta}{2}  \psi(R_j)(u_j - u_{j-1}) +  \frac{(1-\theta) }{2}\psi(\frac{1}{R_{j}})(u_{j+1}-u_{j}),\label{eq:UU}\\
	&u_j^D = u_j  + \frac{\theta}{2}  \psi(\frac{1}{R_{j}})(u_j - u_{j+1}) -   \frac{(1-\theta) }{2}\psi(R_{j})(u_{j} - u_{j-1}), \label{eq:UD} \\
	& R_{i}  = \frac{u_{i+1}-u_{i}}{u_{i}-u_{i-1}}, \quad R_{j}  = \frac{u_{j+1}-u_{j}}{u_{j}-u_{j-1}},
\end{align}
this reconstruction imposes upwind bias flux limiting. The right, left, up and down of cell $(i,j)$ are denoted with $R,L,U,D$ superscripts. We are using the ratio of successive gradients defined in \cite{S_1987} as $R$ when $\theta =1$, we call this the Roe gradient. When $\theta =0$ we are using the inverse ratio defined in \cite{S_1984} $r_{i} = 1/R_{i}$, we call the Sweby gradient. The choice of $\theta$ does not matter when the limiter is symmetric $\psi(1/R) = \psi(R)/R$. We only investigate $\theta = \lbrace 0,1 \rbrace$ due to the reduced computational cost in evaluating the expressions in \cref{eq:UD,eq:UL,eq:UR,eq:UU}. When the symmetry condition is broken the two limiter functions should be distinguished differently, as they have their own respective limiter regions, however this should be apparent from the context of the limiter, so we do not include this notationally.


\item Riemann: the donor cell numerical flux function 
	\begin{align}
&F(u^{R}_{i},u^{L}_{i+1},c_{i+0.5} )  = c_{i+0.5}^{+}u^{R}_{i} + c_{i+0.5}^{-}u^{L}_{i+1},
	\end{align}
	resolves the Riemann problems at the cell boundaries. Here we define the notation $(\cdot)^{+} := \max (\cdot,0)$, $(\cdot)^{-} := \min (\cdot,0)$, to mean the positive or negative component of the argument respectively, to be used throughout this paper.
	$c_{i+1/2}$ denotes the $x$ component of the velocity at position $(i+1/2,j)$, and $c_{j+1/2}$ denotes the $y$ component of the velocity at position $(i,j+1/2)$. We often absorb constants like timestep and mesh width into this constant because the numerical flux function is linear. This numerical flux function is consistent, monotone, and Lipschitz continuous as defined in \cite{EGH_2000}. This type of flux is known as a state interpolated flux as defined in \cite{koren1993robust, hundsdorfer1994method,hundsdorfer1995positive}, rather than a flux interpolated flux. In the finite volume setting we approximate the integral of flux through the face by second order Gauss quadrature over each face (midpoint rule), whereas in the finite difference setting we interpret as a point valued flux. 
	\item Evolve: the semi discrete evolution operator is given by the flux form method
	\begin{align}
	\begin{split}
		\frac{\partial u}{\partial t} + [F_{i+0.5}(u^{R}_{i},u^{L}_{i+1},c_{i+0.5}) - F_{i-0.5}(u^{R}_{i-1},u^{L}_{i},c_{i-0.5})]+  [F_{j+0.5}(u^{U}_{j},u^{D}_{j+1},c_{j+0.5}) - F_{j-0.5}(u^{U}_{j-1},u^{D}_{j},c_{j-0.5})] = 0.
	\end{split}
	\label{eq:EVOLVE_METHOD}
	\end{align}
We have absorbed the mesh spacing into the coefficients $c$ by linearity of the donor cell flux function, and plan to do so with the time-step $\Delta t$.
\end{enumerate}
For context see \cite{S_1987},\cite{zijlema1998higher},\cite{hundsdorfer1995positive}, and for an example of its generalisation to unstructured meshes see \cite{wilders2002positive}.
The above scheme is a flux form scheme with velocities defined at cell edges, see \cref{fig:flux form stencil}. In contrast, an advective form scheme uses a cell-centred velocity $c_{i},c_{j}$ (see \cref{fig: advective form stencil}) such that the scheme has the following form
\begin{align}
\frac{\partial u}{\partial t} +  [F_{i+0.5}(u^R_{i},u^{L}_{i+1},c_i)-F_{i-0.5}(u^R_{i-1},u^{L}_{i},c_i)] + [F_{j+0.5}(u^R_{j},u^{L}_{j+1},c_j)-F_{j-0.5}(u^R_{j-1},u^{j}_{i},c_j)].\label{eq:advective form}
\end{align}
The advective form scheme is provably monotone using the Spekreijse's limiter region and framework \cite{S_1987}. However, this argument is not always valid in the flux form multidimensional case, as will be explained in \cref{sec:Applicability of the extended Spekreijse region}. We present the theoretical results regarding flux form scheme in the next section.  


\begin{figure}
\centering
\begin{subfigure}[b]{0.30\textwidth}
\begin{tikzpicture}
\draw[fill=gray!5,line width=1] (5,3) rectangle (7,5);
\draw[fill=gray!20,line width=1] (5,5) rectangle (7,7);
\draw[fill=gray!20,line width=1] (5,1) rectangle (7,3);
\draw[fill=gray!20,line width=1] (3,3) rectangle (5,5);
\draw[fill=gray!20,line width=1] (7,3) rectangle (9,5);
\draw[-stealth,line width=0.25mm] (7,4)--(7.75,4);
\draw[-stealth,line width=0.25mm] (6,5)--(6,5.75);
\node at (6.7,4)  {\small{$u^{R}_{i,j}$}};
\node at (6,4.2) {\small{$u_{i,j}$}};
\node at (5.3,4) {\small{$u^{L}_{i,j}$}};
\node at (6,4.75) {\small{$u^{U}_{i,j}$}};
\node at (6,3.3) {\small{$u^{D}_{i,j}$}};
\filldraw (6,4) circle[radius=1.5pt];
\filldraw (8,4) circle[radius=1.5pt];
\filldraw (4,4) circle[radius=1.5pt];
\filldraw (6,6) circle[radius=1.5pt];
\filldraw (6,2) circle[radius=1.5pt];
\node at (8,4.2) {\small{$u_{i+1,j}$}};
\node at (4,4.2) {\small{$u_{i-1,j}$}};
\node at (6,6.2) {\small{$u_{i,j+1}$}};
\node at (6,2.2) {\small{$u_{i,j-1}$}};
\node [draw, diamond, fill, scale=0.125] at (6,3) {$ok$};
\node [draw, diamond, fill, scale=0.125] at (6,5) {$ok$};
\node [draw, diamond, fill, scale=0.125] at (7,4) {$ok$};
\node [draw, diamond, fill, scale=0.125] at (5,4) {$ok$};
\node at (6.5,5.25) {\small{$c_{j+1/2}$}};
\node at (7.5,3.75) {\small{$c_{i+1/2}$}};
\end{tikzpicture}
\caption{Flux form stencil.}
\label{fig:flux form stencil}
\end{subfigure}
\hspace{20mm}
\begin{subfigure}[b]{0.30\textwidth}\centering
\begin{tikzpicture}
\draw[fill=gray!20,line width=1] (5,3) rectangle (7,5);
\draw[dotted,-stealth,line width=0.25mm] (6,4)--(7.75,4);
\draw[dotted,-stealth,line width=0.25mm] (7.75,4)--(7.75,5.25);
\draw[-stealth,line width=0.25mm] (6,4)--(7.75,5.25);
\filldraw (6,4) circle[radius=1.5pt];
\node at (6.7,4)  {\small{$u^{R}_{i,j}$}};
\node at (6,4.2) {\small{$u_{i,j}$}};
\node at (5.3,4) {\small{$u^{L}_{i,j}$}};
\node at (6,4.75) {\small{$u^{U}_{i,j}$}};
\node at (6,3.3) {\small{$u^{D}_{i,j}$}};
\node [draw, diamond, fill, scale=0.125] at (6,3) {$ok$};
\node [draw, diamond, fill, scale=0.125] at (6,5) {$ok$};
\node [draw, diamond, fill, scale=0.125] at (7,4) {$ok$};
\node [draw, diamond, fill, scale=0.125] at (5,4) {$ok$};
\node at (8.12,4.5) {\small{$c_{j}$}};
\node at (7.25,3.75) {\small{$c_{i}$}};
\end{tikzpicture}
\caption{Advective form stencil.}
\label{fig: advective form stencil}
\end{subfigure}
\caption{ \Cref{fig:flux form stencil} Shows variable placement for a flux form scheme, consistent with numerical method \cref{EQ:SemiDiscrete Method}, where the velocity field is located at cell faces/edges. \Cref{fig: advective form stencil} shows variable placement for an advective form scheme, consistent with \cref{eq:advective form}, where the velocity is located at the cell centre. One can use Spekreijse's limiter region using the advective form stencil, but not when using the flux form stencil (\cref{thm:necessary}). When using a flux form scheme, one can instead use the limiter regions presented in this paper \cref{region:Woodfield}, \cref{thm: divergence}.}
\end{figure}

\section{Theoretical monotonicity properties of the scheme}\label{sec:theoretical monotonicity properties of the scheme}

We present and derive less restrictive sufficient conditions on the limiters functions $\psi$ to prove that the semidiscrete numerical method (\cref{EQ:SemiDiscrete Method})
\begin{equation}
  \begin{aligned}
&\frac{\partial u}{\partial t} + [F_{i+0.5}(u^{R}_{i},u^{L}_{i+1},c_{i+0.5}) - F_{i-0.5}(u^{R}_{i-1},u^{L}_{i},c_{i-0.5})]+  [F_{j+0.5}(u^{U}_{j},u^{D}_{j+1},c_{j+0.5}) - F_{j-0.5}(u^{U}_{j-1},u^{D}_{j},c_{j-0.5})] = 0.\\
&F(u^{R}_{i},u^{L}_{i+1},c_{i+0.5} )  = c_{i+0.5}^{+}u^{R}_{i} + c_{i+0.5}^{-}u^{L}_{i+1},\quad F(u^{R}_{i-1},u^{L}_{i},c_{i-0.5} )  = c_{i-0.5}^{+}u^{R}_{i-1} + c_{i-0.5}^{-}u^{L}_{i},\\
&F(u^{U}_{j},u^{D}_{j+1},c_{j+0.5} )  = c_{j+0.5}^{+}u^{U}_{j} + c_{j+0.5}^{-}u^{D}_{j+1},\quad F(u^{U}_{j-1},u^{D}_{j},c_{j-0.5} )  = c_{j-0.5}^{+}u^{U}_{j-1} + c_{j-0.5}^{-}u^{D}_{j-1}, \\
&u^{R}_{i} = u_{i} +  \frac{\theta}{2}\psi(R_{i})(u_{i}-u_{i-1}) +  \frac{(1-\theta) }{2}\psi(\frac{1}{R_{i}})(u_{i+1}-u_{i}), \\
&u_i^L = u_i  +  \frac{\theta}{2} \psi(\frac{1}{R_{i}})(u_i - u_{i+1}) -   \frac{(1-\theta) }{2}\psi(R_{i})(u_{i} - u_{i-1}),  \\
&u_j^U = u_j + \frac{\theta}{2}  \psi(R_j)(u_j - u_{j-1}) +  \frac{(1-\theta) }{2}\psi(\frac{1}{R_{j}})(u_{j+1}-u_{j}),\\
&u_j^D = u_j  + \frac{\theta}{2}  \psi(\frac{1}{R_{j}})(u_j - u_{j+1}) -   \frac{(1-\theta) }{2}\psi(R_{j})(u_{j} - u_{j-1}),  \\
& R_{i}  = \frac{u_{i+1}-u_{i}}{u_{i}-u_{i-1}}, \quad R_{j}  = \frac{u_{j+1}-u_{j}}{u_{j}-u_{j-1}}, 
\end{aligned}\label{EQ:SemiDiscrete Method}
\end{equation}

is Local Extrema Diminishing (LED) for incompressible flow for $\theta \in \lbrace 0,1 \rbrace$ (We introduced all necessary notation and explanation of this scheme in \cref{sec:semi discrete scheme}). 
This can be used to identify less strict conditions in which the forward Euler numerical flow map $u^{n+1} = \text{Forward Euler} ( u^n, c^n), $
defined by
\begin{align}
	\begin{split}
u_{i,j}^{n+1}  =  u_{i,j}^n &- [F_{i+0.5}(u^{R}_{i},u^{L}_{i+1},c_{i+0.5}^n) - F_{i-0.5}(u^{R}_{i-1},u^{L}_{i},c_{i-0.5}^n)] - [F_{j+0.5}(u^{U}_{j},u^{D}_{j+1},c_{j+0.5}^n) - F_{j-0.5}(u^{U}_{j-1},u^{D}_{j},c_{j-0.5}^n)], \label{eq:FE}
\end{split}
\end{align}
can satisfy a Discrete Local Maximum Principle (DLMP).
The time step, $\Delta t$, and mesh spacing, $\Delta x$, can be absorbed into the face defined velocity field $c$, as and when appropriate. This scheme generalises the finite volume scheme in \cite{zijlema1998higher}, and the finite difference scheme in \cite{hundsdorfer1995positive}, and can be used with a Strong Stability Preserving timestepping scheme to achieve higher order accuracy whilst retaining monotonic properties.

\subsection{New limiter regions suitable for incompressible flow} \label{sec:theory}

In this section we will show that under new considerations of the limiter function we can put the scheme in a positive coefficient type form \cref{Def:Positive Coefficient} as defined by Spekreijse \cite{S_1987}. This definition allows use of \cref{Lemma:spekreijse} also given by Spekreijse \cite{S_1987} which can be used to prove a discrete local maximum principle for the forward Euler flow map. And also allows the use of a Lemma by Jameson \cite{jameson1995} used to prove the the Local Extrema Diminisishing (LED) property of the semi-discrete scheme.

\begin{definition}[Spekreijse, Positive-coefficient type scheme \cite{S_1987}]\label{Def:Positive Coefficient} The semidiscrete scheme given by 
	\begin{align}
		\label{eq:semidiscrete form}
		\begin{split}
		\frac{\partial  u_{i,j}}{ \partial t} &+ A_{i-1/2}(u_{i} - u_{i-1})+ B_{i+1/2}(u_{i} - u_{i+1})+C_{j-1/2}(u_{j} - u_{j-1})+ D_{j+1/2}(u_{j} - u_{j+1}) = 0, 
		\end{split}
	\end{align}
	is a positive coefficient scheme when all the nonlinear leading coefficients are non-negative
	\begin{align}
 		A_{i-1/2}, B_{i+1/2}, C_{j-1/2}, D_{j+1/2}\geq 0. \label{ineq:Positive coefficients}
	\end{align}
\end{definition}

\begin{lemma} [Jameson \cite{jameson1995}]\label{Lemma:Jameson LED}
Semi-discrete positive coefficient type schemes (satisfying \cref{Def:Positive Coefficient}) are local extrema diminishing \cite{jameson1995}
\end{lemma}
\begin{proof}
 If $u_{i,j}\geq u_{i+1},u_{i-1},u_{j+1},u_{j-1}$, then $\partial_t u_{i,j} \leq 0$, if $u_{i,j}\leq u_{i+1},u_{i-1},u_{j+1},u_{j-1}$, then $\partial_t u_{i,j} \geq  0$, local maxima do not increase and local minima do not decrease, this is the definition of Local Extrema Diminishing (LED). 
\end{proof}

\begin{lemma}[Spekreijse~\cite{S_1987}]\label{Lemma:spekreijse}

When the forward Euler temporal discretisation of a positive coefficient type scheme  satisfies the time step restriction,
	\begin{align}
		\left(A^{n}_{i-1/2}+B^{n}_{i+1/2}+C^{n}_{j-1/2}+ D^{n}_{j+1/2}\right)\leq 1, \label{timestepping restriction}
	\end{align} 
	then $u^{n+1}_{i,j}$ is a convex combination of $u_{i+1}^n,u_{i-1}^n,u^n,u_{j+1}^n,u_{j-1}^n$, trivially implying the following local discrete maximum principle with respect to edge sharing neighbours cell mean (or cell centre pointwise) values, 
	\begin{align}
		u^{n+1} \in [ \min \lbrace u_{i+1}^n,u_{i-1}^n,u^n,u_{j+1}^n,u_{j-1}^n \rbrace, \max \lbrace u_{i+1}^n,u_{i-1}^n,u^n,u_{j+1}^n,u_{j-1}^n \rbrace]. \label{eq:discrete maximum principle}
	\end{align}
	This in turn this implies weaker properties such as boundedness $||u^{n+1}||_{\infty} \leq ||u^{n}||_{\infty}$, and sign preservation. 
\end{lemma}
\begin{proof}
The forward Euler discretisation of \cref{eq:semidiscrete form} can be written as 
\begin{align}
	\begin{split}
		u^{n+1}_{i,j} &= (1- [A^{n}_{i-1/2}+B^{n}_{i+1/2}+C^{n}_{j-1/2}+ D^{n}_{j+1/2}])u_{i,j} + A^{n}_{i-1/2}u_{i-1}+ B^{n}_{i+1/2} u_{i+1}+C^{n}_{j-1/2}u_{j-1}+ D^{n}_{j+1/2}u_{j+1}. 
	\end{split}
	\label{eq: convex combo}
\end{align}
Which is a convex combination of neighbouring cells under the condition \cref{timestepping restriction}. 
\end{proof}
\begin{remark}
The above formulation of the maximum principle is pessimistic, in practice the upwind biasing imposes a more physically motivated maximum principle by setting some of the coefficients $A^{n}_{i-1/2},B^{n}_{i+1/2},C^{n}_{j-1/2},D^{n}_{j+1/2}$ to zero depending on the direction of flow. 
\end{remark}

\begin{theorem}[In-compressible flow] \label{thm: divergence} The Forward Euler discretisation of the numerical method described in \cref{sec:semi discrete scheme}, can be written as convex combination of neighbour cell mean values \cref{eq: convex combo} and as a result satisfies the discrete maximum principle \cref{eq:discrete maximum principle} when the \cref{assumptions for my theorem} hold. 
\end{theorem}

\begin{assumptions}
[\cref{thm: divergence} assumptions]
\label{assumptions for my theorem}
\phantom{ok}\newline
\begin{enumerate} 
	\item The mesh scaled velocity satisfies a discrete divergence free condition 
	\begin{align}c_{i+1/2} - c_{i-1/2}  + c_{j+1/2} -c_{j-1/2} = 0,\quad  \forall (i,j).
	\end{align} 
	\item When $\theta = 1$, the limiter function satisfies
\begin{align}
	\psi(R) \in[0,M_{\psi}], \quad \psi(S)/S \in [m_{\psi} , 2 ]. \quad M_{\psi}\in [0,\infty), \quad m_{\psi}\in (-\infty, 0]. \label{eq: new limiter region}
\end{align} 
When $\theta = 0$, the limiter function satisfies
\begin{align}
\psi(1/R)\in [m_{\psi},2], \quad  S\psi(\frac{1}{S}) \in [0,M_{\psi}], \quad M_{\psi}\in [0,\infty),\quad m_{\psi}\in (-\infty, 0].
\end{align}
\item The timestep restriction
\begin{align}
C:= c_{i+1/2}^{+}-c_{i-1/2}^{-}+c_{j+1/2}^{+}-c_{j-1/2}^{-} \leq C_{FE}:= \frac{2}{2+M_{\psi}- m_{\psi}}
\end{align}
  holds (written in terms of the Courant number). 
\end{enumerate}
\end{assumptions}

\begin{proof} of \cref{thm: divergence}.
Expand the numerical method
in terms of their positive and negative components, and taking away the divergence free condition gives 
\begin{align}
\begin{split}
	\frac{\partial u}{\partial t}& + c_{i+0.5,j}^{+}(u^{R}_{i}-u_{i})+c_{i+0.5}^{-}(u^{L}_{i+1} -u_{i}) -c_{i-0.5}^{+}(u^{R}_{i-1}-u_{i}) -c_{i-0.5}^{-}(u^{L}_{i} -u_{i})\\& + c_{j+0.5}^{+}(u^{U}_{j}-u_{i}) +c_{j+0.5}^{-}(u^{D}_{j+1}-u_{i}) -c_{j-0.5}^{+}(u^{U}_{j-1}-u_{i})  -c_{j-0.5}^{-}(u^{D}_{j}-u_{i}) = 0.
\end{split}
\end{align}
We have rearranged into the following positive coefficient representation only using the definition of $R$
	\begin{align}
	\begin{split}
		\frac{\partial u}{\partial t} &+ \bigg( c_{i+1/2}^{+}\big[ \frac{\theta \psi( R_{i} ) }{2}  + \frac{(1-\theta)}{2} R_{i}\psi( \frac{1}{R_{i}} ) \big] + c_{i-1/2}^{+}\big[ 1 - \frac{\theta \psi(R_{i-1}}){2 R_{i-1}} - \frac{(1-\theta)}{2}\psi(\frac{1}{R_{i-1}}) \big] \bigg)[u_{i}-u_{i-1}]\\
		& + \bigg( -c^{-}_{i+1/2}\big[ 1-\frac{\theta}{2}\psi(\frac{1}{R_{i+1}})R_{i+1} - \frac{1-\theta}{2}\psi(R_{i+1}) \big] - c_{i-1/2}^{-}\big[ \frac{\theta}{2}\psi(\frac{1}{R_{i}}) + \frac{1-\theta}{2} \frac{\psi(R_{i})}{R_{i}}\big] \bigg)[u_{i}-u_{i+1}]\\
		& + \text{y-direction}. 
	\end{split}
	\end{align}

When $\theta = 1$
the conditions
\begin{align}
0 \leq \psi(R)  \leq M_{\psi}, \quad m_{\psi} \leq \psi(S)/S \leq 2, \quad \forall R,S \in \mathbb{R},
\end{align} are sufficient for the following bounds
\begin{align}
\begin{split}
	A_{i-1/2} &\in [0, c_{i+0.5}^{+}M_{\psi}/2 + c_{i-0.5}^{+} (1-m_{\psi}/2) ] , \\
	B_{i+1/2} & \in [0 , -c_{i+0.5}^{-}[1 - m_{\psi}/2 ]-c_{i-0.5}^{-} M_{\psi}/2] ,\\
	C_{j-1/2} &\in [0, c_{j+0.5}^{+}M_{\psi}/2 + c_{j-0.5}^{+} (1- m_{\psi}/2) ] , \\
	D_{j+1/2} & \in [0 , -c_{j+0.5}^{-}[1 - m_{\psi}/2 ]-c_{j-0.5}^{-} M_{\psi}/2],
\end{split}
\end{align}
for some $M_{\psi}\geq 0$, and $m_{\psi} \leq 2$. This shows the semi-discrete scheme is written in a positive coefficient representation sufficient for \cref{Lemma:Jameson LED,Lemma:spekreijse}. When $\theta =0$, the similar conditions
\begin{align}
	m_{\psi}\leq \psi(1/R) \leq 2, \quad 0 \leq S\psi(1/S) \leq M_{\psi}.
\end{align} are also sufficient for similar bounds. 
A sufficient timestep restriction $(1- [A^{n}_{i-1/2}+B^{n}_{i+1/2}+C^{n}_{j-1/2}+ D^{n}_{j+1/2}])\leq 1$ for both representations to maintain a discrete maximum principle is $C\leq \frac{2}{2+M_{\psi}-m_{\psi}},$ and appears in \cite{S_1984}. The conditions on the timestep in this proof are slightly more convenient than they are sharp, for additional details see \cref{sec:Details of proof}.
\end{proof}

\begin{remark}
We have yet to put some additional sensible constraints on the limiter functions region. $m_{\psi} \leq 0$ is a design principle that should be adhered to, otherwise $\psi(S) \geq m_{\psi}S$, runs into conflict with the conditions $\psi(R) \leq M_{\psi}$, $\psi(R)\leq 0$, when $M_{\psi}$ is finite and $R$ is large. Good design requires passing through $\psi(1)=1$ for second order accuracy, so that $M_{\psi}\in [1,\infty)$ is a sensible construction. We also impose $\psi(0)=0$ as we have not found a case of a limiter function not satisfying this condition. 
\end{remark}

Previously we considered sufficient conditions on the limiter function for the forward Euler scheme to have a discrete maximum principle.  The next theorem pertains to some necessary conditions on the limiter for the same forward Euler scheme to be positivity preserving.

\begin{theorem}[Necessary conditions]\label{thm:necessary}
As can be shown through contradiction using case-by-case considerations of the velocity field in the forward Euler representation one can deduce that
\begin{align}
   \text{when}\quad \theta=1,\quad  \psi(R)\geq 0, \frac{\psi(S)}{S}\leq 2,\quad \forall R\in (-\infty,0],\quad \forall S \in (0,\infty)\\
   \text{when}\quad \theta=0,\quad \psi(T)\leq 2, \frac{\psi(S)}{S} \geq 0 \quad \forall S\in (-\infty,0], \quad \forall T\in (0,\infty)
\end{align} are necessary conditions for the positivity preservation of the forward Euler scheme. These are plotted in \cref{region:Necessary}.
\end{theorem}
\begin{proof}
Case-by-case considerations are presented in \cref{proof:necessary}. 
\end{proof}

\begin{remark}
In the $\theta=1$ framework, the tail (values of $\psi$, for $R<0$) has to be non-negative, in the $\theta=0$ framework the tail (values of $\psi$, for $r<0$) has to be non-positive. In the $\theta=1$ framework in the positive region $R\geq 0$, we must have $\psi(R)\leq 2R$, and for $\theta=0$ we must have $\psi(t)\leq 2$ in the positive region $t\geq 0$. \Cref{region:Necessary} contains a diagram of some of the necessary conditions.
\end{remark}

\begin{figure}
\centering
	\begin{subfigure}{0.475\textwidth}
		\begin{tikzpicture}[scale=1]
			\pgfdeclarelayer{pre main}
			\pgfsetlayers{pre main,main}
			\begin{axis}[
				axis lines = middle,
				axis equal,
				enlargelimits,
				domain  = -4:4,
				xlabel  = {$r$},
				ylabel  = {$y$},
				xmin    = -2,
				xmax    = 2.9,
				ymin    = -1,
				ymax    = 3.1,
				samples = 100,
				mark    = none,
				yticklabels={,,}
				ticks=none
				]
				\pgfonlayer{pre main}
				\endpgfonlayer
				\addplot[name path = sou,  thin,domain=-3:5, color = black ,dashed] {x};
				\addplot[name path = cds,  thin,domain=-3:5, color = black ,dashed] { 1};
				\addplot[mark=*] coordinates {(1,1)};
				\addplot [name path=zero, thin]  {(0)};
				\addplot [name path=ninf, thin]  {(-4)};
				\addplot [name path=inf, thin]  {(5)};
				\addplot [name path=mpsi, solid]  {(-0.5)};
				\addplot [name path=twor, thin,dashed]  {(2*x)};
				\addplot [name path=Ms, thin]  {(3*x)};
				\addplot [name path=two, solid]  {(2)};
                \addplot[color=black!10] fill between[of=mpsi and zero,
		        soft clip={domain=0:4}];
				\addplot[color=black!10] fill between[of=mpsi and zero,
				soft clip={domain=-4:-0.3333*0.5}];
				\addplot[color=black!10] fill between[of=Ms and zero,
				soft clip={domain=-0.3333*0.5:0}];
				\addplot[color=black!10] fill between[of= zero and Ms,
				soft clip={domain=0:1.3333*0.5}];
				\addplot[color=black!10] fill between[of= zero and two,
				soft clip={domain=1.3333*0.5:4}];
				\addplot[color=black!10] fill between[of=mpsi and zero,
				soft clip={domain=-4:-0.5}];
				\addplot[color=black!10] fill between[of=zero and sou,
				soft clip={domain=-0.5:0.0}];
				\addplot[color=black!10] fill between[of=sou and Ms,
				soft clip={domain=0.0:0.333333}];
				\addplot[color=black!10] fill between[of=sou and cds,
				soft clip={domain=0.333333:1}];
				\addplot[color=black!10] fill between[of=sou and cds,
				soft clip={domain=1:2}];
				\addplot[color=black!10] fill between[of=two and cds,
				soft clip={domain=2:5}];
				\node at (0.8,0.5) {$\mathcal{D}_{3}$};
				\node at (0.7,2.8) {$y = M_{\psi} r$};
				\node at (1.8,2.8) {$y =2 r$};
				\node at (-1.2,1.8) {$y= 2$};
				\node at (1.4,-0.75) {$y = m_{\psi}$};
			\end{axis}
		\end{tikzpicture}
		\caption{ $\theta = 0$}
	\end{subfigure}
	\quad
	\begin{subfigure}{0.475\textwidth}
		\begin{tikzpicture}[scale=1.0]
			\pgfdeclarelayer{pre main}
			\pgfsetlayers{pre main,main}
			\begin{axis}[
				axis lines = middle,
				axis equal,
				domain  = -3:5,
				xlabel  = {$R$},
				ylabel  = {$y$},
				xmin    = -2.5,
				xmax    = 4.9,
				ymin    = -0.5,
				ymax    = 3.9,
				samples = 300,
				mark    = none,
				]
				\pgfonlayer{pre main}
				\endpgfonlayer
				\addplot [name path=mMr, solid]  {(-2.9*x)};
				\addplot [name path=twor, solid]  {(2*x)};
				\addplot [name path=r, dashed]  {(x)};
				\addplot [name path=one, dashed]  {(1)};
				\addplot [name path=zero, thin]  {(0)};
				\addplot [name path=ninf, thin]  {(-2)};
				\addplot [name path=inf, thin]  {(6)};
				\addplot [name path=two, dashed]  {(2)};
				\addplot [name path=M, solid]  {(2.9)};
                \addplot[color=black!10] fill between[of=zero and ninf,
				soft clip={domain=0:6}];
    
				\addplot[color=black!10] fill between[of=zero and one,
				soft clip={domain=-3:-.5}];
				\addplot[color=black!10] fill between[of=mMr and zero,
				soft clip={domain=-0.5:0}];
				\addplot[color=black!10] fill between[of=one and M,
				soft clip={domain=-3:-1}];
				\addplot[color=black!10] fill between[of=zero and mMr,
				soft clip={domain=-1:0}];
				\addplot[color=black!10] fill between[of=zero and one,
				soft clip={domain=-1:-0.3448275862}];
				\addplot[color=black!10] fill between[of=zero and mMr,
				soft clip={domain=-0.3448275862:0}];
				\addplot[color=black!10] fill between[of=r and twor,
				soft clip={domain=0:0.5}];
				\addplot[color=black!10] fill between[of=r and one,
				soft clip={domain=0.5:1}];
				\addplot[color=black!10] fill between[of=one and twor,
				soft clip={domain=0.5:1.45}];
				\addplot[color=black!10] fill between[of=M and zero,
				soft clip={domain=1.45:5}];
				\addplot[color=black!10] fill between[of=zero and r,
				soft clip={domain=0:1}];
				\addplot[color=black!10] fill between[of=zero and one,
				soft clip={domain=1:5}];
				\addplot[color=black!10] fill between[of=one and r,
				soft clip={domain=1:2.9}];
				\addplot[color=black!10] fill between[of=one and M,
				soft clip={domain=2.9:5}];
				\addplot[mark=*] coordinates {(1,1)};
				\node at (0.8,0.5) {$\mathcal{D}_{4}$};
				\node at (2.2,4.3) {$y = 2R$};
				\node at (3.9,4.3) {$y = R$};
				\node at (0.5,3.0) {$y = M_{\psi}$};
				\node at (-1.3,4.3) {$y = m_{\psi}R$};
			\end{axis}
		\end{tikzpicture}
		\caption{ $\theta = 1$ }
	\end{subfigure}
\caption{It is necessary (but not proven sufficient), that limiters for incompressible flux form advection are contained in the above regions for positivity preservation. For $\theta=0$, $\theta = 1$, respectively where $M_{\psi}\in[0,\infty]$, $m_{\psi}\in[-\infty,0]$.} 
\label{region:Necessary}
\end{figure}
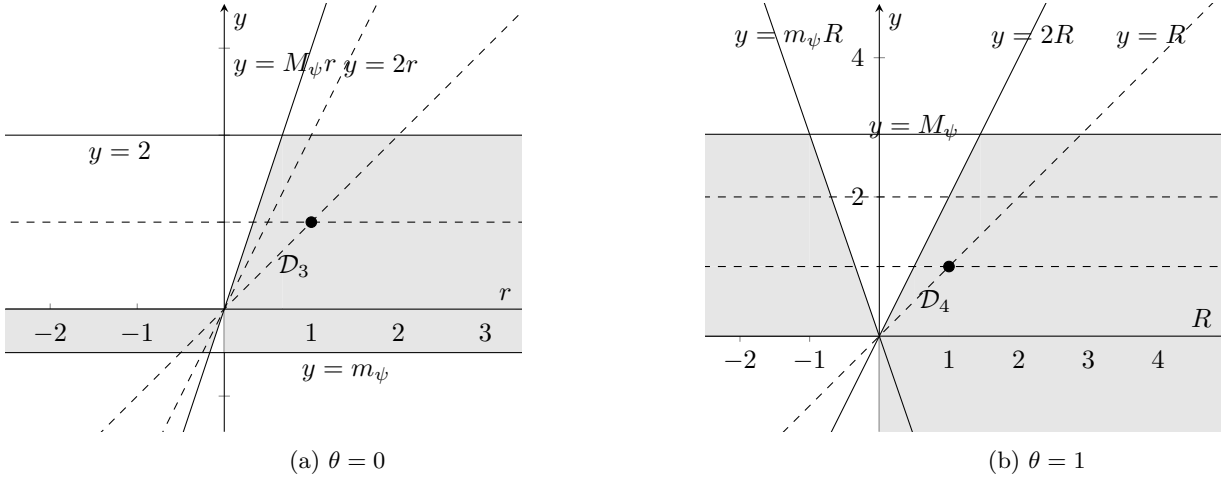

\begin{figure}
\centering
	\begin{subfigure}{0.475\textwidth}
		\begin{tikzpicture}[scale=1]
			\pgfdeclarelayer{pre main}
			\pgfsetlayers{pre main,main}
			\begin{axis}[
				axis lines = middle,
				axis equal,
				enlargelimits,
				domain  = -4:4,
				xlabel  = {$r$},
				ylabel  = {$y$},
				xmin    = -2,
				xmax    = 2.9,
				ymin    = -1,
				ymax    = 3.1,
				samples = 100,
				mark    = none,
				yticklabels={,,}
				ticks=none
				]
				\pgfonlayer{pre main}
				\endpgfonlayer

				\addplot[name path = sou,  thin,domain=-3:5, color = black ,dashed] {x};
				\addplot[name path = cds,  thin,domain=-3:5, color = black ,dashed] { 1};
				\addplot[mark=*] coordinates {(1,1)};
				\addplot [name path=zero, thin]  {(0)};
				\addplot [name path=ninf, thin]  {(-4)};
				\addplot [name path=inf, thin]  {(5)};
				\addplot [name path=mpsi, solid]  {(-0.5)};
				\addplot [name path=twor, thin,dashed]  {(2*x)};
                \addplot [name path=third order, thin,dotted] {(1/3*x+2/3)};
				\addplot [name path=Ms, thin]  {(3*x)};
				\addplot [name path=two, solid]  {(2)};
				\addplot[color=black!10] fill between[of=mpsi and zero,
				soft clip={domain=-4:-0.3333*0.5}];
				\addplot[color=black!10] fill between[of=Ms and zero,
				soft clip={domain=-0.3333*0.5:0}];
				\addplot[color=black!10] fill between[of= zero and Ms,
				soft clip={domain=0:1.3333*0.5}];
				\addplot[color=black!10] fill between[of= zero and two,
				soft clip={domain=1.3333*0.5:4}];
				\addplot[color=black!30] fill between[of=mpsi and zero,
				soft clip={domain=-4:-0.5}];
				\addplot[color=black!30] fill between[of=zero and sou,
				soft clip={domain=-0.5:0.0}];
				\addplot[color=black!30] fill between[of=sou and Ms,
				soft clip={domain=0.0:0.333333}];
				\addplot[color=black!30] fill between[of=sou and cds,
				soft clip={domain=0.333333:1}];
				\addplot[color=black!30] fill between[of=sou and cds,
				soft clip={domain=1:2}];
				\addplot[color=black!30] fill between[of=two and cds,
				soft clip={domain=2:5}];
				\node at (0.8,0.5) {$\mathcal{D}_{3}$};
				\node at (0.7,2.8) {$y = M_{\psi} r$};
				\node at (1.8,2.8) {$y =2 r$};
				\node at (-1.2,1.8) {$y= 2$};
				\node at (1.4,-0.75) {$y = m_{\psi}$};
			\end{axis}
		\end{tikzpicture}
		\caption{ $\theta = 0$}
	\end{subfigure}
	\quad
	\begin{subfigure}{0.475\textwidth}
		\begin{tikzpicture}[scale=1.0]
			\pgfdeclarelayer{pre main}
			\pgfsetlayers{pre main,main}
			\begin{axis}[
				axis lines = middle,
				axis equal,
				domain  = -3:5,
				xlabel  = {$R$},
				ylabel  = {$y$},
				xmin    = -2.5,
				xmax    = 4.9,
				ymin    = -0.5,
				ymax    = 3.9,
				samples = 300,
				mark    = none,
				]
				\pgfonlayer{pre main}
				\endpgfonlayer
				\addplot [name path=mMr, solid]  {(-2.9*x)};
				\addplot [name path=twor, solid]  {(2*x)};
				\addplot [name path=thirdorder, dotted, domain={-2:5}]  {(2/3*x+1/3)};
				\addplot [name path=r, dashed]  {(x)};
				\addplot [name path=one, dashed]  {(1)};
				\addplot [name path=zero, thin]  {(0)};
				\addplot [name path=ninf, thin]  {(-2)};
				\addplot [name path=inf, thin]  {(6)};
				\addplot [name path=two, dashed]  {(2)};
				\addplot [name path=M, solid]  {(2.9)};
				\addplot[color=black!30] fill between[of=zero and one,
				soft clip={domain=-3:-.5}];
				\addplot[color=black!30] fill between[of=mMr and zero,
				soft clip={domain=-0.5:0}];
				\addplot[color=black!10] fill between[of=one and M,
				soft clip={domain=-3:-1}];
				\addplot[color=black!10] fill between[of=zero and mMr,
				soft clip={domain=-1:0}];
				\addplot[color=black!30] fill between[of=zero and one,
				soft clip={domain=-1:-0.3448275862}];
				\addplot[color=black!30] fill between[of=zero and mMr,
				soft clip={domain=-0.3448275862:0}];
				\addplot[color=black!30] fill between[of=r and twor,
				soft clip={domain=0:0.5}];
				\addplot[color=black!30] fill between[of=r and one,
				soft clip={domain=0.5:1}];
				\addplot[color=black!10] fill between[of=one and twor,
				soft clip={domain=0.5:1.45}];
				\addplot[color=black!10] fill between[of=M and zero,
				soft clip={domain=1.45:5}];
				\addplot[color=black!10] fill between[of=zero and r,
				soft clip={domain=0:1}];
				\addplot[color=black!10] fill between[of=zero and one,
				soft clip={domain=1:5}];
				\addplot[color=black!30] fill between[of=one and r,
				soft clip={domain=1:2.9}];
				\addplot[color=black!30] fill between[of=one and M,
				soft clip={domain=2.9:5}];
				\addplot[mark=*] coordinates {(1,1)};
				\node at (0.8,0.5) {$\mathcal{D}_{4}$};
				\node at (2.2,4.3) {$y = 2R$};
				\node at (3.9,4.3) {$y = R$};
				\node at (0.5,3.0) {$y = M_{\psi}$};
				\node at (-1.3,4.3) {$y = m_{\psi}R$};
			\end{axis}
		\end{tikzpicture}
		\caption{ $\theta = 1$ }
	\end{subfigure}
\caption{New incompressible flow limiter regions, for $\theta=0,1$ respectively. Light gray regions are sufficient for a local maximum principle (bounded by solid lines), and dark grey indicates a desirable second order region (bounded by dashed lines). Light dotted line denotes a 3rd order region in the finite difference sense for uniform flow.}
\label{region:Woodfield}
\end{figure}
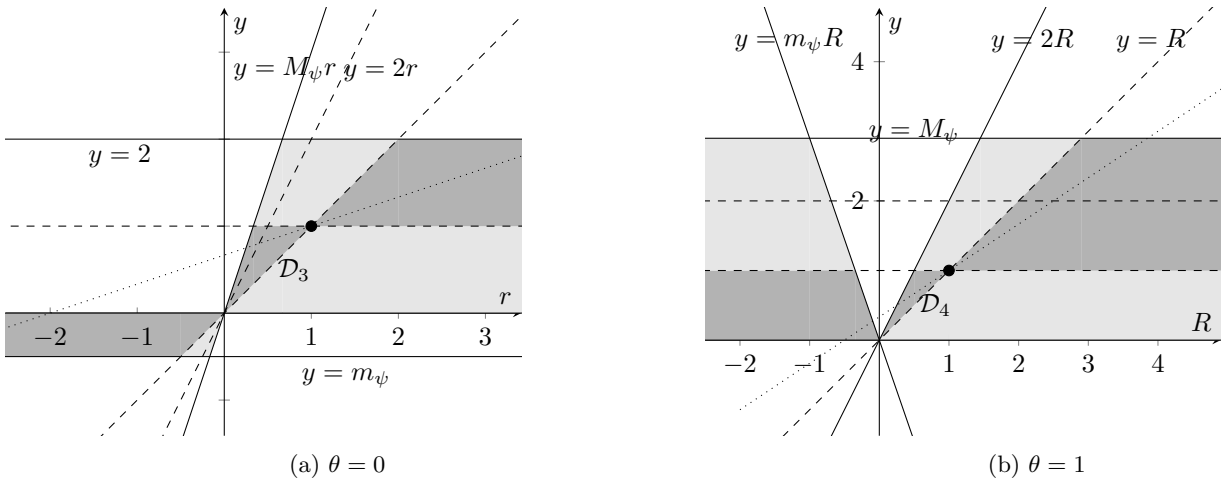

\subsection{Applicability of Spekreijse's limiter region.}\label{sec:Applicability of the extended Spekreijse region}
In this subsection we will discuss why Spekreijse’s theory may not always apply to flux form schemes. It is important to note that Spekreijse's theory is correct as a flux splitting framework for hyperbolic PDE's, and if you were to directly use Spekreijse's flux splitting method as defined in \cite{S_1987} on the advection equation, you would be using an advective form method, where velocity is located at the cell centres, (see \cref{fig: advective form stencil}, \cref{eq:advective form}). In this circumstance, not only will the extended limiter of Spekreijse provably work, but for the advection equation, the monotonicity criteria could be generalised further as the uniform boundedness condition requiring $\psi(R)-\psi(S)/S<2M$ (eq(2.1)\cite{S_1987}) can be relaxed.

However, using an advective form scheme can lead to the loss of mass conservation for non-uniform flow. We instead investigate whether the positive coefficient argument applies to flux form methods in the multidimensional case where there is an additional dependence on an edge-defined velocity field (see \cref{fig:flux form stencil}), as local mass conservation is essential for many applications. To use the positive coefficient lemma of Spekreijse, it is natural to write the flux form scheme in the following positive coefficient form
\begin{align}
\begin{split}
	&\frac{\partial u_{i,j}}{\partial t} + \frac{ (u^{R}_{i}c_{i+1/2}^+ - u^{R}_{i-1}c_{i-1/2}^{+})}{(u^{R}_{i} - u^{R}_{i-1})}\frac{u^{R}_{i} - u^{R}_{i-1}}{u_i - u_{i-1}} (u_i - u_{i-1}) - \frac{ (u^{L}_{i+1}c_{i+1/2}^{-} - u^{L}_{i}c_{i-1/2}^{-} ) }{(u^{L}_{i+1} - u^{L}_{i})}\frac{u^{L}_{i+1} - u^{L}_{i}}{ u_{i+1} - u_{i}} (u_{i} - u_{i+1})+ \text{y-direction}= 0. \label{eq: Spekreijse advective semidiscrete}
\end{split}
\end{align}
In order for this to be of positive coefficient type, we require that the leading terms (as defined in \cref{Def:Positive Coefficient}) are positive. This is traditionally done \cite{S_1987} by the use of limiting to establish $(u^{R}_{i} - u^{R}_{i-1})/(u_i - u_{i-1}) \geq 0$, then this is subsequently used in a mean value theorem on a numerical flux function (typically arising from the flux splitting framework \cite{S_1987} on a flux $f^{+}$ with respect to $u$ to establish $(f^+(u^{R}_{i}) - f^{+}(u^{R}_{i-1}))/(u^{R}_{i} - u^{R}_{i-1}) \geq 0$. However, we are not in the same flux splitting framework as the numerical flux function here has additional dependence on the velocity field and the mean value theorem won't necessarily apply to a term of the form  $(u^{R}_{i}c_{i+1/2}^+ - u^{R}_{i-1}c_{i-1/2}^{+})/(u^{R}_{i} - u^{R}_{i-1})$. It is here one could say that the Spekreijse region is inappropriate for incompressible flow, by choosing an incompressible flow with $(u^{R}_{i},c_{i+1/2}^+,u^{R}_{i-1}, c_{i-1/2}^{+} )= (1,0,0.5,1)$ giving a negative leading coefficient. However this argument is only an indication that the theory doesn't hold, it is not sufficient, one would have to rule out the existence of a transform in which the scheme can be put in a different positive coefficient representation. Instead, we will use a counter example described in \cref{sec:Setup: monotonicity tests}, and rely on a numerical demonstration in \cref{table: minimum values}, \cref{sec:results:positivity}, to show that this region is not monotone for the flux form scheme. Analytically, we also prove some additional necessary conditions are required of the limiter function to remain positivity preserving for incompressible flux form advection (in \Cref{thm:necessary}), limiters in the Spekreijse limiter region can violate these strictly necessary conditions.

This gives us the first contribution of this paper: the extended region of Spekreijse does not give schemes that have a discrete maximum principle or even positivity preserving when applied to the multidimensional flux form incompressible advection equation. 
\begin{remark}
Despite this there are circumstances in which you can still use Spekreijse's limiter region. A sufficient condition for the leading coefficients in \cref{eq: Spekreijse advective semidiscrete} to be positive is that the $x$-component of velocity $v_{1}$, is independent of $x$, so that $v_1 = v_1(y,t)$ and the $y$-component of velocity $v_2$, is independent of $y$ so that $v_2= v_2(x,t)$. This implies that the flow is directionally constant $c_{i+1/2} = c_{i-1/2} = c_{i}$, $c_{j+1/2} = c_{j-1/2} = c_{j}$ property, and we can show that the leading terms will be positive directly
\begin{align}
	\frac{ (u^{R}_{i}c_{i+1/2}^+ - u^{R}_{i-1}c_{i-1/2}^{+})}{(u^{R}_{i} - u^{R}_{i-1})} = c_i^{+} \geq 0.\label{eq: direct computation}
\end{align} 
This can also be inferred from the mean value theorem, and more generally we will call a flow directional mean value theorem satisfying (MVTS), when one has sufficient conditions on the flow to directionally apply a mean value theorem and attain a positive coefficient representation of the numerical scheme. A flow without this property we will call mean value theorem violating (MVTV). 
\end{remark}

\section{Other theoretical properties of the scheme.}\label{Other theoretical properties of the scheme}
\subsection{Temporal discretisation}\label{sec: temporal discretisation}
The strong stability literature (e.g. \cite{bolley1978conservation,spijker1983contractivity,kraaijevanger1991contractivity,shu1988efficient,ferracina2005extension}) allows boundedness and contractive behaviour in arbitrary seminorms of the semi-discrete/(forward Euler scheme) scheme to be to be directly translated into a wide variety of SSP methods, allowing the order of the scheme to be increased. The local maximum principle is not preserved in exactly the same way, we actually ensure that each substage satisfies a local maximum principle with respect to the previous substage, this is apparent from the Shu Osher representation, see \cref{sec:hidden maximum principle}, for a more detailed explanation.
We can recall the forward Euler numerical flow map $u^{n+1} = \text{Forward Euler} ( u^n, c^n), $ defined as
\begin{align}
u_{i,j}^{n+1}  =  u_{i,j}^n &- [F_{i+0.5}(u^{R}_{i},u^{L}_{i+1},c_{i+0.5}^n) - F_{i-0.5}(u^{R}_{i-1},u^{L}_{i},c_{i-0.5}^n)] - [F_{j+0.5}(u^{U}_{j},u^{D}_{j+1},c_{j+0.5}^n) - F_{j-0.5}(u^{U}_{j-1},u^{D}_{j},c_{j-0.5}^n)]. \label{eq:FE}
\end{align}
Where the time step, $\Delta t$, and mesh spacing, $\Delta x$, have been absorbed into the face defined velocity $c$. The SSP33 Runge Kutta scheme, will be used in the numerical experiments and can be implemented in the following  memory efficient (2 register) Shu Osher representation
\begin{align}
	k^1 &= \text{Forward Euler}(u^n, c^{n}),\\
	k^2 &= 3 / 4  \cdot u^n + 1 / 4 \cdot \text{Forward Euler}(k^1, c^{n+1}),\\
	u^{n+1}&= 1 / 3  \cdot u^{n} + 2 / 3 \cdot \text{Forward Euler}(k^2, c^{n+1/2}).\label{eq:ssp33}
\end{align}
This scheme is a third order, three stage Runge Kutta method with radius of monotonicity of 1. This method preserves convex semi-norms (such as $||\cdot||_{\infty}$) under the same timestep restriction forward Euler does. 
In \cite{koren1993robust} it is remarked that for small enough Courant numbers the semi-discrete scheme in this paper is monotone when discretised with RK4. The strong stability literature has since shown that the RK4 method has zero radius of monotonicity, in \cref{sec:RK4} we remark on this particular choice of timestepping method.
\subsection{Accuracy}
\label{accuracy}
The numerical method viewed as a finite volume method approximating the integral form of the equation is second order or less in dimensions two or greater, and third order or less when in one dimension. The numerical method as viewed as a finite difference method approximating the point valued equation is second order, but can be third order for uniform flow. 
See \cref{sec:accuracy of un-limited} for a more detailed discussion. Spekreijse, showed that sufficient conditions on the limiter function $\psi(R)$ for second order accuracy are $\psi(1) = 1$, $\psi \in C^{2}$ near 1, \cite{S_1987}. Sufficient and necessary conditions on the limiter function $\psi(r)$ were shown by Hua-mo \cite{hua1992possible} to be $\psi(1) = 1$, and $\psi$ is Lipschitz continuous. More specific conditions on accuracy of finite difference accuracy of uniform flow in the neighbourhood of extrema can be found in \cite{hua1992possible}. See \cref{sec:accuracy of limited} for a more detailed discussion. It is also often mistaken that TVD schemes are first order at noncritical extrema because of a slight oversight in the truncation analysis in Osher \cite{osher1984high}. This is addressed by Hua-mo for $\theta =1$ \cite{hua1992possible}, but often attributed to \cite{zijlema1998higher}, who make similar mistakes in the truncation error analysis as Osher, for $\theta = 0$. See \cref{sec:accuracy of limited} for a more detailed discussion on this.
\subsection{Linear Invariance}
Linear invariance implies constants are preserved under incompressible flow and essential in atmospheric advection \cite{lin1996multidimensional} \cite{leonard1996conservative}. In \cite{hundsdorfer1995positive} it is incorrectly stated that the scheme is not linear invariant whether in state or flux interpolated form, the state interpolated scheme in this paper is linear invariant for incompressible flow after one takes away a discrete incompressibility condition, see \cref{sec:Linear invariance}.

\subsection{Symmetric limiters: Old and New}\label{sec:Symmetric limiters}

\begin{figure}
\centering
	\begin{subfigure}{0.45\textwidth}
		\begin{tikzpicture}[scale=0.9]
		\pgfdeclarelayer{pre main}
			\pgfsetlayers{pre main,main}
			\begin{axis}[
				axis lines = middle,
				axis equal,
				domain  = -4:4,
				xlabel  = {$r$},
				ylabel  = {$y$},
				xmin    = -3.5,
				xmax    = 3.9,
				ymin    = -1.1,
				ymax    = 2.1,
				samples = 200,
				mark    = none,
				]
				\pgfonlayer{pre main}
				\endpgfonlayer
				\addplot [name path=Vanalbada, thick,color=red,width=1]  {( (x*x + x)/(x*x+1))};
				\addplot [name path=Ospre, thick,color=blue]  {(1.5* (x*x + x)/(x*x+x+1))};
				\addplot [name path=Eno2, thick,color=green,domain=-1:1]  { (x) };
				\addplot [name path=Eno2, thick,color=green,domain=-4:-1]  { -1 };
				\addplot [name path=Eno2, thick,color=green,domain=1:4]  { 1 };
				\addplot [name path=twor, dashed]  {(2*x)};
				\addplot [name path=zero, thick]  {(0)};
				\addplot [name path=ninf, thin]  {(-3)};
				\addplot [name path=inf, thin]  {(4)};
				\addplot [name path=two, dashed]  {(2)};
				\addplot[color=black!30] fill between[of=zero and twor,
				soft clip={domain=0:1}];
				\addplot[color=black!30] fill between[of=zero and two,
				soft clip={domain=1:3.9}];
				\node at (0.7,0.5) {$\mathcal{D}_{1}$};
				\node at (1.3,2.5) {$y = 2r$};
			\end{axis}
		\end{tikzpicture}
		\caption{Sweby admissible region, as well as the ENO2 scheme in green, Ospre in blue, and van Albada in red}
		\label{fig:Eno2 green}
	\end{subfigure}
	\quad
	\begin{subfigure}{0.4\textwidth}
		\begin{tikzpicture}[scale=0.9]
			\pgfdeclarelayer{pre main}
			\pgfsetlayers{pre main,main}
			\begin{axis}[
				axis lines = middle,
				axis equal,
				domain  = -2:3,
				xlabel  = {$R$},
				ylabel  = {$y$},
				xmin    = -2,
				xmax    = 2.9,
				ymin    = -1,
				ymax    = 2.9,
				samples = 100,
				mark    = none,
				yticklabels={,,}
				xticklabel=empty,
				ticks=none
				]
				\pgfonlayer{pre main}
				\endpgfonlayer
				\addplot [name path=Ospre, thick,color=blue]  {(1.5* (x*x + x)/(x*x+x+1))};
				\addplot [name path=mMr, dashed]  {(-1.5*x)};
				\addplot [name path=twor, dashed]  {((2)*x)};
				\addplot [name path=twoalphar, dashed]  {((2-0.5)*x)};
				\addplot [name path=zero, thick]  {(0)};
				\addplot [name path=ninf, thin]  {(-4)};
				\addplot [name path=inf, thin]  {(4)};
				\addplot [name path=alpha, thick]  {(-0.5)};
				\addplot [name path=M, thick]  {(1.5)};
				\addplot[color=blue!10] fill between[of=mMr and twoalphar,
				soft clip={domain=0:0.333}];
				\addplot[color=blue!10] fill between[of=alpha and twoalphar,
				soft clip={domain=0.333:1}];
				\addplot[color=blue!10] fill between[of=alpha and M,
				soft clip={domain=1:3}];
				\addplot[color=blue!10] fill between[of=mMr and twoalphar,
				soft clip={domain=-0.333:0}];
				\addplot[color=blue!10] fill between[of=mMr and alpha,
				soft clip={domain=-1:-0.333}];
				\addplot[color=blue!10] fill between[of=M and alpha,
				soft clip={domain=-2:-1}];
				\node at (0.7,0.5) {$\mathcal{D}_{2}$};
				\node at (1.8,2.5) {$y = (2+\alpha)R$};
				\node at (2.4,1.7) {$y= M$};
				\node at (-1.0,2.5) {$y= -M R$};
				\node at (2,-0.7) {$y = \alpha$};
			\end{axis}
		\end{tikzpicture}
		\caption{Spekreijse limiter region for the Ospre scheme.}
		\label{fig:spekOspre}
	\end{subfigure}
	\caption{$\mathcal{D}_{1}$ is the Sweby region, $\mathcal{D}_{2}$, is the Spekreijse region, which has two free parameters $\alpha \in [-\infty,0]$, $M \in (0,\infty)$. We can see that the limiters, van Albada (red) \cref{limiter: Van Albada}, Ospre (blue) \cref{limiter: Ospre}, and ENO2 (green) \cref{limiter: Eno2} are not contained in the Sweby region $\mathcal{D}_{1}$ (or even the new incompressible flow limiter region $\mathcal{D}_{4}$), but are in the Spekreijse region $\mathcal{D}_{2}$ with values for $[M,\alpha]$ given by $[1.5,-0.5]$, $[1/2(1+\sqrt{2}),1/2(1-\sqrt{2})]$, $[1,-1]$ respectively. Only the Ospre scheme and its Spekreijse Region are shown in \cref{fig:spekOspre}. } 
	\label{fig: ospre, van, eno}
\end{figure}
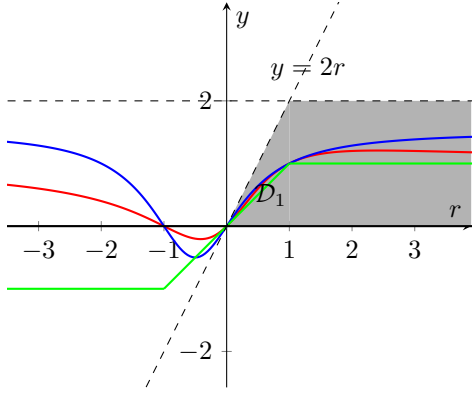
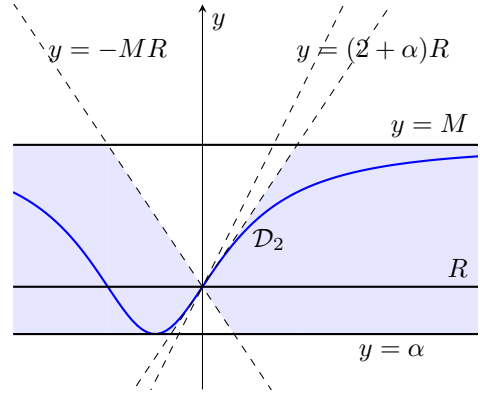

We first introduce some important symmetric flux limiter functions in Speckreijse admissible region, but not in the Sweby admissible region
\begin{align}
	\text{van-Albada$(R)$} &= \frac{R^2+R}{R^2+1}, \label{limiter: Van Albada}\\
	\text{Ospre$(R)$} &= \frac{3}{2}\frac{R^2+R}{R^2+R+1}, \label{limiter: Ospre}\\
	\text{ENO2$(R)$} &=
	\begin{cases}
		R \quad \text{where}\quad |R|\leq 1,\\
		1  \quad \text{where}\quad |R|\geq 1.
	\end{cases}\label{limiter: Eno2}
\end{align}
These are introduced in the respective papers \cite{van1997comparative,waterson1995unified,harten1987eno} and are plotted in \cref{fig:Eno2 green}. 
We also introduce the subscript $P$ limiters, as those same limiters who are ``pushed" into being in the Sweby diagram by the removal of the tail, these limiters are designed keep the symmetry property and fix the monotonicity problems associated with using Spekreijse's region for multidimensional flux form incompressible flow.
\begin{align}
	\text{van-Albada$_P(R)$} &=	
	\begin{cases} 
		\frac{R^2+R}{R^2+1}, \quad R\geq 0 \\
		0, \quad R<0 \\
	\end{cases},\label{limiter: Van Albaba P} \\
	\quad 	\text{Ospre$_P(R)$} &=
	\begin{cases} 
		\frac{3}{2}\frac{R^2+R}{R^2+R+1}, \quad R \geq 0\\
		0, \quad R<0 \label{limiter: OspreP}
	\end{cases},	\\
	\quad 
	\text{ENO$_P(R)$} &= \text{minmod}(R) = \max (0,\min(R,1)).\label{limiter: Eno2P}
\end{align}


\subsection{Non symmetric limiters and symmetry breaking}\label{sec:Non symmetric limiters and symmetry breaking}

\begin{figure}
\centering
	\begin{subfigure}{0.45\textwidth}
		\begin{tikzpicture}[scale=1]
			\pgfdeclarelayer{pre main}
			\pgfsetlayers{pre main,main}
			\begin{axis}[
				axis lines = middle,
				axis equal,
				enlargelimits,
				domain  = -4:4,
				xlabel  = {$r$},
				ylabel  = {$y$},
				xmin    = -2,
				xmax    = 2.9,
				ymin    = -1,
				ymax    = 3.9,
				samples = 100,
				mark    = none,
				yticklabels={,,}
				ticks=none
				]
				\pgfonlayer{pre main}
				\endpgfonlayer
				\addplot[name path = WoofR, thick, domain=-7:-5, color = red] {1/(x+4) };
				\addplot[name path = WoofR, thick, domain=-5:-2, color = red] {1/3*x+2/3};
				\addplot[name path = WoofR, thick, domain=-2:0, color = red] {0};
				\addplot[name path = WoofR, thick, domain=0:5, color = red] {min(3*x,1/3*x+2/3,2)};
				\addplot[name path = analytic, thick, domain=-5:5, color = green] { tanh(x*exp(x))};
				
				\addplot[name path = diffble2, thick, domain=-5:0, color =blue] { tanh(x)*exp(x)};
				\addplot[name path = diffble2, thick, domain=0:0.5, color =blue] { -8*x*x*x+16/3*x*x+x};
				\addplot[name path = diffble2, thick, domain=0.5:3, color =blue] { 1/3*x+2/3};
				\addplot[name path = diffble2, thick, domain=3:5, color =blue] { 1/3*tanh(x-3)+5/3};
				
				\addplot[name path = sou,  thin,domain=-3:5, color = black ,dashed] {x};
				\addplot[name path = cds,  thin,domain=-3:5, color = black ,dashed] { 1};
				\addplot[mark=*] coordinates {(1,1)};
				\addplot [name path=zero, thin]  {(0)};
				\addplot [name path=ninf, thin]  {(-4)};
				\addplot [name path=inf, thin]  {(5)};
				\addplot [name path=mpsi, thick]  {(-0.5)};
				\addplot [name path=twor, thin,dotted]  {(2*x)};
				\addplot [name path=Ms, thin]  {(3*x)};
				\addplot [name path=two, thick]  {(2)};
				\addplot[color=black!10] fill between[of=mpsi and zero,
				soft clip={domain=-4:-0.3333*0.5}];
				\addplot[color=black!10] fill between[of=Ms and zero,
				soft clip={domain=-0.3333*0.5:0}];
				\addplot[color=black!10] fill between[of= zero and Ms,
				soft clip={domain=0:1.3333*0.5}];
				\addplot[color=black!10] fill between[of= zero and two,
				soft clip={domain=1.3333*0.5:4}];
				\addplot[color=black!30] fill between[of=mpsi and zero,
				soft clip={domain=-4:-0.5}];
				\addplot[color=black!30] fill between[of=zero and sou,
				soft clip={domain=-0.5:0.0}];
				\addplot[color=black!30] fill between[of=sou and Ms,
				soft clip={domain=0.0:0.333333}];
				\addplot[color=black!30] fill between[of=sou and cds,
				soft clip={domain=0.333333:1}];
				\addplot[color=black!30] fill between[of=sou and cds,
				soft clip={domain=1:2}];
				\addplot[color=black!30] fill between[of=two and cds,
				soft clip={domain=2:5}];
				\node at (0.8,0.5) {$\mathcal{D}_{4}$};
				\node at (1.4,2.8) {$y = M_{\psi} r$};
				\node at (2.2,3.8) {$y =2 r$};
				\node at (-1.2,1.8) {$y= 2$};
				\node at (1.4,-0.75) {$y = m_{\psi}$};
			\end{axis}
		\end{tikzpicture}
		\caption{The new Differentiable$(r)$ limiter in blue, the analytic limiter $\tanh(r)\exp(r)$ in green, and the Woodfield$(r,M=3,-m=-1)$ limiter in red.}
		\label{fig:continuous limiter}
	\end{subfigure}
	\quad
	\begin{subfigure}{0.45\textwidth}
		\begin{tikzpicture}[scale=1]
			\pgfdeclarelayer{pre main}
			\pgfsetlayers{pre main,main}
			\begin{axis}[
				axis lines = middle,
				axis equal,
				domain  = -3:5,
				xlabel  = {$R$},
				ylabel  = {$y$},
				xmin    = -2.5,
				xmax    = 4.9,
				ymin    = -0.5,
				ymax    = 3.9,
				samples = 300,
				mark    = none,
				]
				\pgfonlayer{pre main}
				\endpgfonlayer
                \addplot[name path = TCDF, thick, domain=-3:0, color = green] { x*(1+x)/(x*x+1) };
                \addplot[name path = TCDF, thick, domain=0:0.5, color = green] { x*x*x-2*x*x +2*x };
                \addplot[name path = TCDF, thick, domain=0.5:2.0, color = green] { 0.75*x+0.25};
                \addplot[name path = TCDF, thick, domain=2.0:5, color = green] { (2*x*x-2*x-2.25) /(x*x -x -1)};
				\addplot[name path = superbeeR, thick, domain=-3:0, color = red] { min(1,-2.9*x) };
				\addplot[name path = superbeeR, thick, domain=0:5, color = red] { max(min(1,2*x),min(x,2.9)) };
				\addplot[name path = woodfield, thick, domain=-3:-1, color = blue] { 0 };
				\addplot[name path = woodfield, thick, domain=-1:-0.5, color = blue] { 0 };
				\addplot[name path = woodfield,  thick,domain=-0.5:-1/(3*2.9+2), color = blue] { 2/3*x + 1/3 };
				\addplot[name path = woodfield, thick, domain=-1/(3*2.9+2):0, color = blue] { -2.9*x  };
				\addplot[name path = woodfield,  thick,domain=0:0.25), color =blue] { 2*x };
				\addplot[name path = woodfield,  thick,domain=0.25:0.5*(3*2.9-1), color =blue] { 2/3*x + 1/3 };
				\addplot[name path = woodfield,  thick,domain=0.5*(3*2.9-1):5, color =blue] { 2.9};
				\addplot [name path=mMr, dashed]  {(-2.9*x)};
				\addplot [name path=twor, dashed]  {(2*x)};
                \addplot [name path=1p5r, dotted,color=blue]  {(1.5*x)};
				\addplot [name path=thirdorder, dotted, domain={-2:5}]  {(2/3*x+1/3)};
				\addplot [name path=r, dashed]  {(x)};
				\addplot [name path=one, dashed]  {(1)};
				\addplot [name path=zero, thin]  {(0)};
				\addplot [name path=ninf, thin]  {(-2)};
				\addplot [name path=inf, thin]  {(6)};
				\addplot [name path=two, dashed]  {(2)};
				\addplot [name path=M, dotted]  {(2.9)};
				\addplot[color=black!30] fill between[of=zero and one,
				soft clip={domain=-3:-.5}];
				\addplot[color=black!30] fill between[of=mMr and zero,
				soft clip={domain=-0.5:0}];
				\addplot[color=black!10] fill between[of=one and M,
				soft clip={domain=-3:-1}];
				\addplot[color=black!10] fill between[of=zero and mMr,
				soft clip={domain=-1:0}];
				\addplot[color=black!30] fill between[of=zero and one,
				soft clip={domain=-1:-0.3448275862}];
				\addplot[color=black!30] fill between[of=zero and mMr,
				soft clip={domain=-0.3448275862:0}];
				\addplot[color=black!30] fill between[of=r and twor,
				soft clip={domain=0:0.5}];
				\addplot[color=black!30] fill between[of=r and one,
				soft clip={domain=0.5:1}];
				\addplot[color=black!10] fill between[of=one and twor,
				soft clip={domain=0.5:1.45}];
				\addplot[color=black!10] fill between[of=M and zero,
				soft clip={domain=1.45:5}];
				\addplot[color=black!10] fill between[of=zero and r,
				soft clip={domain=0:1}];
				\addplot[color=black!10] fill between[of=zero and one,
				soft clip={domain=1:5}];
				\addplot[color=black!30] fill between[of=one and r,
				soft clip={domain=1:2.9}];
				\addplot[color=black!30] fill between[of=one and M,
				soft clip={domain=2.9:5}];
				\addplot[mark=*] coordinates {(1,1)};
				\node at (0.8,0.5) {$\mathcal{D}_{3}$};
				\node at (2.2,4.3) {$y = 2R$};
				\node at (3.9,4.3) {$y = R$};
				\node at (0.5,3.0) {$y = M_{\psi}$};
				\node at (-1.3,4.3) {$y = m_{\psi}R$};
			\end{axis}
		\end{tikzpicture}
	\caption{ The Woodfield$(R,M=2.9,m=2.9)$ limiter in blue, the SuperbeeR$(R,M=2.9,-m=-2.9)$ limiter in red, and the Ultimate TCDF$(R)$ limiter \cite{zhang2015review} in green.}
	\label{fig:woodfieldR}
    \end{subfigure}
\caption{New limiters in their respective $\theta =0,1$ regions.}
\label{region: Limiter regions with examples}
\end{figure}
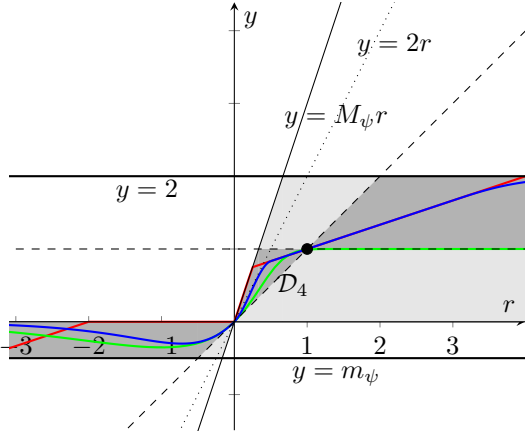
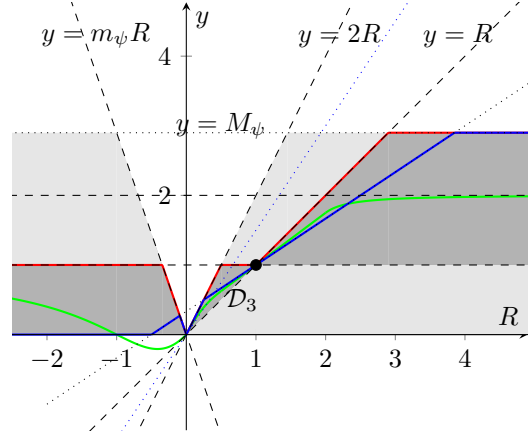

We also introduce the Koren \cite{koren1993robust} limiter,
\begin{align}
	\text{Koren$(R)$} &= \max(0,\min(2,2R,2R/3 + 1/3)). \label{limiter: Koren}
\end{align}
It is not symmetric, but is an accurate limiter consisting of restricting the third order upwind region to the Sweby region. We introduce a new limiter with free parameters $M, m$ defined by 
\begin{align}
	\text{Woodfield}(R,M,m)=
	\begin{cases} 
		&\left\{ R\leq -1 :0  \right\} \\
		&\left\{-1<R<-\frac{1}{2}: 0, \right\}\\
		&\left\{-\frac{1}{2}<R<-\frac{1}{-3m-2}:\ \frac{2}{3}R+\frac{1}{3}\right\}\\
		&\left\{\frac{1}{-3m-2}<R<0:\ m_{\psi}R\right\}\\
		&\left\{0<R<\frac{1}{4}:\ 2R\right\}\\
		&\left\{\frac{1}{4}<R<\frac{3M-1}{2}:\ \frac{2}{3}R+\frac{1}{3}\right\}\\
		&\left\{\frac{3M-1}{2}<R:\ M\right\}
	\end{cases} \label{limiter: Woodfield}
\end{align}
consisting of restricting the third order upwind region to the new limiter region. We also introduce the SuperbeeR limiter, 
\begin{align}
	\text{SuperbeeR}(R,M,m)=
	\begin{cases} 
		& \max (0, \max( \min(2 R, 1), \min (R, M))) \quad  \text{where}\quad R \geq 0\\
		& \min (m R ,1 ) \quad  \text{where}\quad  R<0
	\end{cases} \label{limiter: superbeeR}
\end{align}
which serves as an extension of traditional Superbee limiter to the new limiter region. This is to test our limiter region but could have application to front tracking. These limiters are plotted in \cref{fig:woodfieldR}. We define a similar extension to the Koren limiter by restricting the third order upwind region to the $\theta =1$ maximum principle limiter region for some $M,m$ creating the Woodfield$(r,M,m)$ limiter, plotted in \cref{fig:continuous limiter}. 


In the review paper \cite{zhang2015review} it is argued on practical experience that using regions in which $R<0$ makes little difference, a similar sentiment is also found in the review \cite{kemm2011comparative}. For example, in \cite{zhang2015review} a limiter called ultimate-TCDF (Third Order Continuously Differentiable) is introduced, it is defined as
\begin{align}
\text{UTCDF}(R)=
\begin{cases} 
&
\frac{R(R+1)}{(R^2+1)},  \quad  \text{where}\quad  R < 0, \\
&R^3-2R^2 +2R , \quad  0\leq R \leq 1/2 ,\\
	& 0.75 R+0.25, \quad (1/2\leq R<2),\\
		& \frac{(2R^2-2R-9/4)}{(R^2-R-1)},  \quad R\geq 2.
	\end{cases} \label{limiter: UTCDF}
\end{align}
First note the UTCDF limiter is not differentiable at $R=0$, and more importantly, has no theoretical guarantees. It is not contained within any Spekreijse limiter region\footnote{This can be verified since the tail $R<0$ is identical to the van-Albada limiter so that $\alpha=-0.5$. However, one can see that at $R=0$, the gradient is $2>2-0.5$, and for $R\in (0,0.2928)$ this limiter is not contained within a Spekreijse region.}, it also leaves the newly defined $\theta=1$ limiter region, as seen in green in \cref{region: Limiter regions with examples}. 
We conjecture that leaving the Spekreijse region may cause problems in practice (fail to be monotonic), in the context of directionally constant flow, as well as flux splitting frameworks. We also conjecture that leaving the newly defined limiter region $\theta=1$ will cause problems in practice for incompressible flux form advection, in light of \cref{thm:necessary}.

Instead one can design limiter functions using the Spekreijse limiter region when in a flux-splitting framework (or when a mean value theorem applies) or using the new $\theta=0$, $\theta = 1$ limiter regions in the case of incompressible flux form advection. For example we define the following two limiter functions
\begin{align}
\text{UTCDF}_P(R)=
\begin{cases} 
&
0\quad \text{where}\quad R < 0\\
&R^3-2R^2 +2R , \quad  0\leq R \leq 1/2 ,\\
& 0.75 R+0.25, \quad (1/2\leq R<2),\\
& \frac{(2R^2-2R-9/4)}{(R^2-R-1)},  \quad R\geq 2.
	\end{cases} \label{limiter: TCDF_P}
\end{align}
\begin{align}
\text{UTCDF}_S(R)=
\begin{cases} 
&
\frac{R(R+1)}{(R^2+1)},  \quad  \text{where}\quad  R < -1, \\
&0\quad \text{where}\quad  -1\leq R < 0\\
&R^3-2R^2 +2R , \quad  0\leq R \leq 1/2 ,\\
& 0.75 R+0.25, \quad (1/2\leq R<2),\\
& \frac{(2R^2-2R-9/4)}{(R^2-R-1)},  \quad R\geq 2.
	\end{cases} \label{limiter: TCDF_S}
\end{align}
The $\text{UTCDF}_P(R)$ limiter is the Ultimate-TCDF restricted to the Sweby region. The $\text{UTCDF}_S(R)$ limiter is formed by restricting the Ultimate-TCDF to the new $\theta=1$ limiter region. For differentiable limiters suitable for flux splitting, one can design limiters in the Spekreijse region e.g. \cref{limiter: Ospre}. 
When $\theta = 1$, there are no globally differentiable limiter functions contained within the second order limiter region, that remain suitable for incompressible flow without a mean value theorem assumption of the flow. We introduce the first globally differentiable limiter function contained entirely within the second order region, suitable for incompressible flow by using our new $\theta = 0$ region, it touches the third order region for accuracy. 
\begin{align}
	\text{Differentiable}(r)=
	\begin{cases} 
		&\tanh(r) \exp(r),  \quad  \text{where}\quad  r \leq 0, \\
		&-8 r^3 + 16/3 r^2 + r , \quad  0<r \leq1/2 \\
		& 1/3 r+2/3, \quad 1/2<r \leq 3\\
		& 1/3 \tanh(r-3) + 5/3 ,  \quad r>3
	\end{cases} \label{limiter: differentiable new}
\end{align}
plotted in \cref{fig:continuous limiter}. This limiter is a placeholder that could be improved on/replaced by piecewise polynomials with the same properties.

\section{Numerical Demonstrations: Test setup and results}\label{sec:Numerical Demonstrations}

We have two distinct types of flow (Mean value theorem satisfying(MVTS), and Mean value theorem violating(MVTV)), three different limiter regions (Spekreijse, Sweby, and our new limiter region(s)). We introduce two monotonicity tests:
\begin{itemize}
    \item Solid body rotation of the LeVeque initial conditions in which a directional mean value theorem holds and all limiter functions should remain monotone \cite{S_1987}.
    \item A sinusoidal deformational flow of the LeVeque initial conditions which is directional mean value theorem violating, used to verify that the Spekreijse limiter region is not monotone for general incompressible flow, but the Sweby region and the new limiter regions are. We also introduce another deformational vector field \cref{test:pathalogical} with the same aim. 
\end{itemize}

\subsection{Setup: monotonicity tests}\label{sec:Setup: monotonicity tests}
The numerical domain is $\Omega = (0,1]\times(0,1]$ discretised by  $128 \times 128$ cells, with periodic boundary conditions. We run the scheme with 4096 timesteps until $T_{max}=1$ resulting in a maximum cell defined Courant number around $0.2$.  We locate the stream function $\Psi(x,y)$ at the cell vertices and use a discrete form of the curl operator to create a divergence free vector field (to machine precision) with velocities located at the midpoints of faces, with only normal components. We use the solid body rotation and the sinusoidal deformation defined from the stream functions 
\begin{align}
    \Psi = - \pi \left( (x - 1/2)^2 + (y - 1/2)^2\right) ,\label{test:solid body rotation} \\
    \Psi = \frac{1}{2}\sin(\pi x)\sin(\pi y) \cos(2 \pi t/T_{max}) \label{test:sin deformation},\\
    \Psi =  1/16 \sin(32 \pi x)  \sin(32 \pi y)  \cos(2\pi t/T_{max})\label{test:pathalogical}.
\end{align}
The solid body rotation flow \cref{test:solid body rotation} is directional mean value theorem satisfying (with the directionally constant property), the sinusoidal deformation \cref{test:sin deformation} does not have a directional mean value property. \Cref{test:pathalogical}
is similar to \cref{test:sin deformation}, it does not have a mean value theorem, it is used to create larger magnitude monotonicity failures.
For the initial condition of the tracer, we use the LeVeque initial conditions \cite{Leveque_test_ic}, defined below
\begin{align}
	u_{0} & = 
	\begin{cases}
		1 &  \text{where}\quad \sqrt{ (x-0.5)^2+(y-0.75)^2} \leq 0.15, \quad \text{and}\quad (x\leq 0.475), \\
		1 &  \text{where}\quad \sqrt{ (x-0.5)^2+(y-0.75)^2} \leq 0.15, \quad \text{and}\quad (x>0.525), \\
		1 &  \text{where}\quad \sqrt{ (x-0.5)^2+(y-0.75)^2} \leq 0.15, \quad (y\geq 0.85), \quad \text{and}\quad (0.475<x\leq 0.525), \\
		(1-\frac{R_{cone}}{0.15}) &  \text{where}\quad (R_{cone} = \sqrt{(x-0.5)^2+(y-0.25)^2}\leq 0.15), \\
		\frac{1}{2}(1+\cos(\pi  \frac{R_{cos}}{0.15}) &  \text{where}\quad (R_{cos} = \sqrt{(x-0.25)^2+(y-0.5)^2}\leq 0.15),\\
		0&\text{else}.
	\end{cases}
	\label{ic: Le veque}
\end{align}


\subsection{Numerical Results}

\Cref{sec:snapshots} will show snapshots of the solutions. Positivity will be demonstrated in \cref{sec:results:positivity} and error norms are presented in \cref{sec:accuracy}.

\subsubsection{Snapshots}\label{sec:snapshots}

\begin{figure}[H]
\centering
	\begin{subfigure}[t]{0.19\textwidth}		\includegraphics[width=\textwidth,trim={5mm 7mm 5mm 0mm},clip]{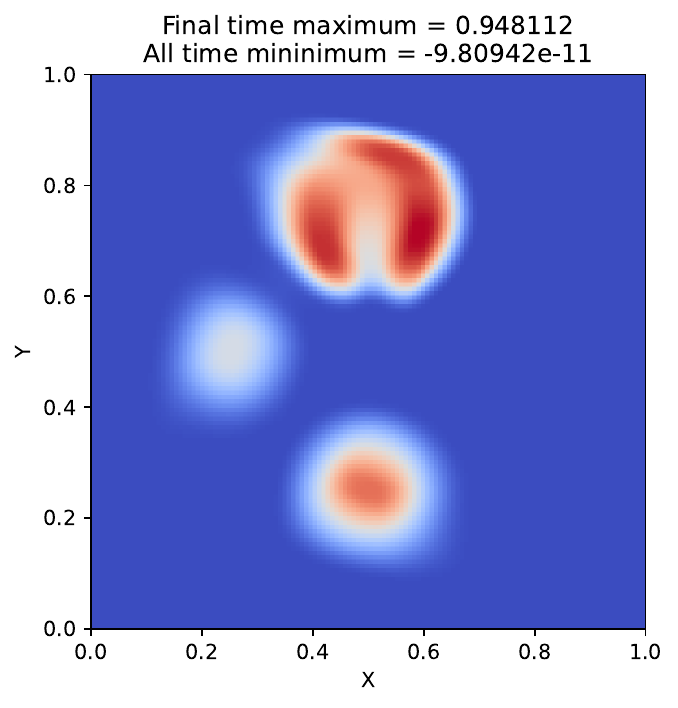}\hspace*{-0.5em} 
		\caption{RK4:Koren$(R)$}
	\end{subfigure}
	\begin{subfigure}[t]{0.19\textwidth}\includegraphics[width=\textwidth,trim={5mm 7mm 5mm 0mm},clip]{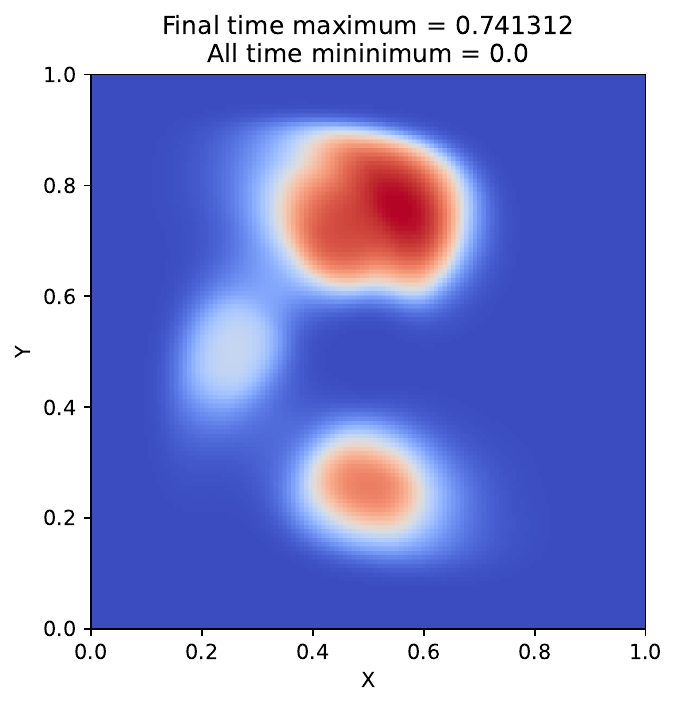}
    \caption{SSP33:ENO2P$(R)$}
	\end{subfigure}
	\begin{subfigure}[t]{0.19\textwidth}
    \includegraphics[width=\textwidth,trim={5mm 7mm 5mm 0mm},clip]{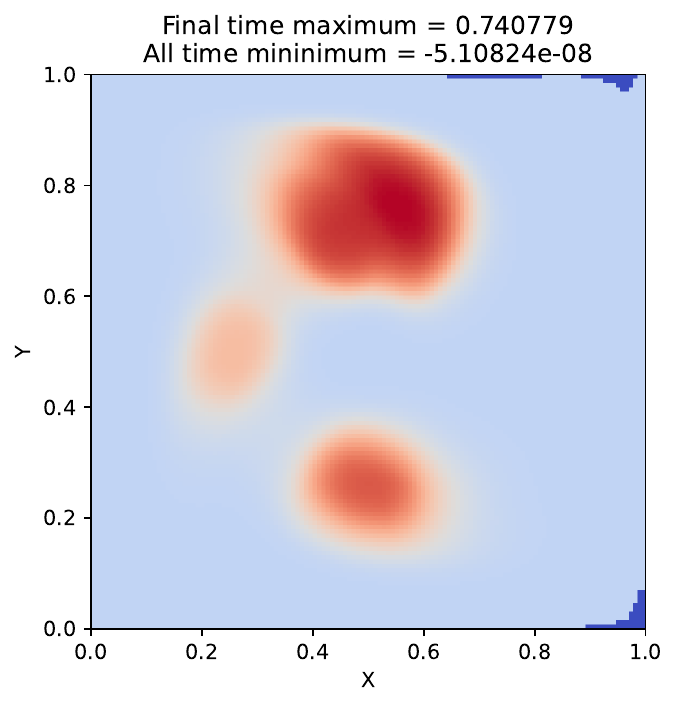}\hspace*{-0.5em}
    \caption{SSP33:ENO2$(R)$}
	\end{subfigure}
	\begin{subfigure}[t]{0.19\textwidth}
    \includegraphics[width=\textwidth,trim={5mm 7mm 5mm 0mm},clip]{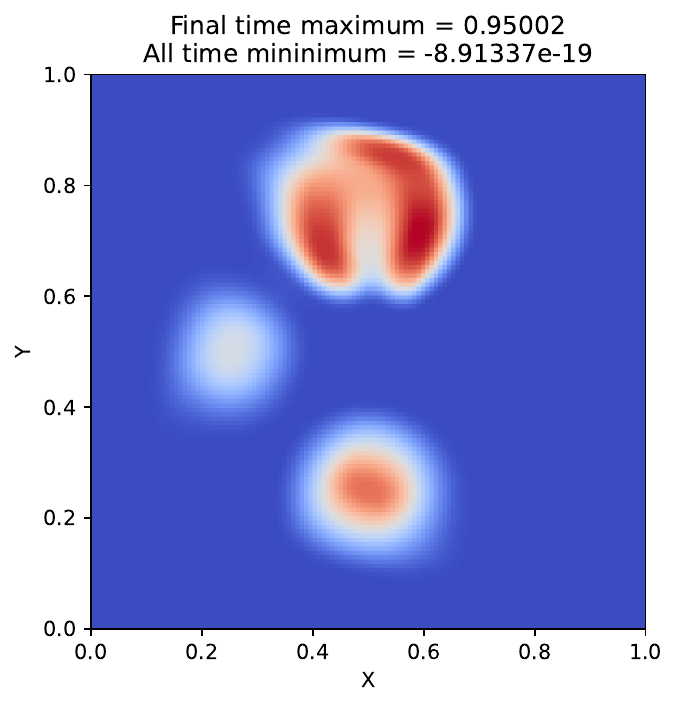}\hspace*{-0.5em}
	\caption{SSP33:Woodfield$(R,2,-2)$}
	\end{subfigure}
	\begin{subfigure}[t]{0.19\textwidth}
	\includegraphics[width=\textwidth,trim={5mm 7mm 5mm 0mm},clip]{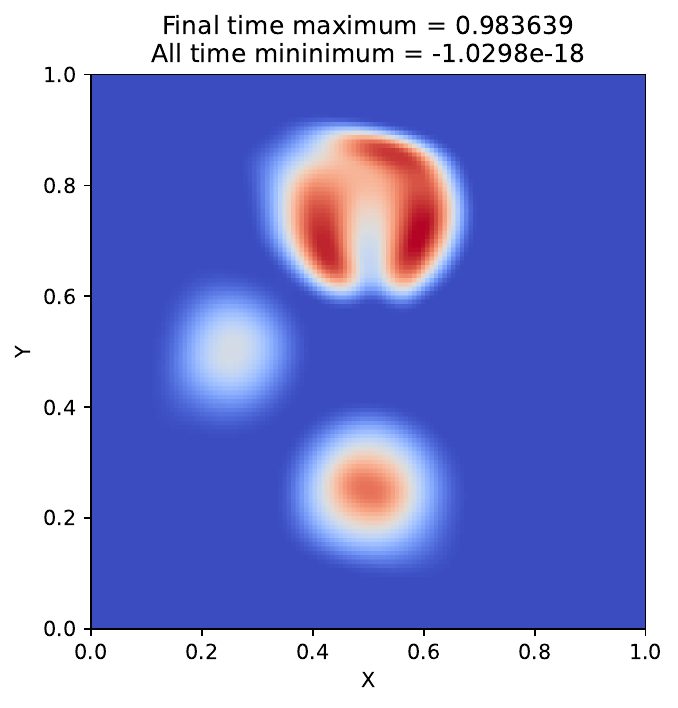}
	\caption{SSP33:Woodfield$(R,4,0)$}
	\end{subfigure} \\
 
 	\begin{subfigure}[t]{0.19\textwidth}
    \includegraphics[width=\textwidth,trim={5mm 7mm 5mm 0mm},clip]{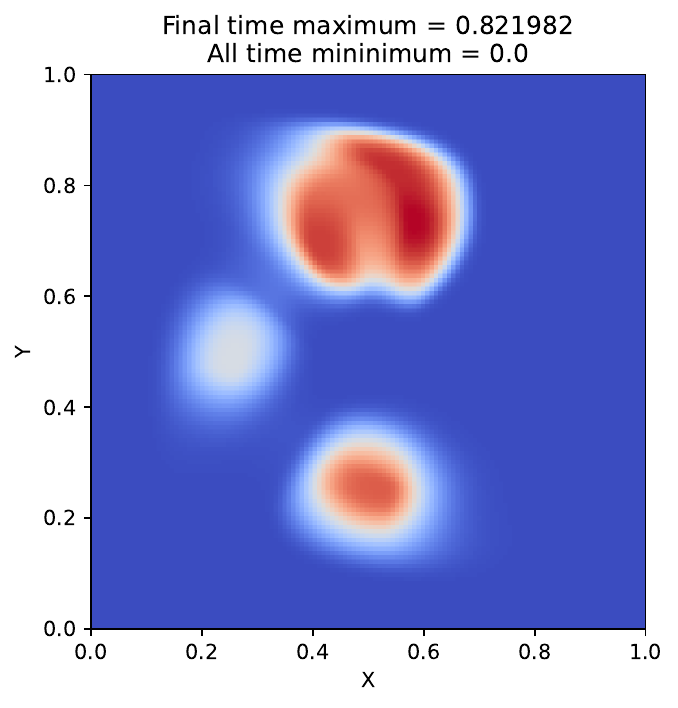}
	\caption{SSP33:Ospre$_P(R)$}
	\end{subfigure}	\begin{subfigure}[t]{0.19\textwidth}
	\includegraphics[width=\textwidth,trim={5mm 7mm 5mm 0mm},clip]{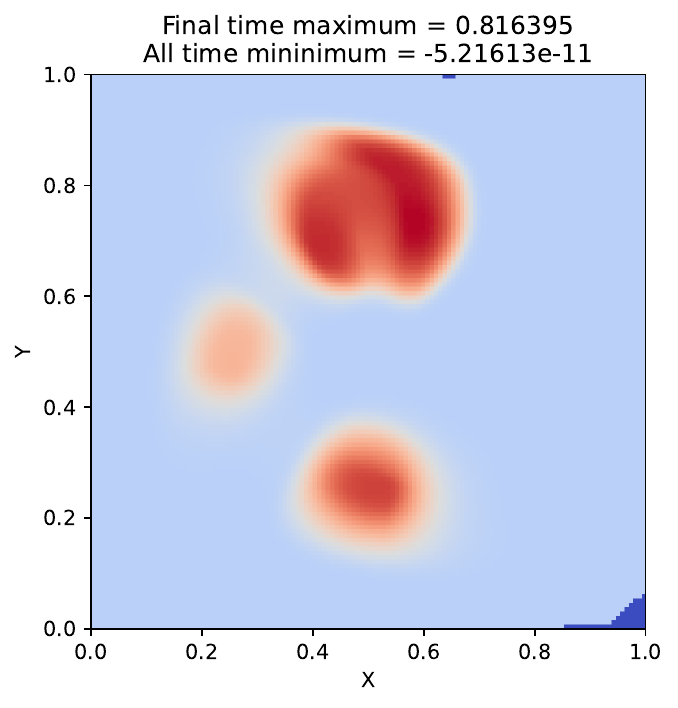}
	\caption{SSP33:Ospre$(R)$}
	\end{subfigure}	\begin{subfigure}[t]{0.19\textwidth}
	\includegraphics[width=\textwidth,trim={5mm 7mm 5mm 0mm},clip]{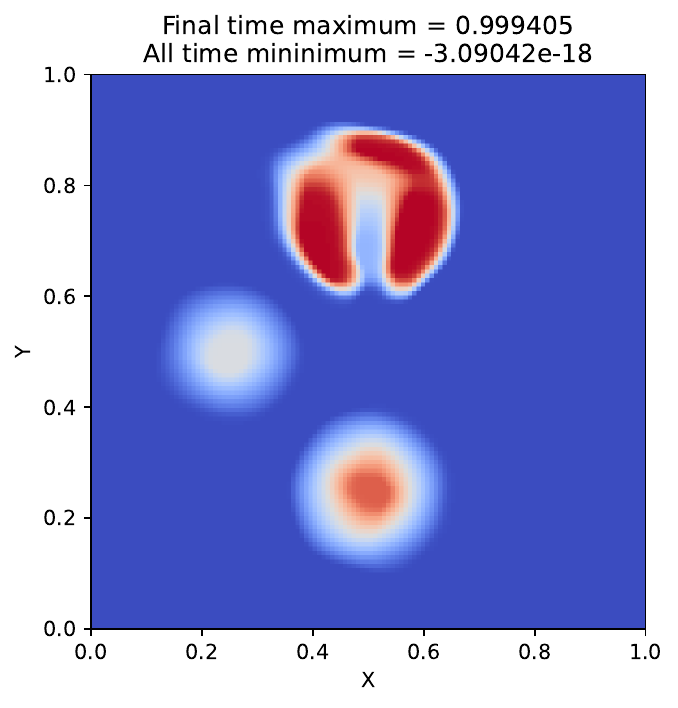}
	\caption{SSP33:superbee$(R)$}
	\end{subfigure}	\begin{subfigure}[t]{0.19\textwidth}
	\includegraphics[width=\textwidth,trim={5mm 7mm 5mm 0mm},clip]{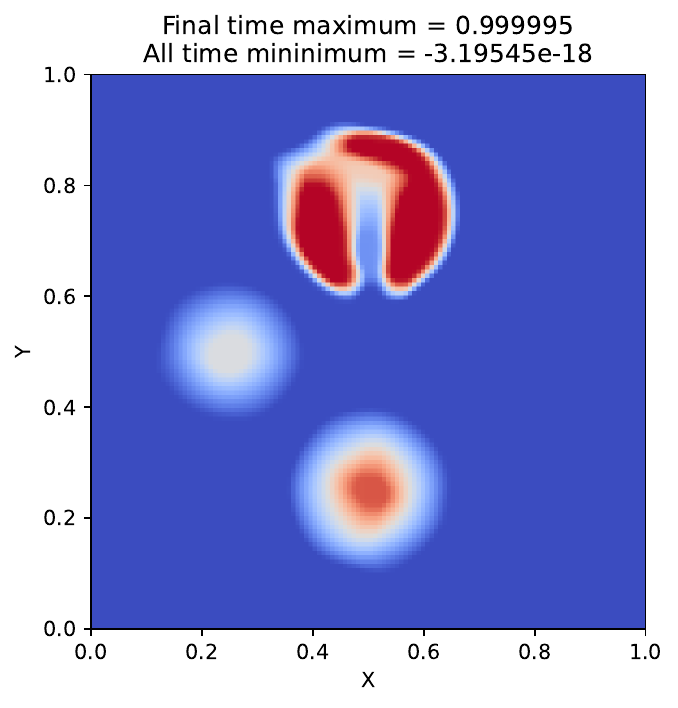}
	\caption{SSP33:superbeeR$(R)$}
	\end{subfigure}	
    \begin{subfigure}[t]{0.19\textwidth}
	\includegraphics[width=\textwidth,trim={5mm 7mm 5mm 0mm},clip]{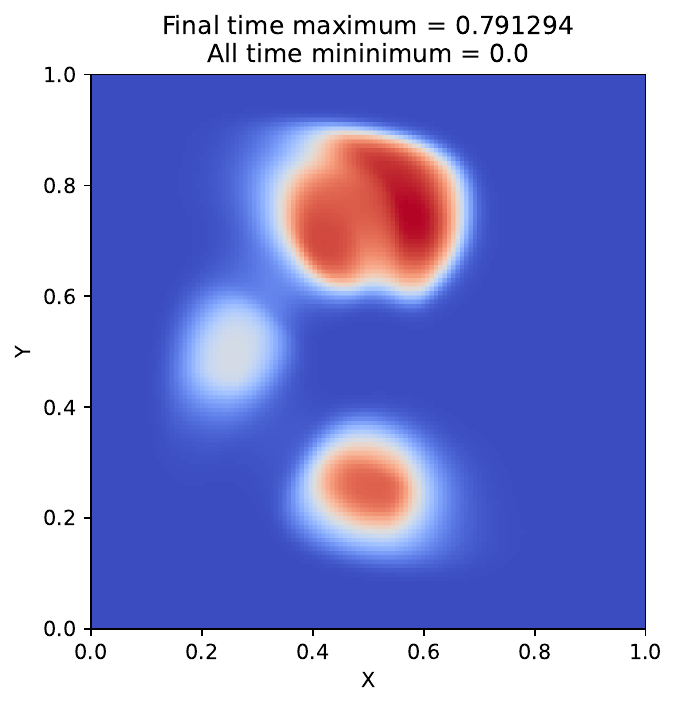}
	\caption{SSP33:vanalbada$_P(R)$}
	\end{subfigure} \\
 
    \begin{subfigure}[t]{0.19\textwidth}
    \includegraphics[width=\textwidth,trim={5mm 7mm 5mm 0mm},clip]{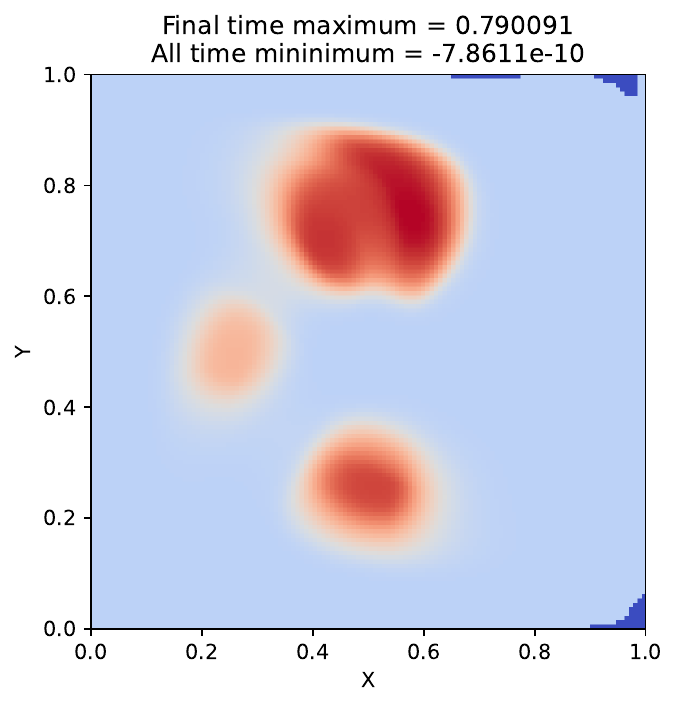}
	\caption{SSP33:VanAlbada$(R)$}
	\end{subfigure} 
    \begin{subfigure}[t]{0.19\textwidth}
    \includegraphics[width=\textwidth,trim={5mm 7mm 5mm 0mm},clip]{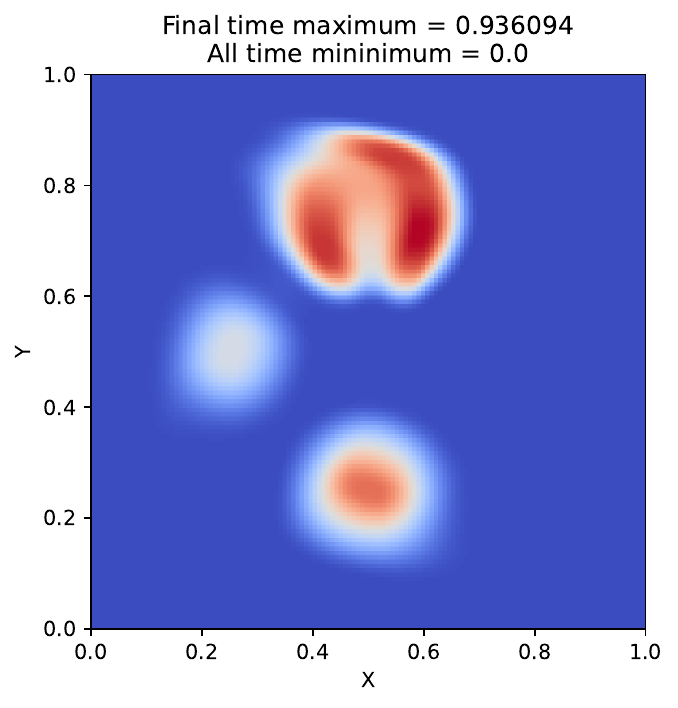}
	\caption{SSP33:Differentiable$(R)$}
	\end{subfigure} 
 \begin{subfigure}[t]{0.19\textwidth}
 \includegraphics[width=\textwidth,trim={5mm 7mm 5mm 0mm},clip]{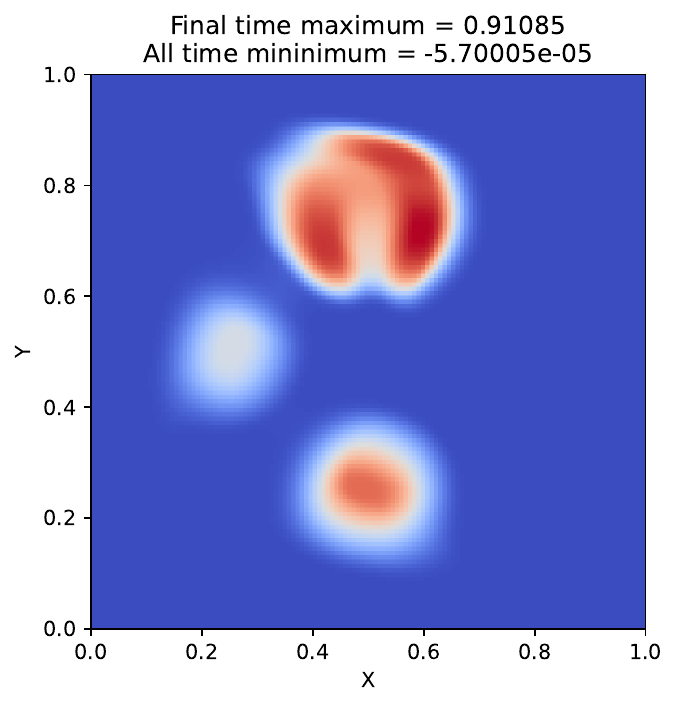}
	\caption{SSP33:UTCDF$(R)$}
	\end{subfigure} 
 \begin{subfigure}[t]{0.19\textwidth}
 \includegraphics[width=\textwidth,trim={5mm 7mm 5mm 0mm},clip]{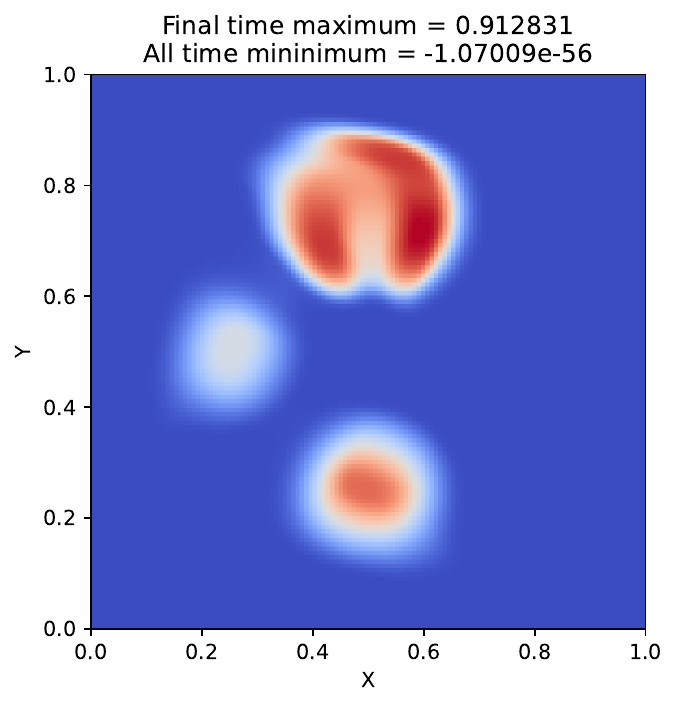}
	\caption{SSP33:UTCDF$_{P}(R)$}
	\end{subfigure} 
 \begin{subfigure}[t]{0.19\textwidth}
 \includegraphics[width=\textwidth,trim={5mm 7mm 5mm 0mm},clip]{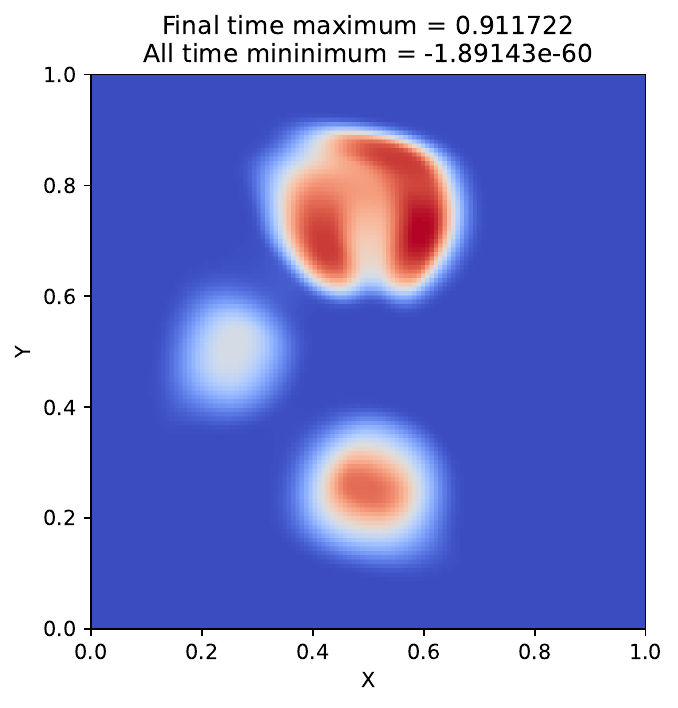}
	\caption{SSP33:UTCDF$_{S}(R)$}
	\end{subfigure} 
 \caption{ Final state of the time reversing sinusoidal flow of LeVeque initial conditions, any negative values below negative $-1\times 10^{-14}$ will be plotted as if $-0.5$ and will shift the entire colour range, so that the colour scheme highlights negatives should they appear. Final time maximum and all time minimum values are displayed on the top of each figure.  }\label{fig: sinusoidal deformational reversing on zal schemes}
\end{figure}

\begin{figure}[h]
\centering
	\begin{subfigure}[t]{0.19\textwidth}
	\includegraphics[width=\textwidth,trim={5mm 7mm 5mm 0mm},clip]{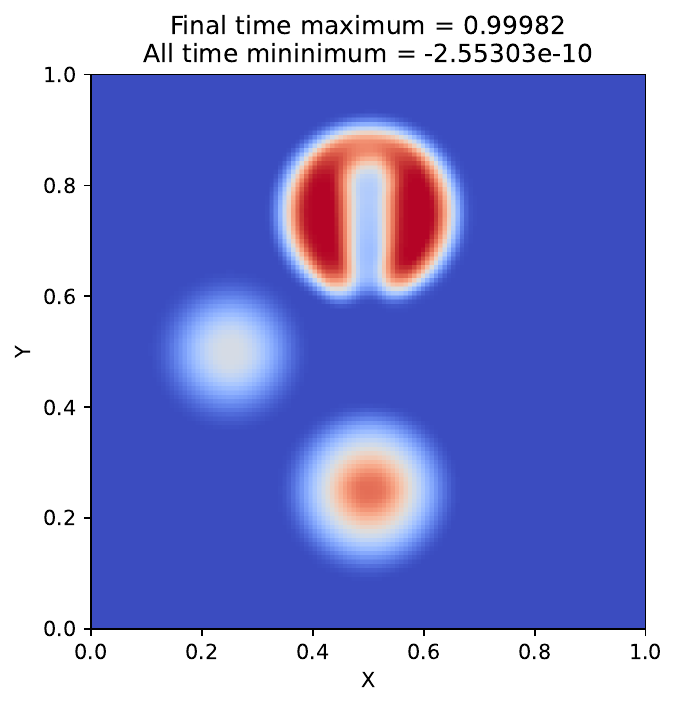}\hspace*{-0.5em} 
		\caption{RK4:Koren$(R)$}
	\end{subfigure}
	\begin{subfigure}[t]{0.19\textwidth}\includegraphics[width=\textwidth,trim={5mm 7mm 5mm 0mm},clip]{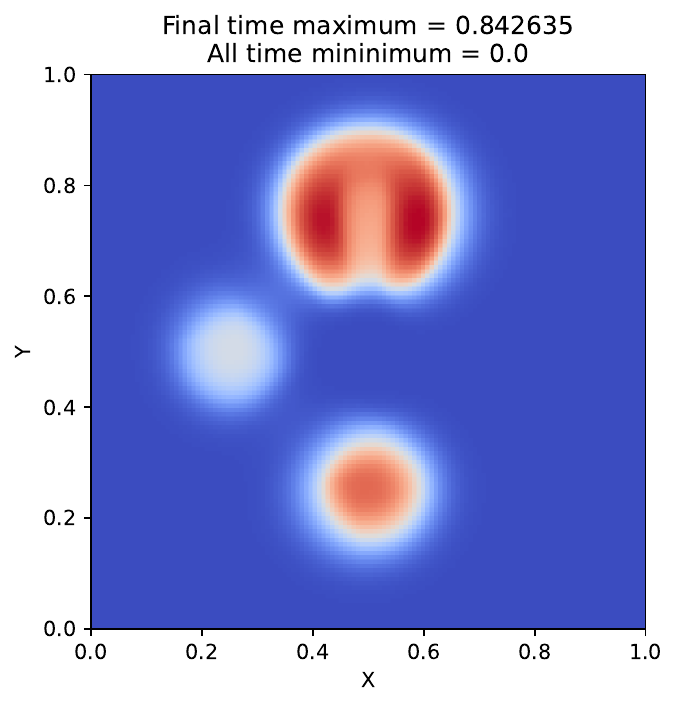}
    \caption{SSP33:ENO2P$(R)$}
	\end{subfigure}
	\begin{subfigure}[t]{0.19\textwidth}
    \includegraphics[width=\textwidth,trim={5mm 7mm 5mm 0mm},clip]{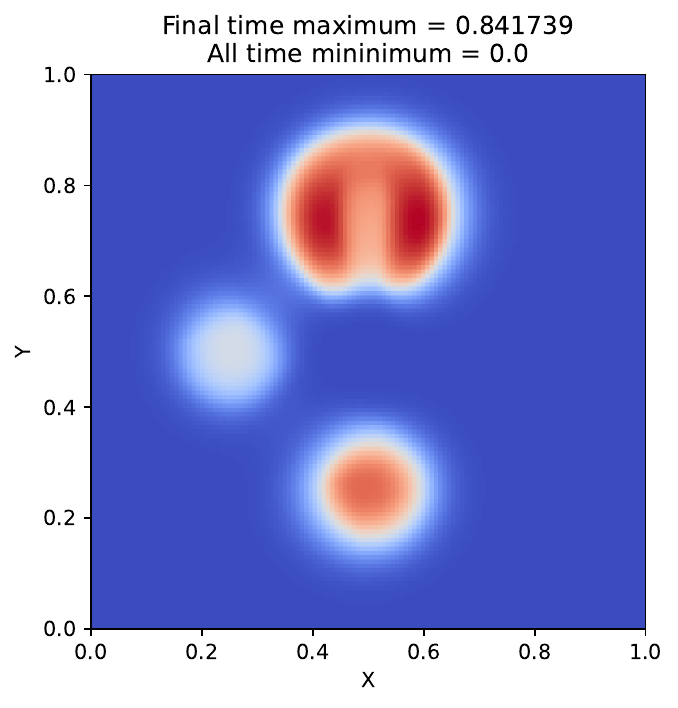}\hspace*{-0.5em}
    \caption{SSP33:ENO2$(R)$}
	\end{subfigure}
	\begin{subfigure}[t]{0.19\textwidth}
    \includegraphics[width=\textwidth,trim={5mm 7mm 5mm 0mm},clip]{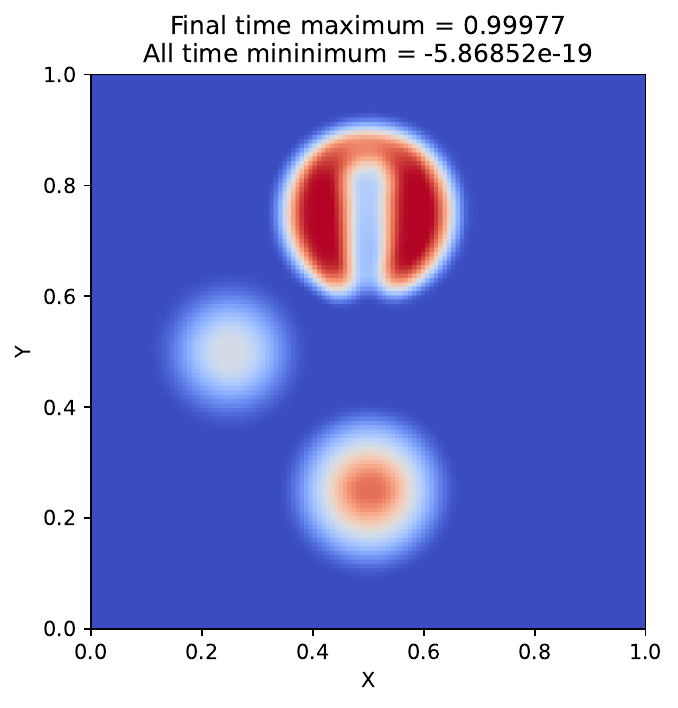}\hspace*{-0.5em}
	\caption{SSP33:Woodfield$(R,2,-2)$}
	\end{subfigure}
\begin{subfigure}[t]{0.19\textwidth}
\includegraphics[width=\textwidth,trim={5mm 7mm 5mm 0mm},clip]{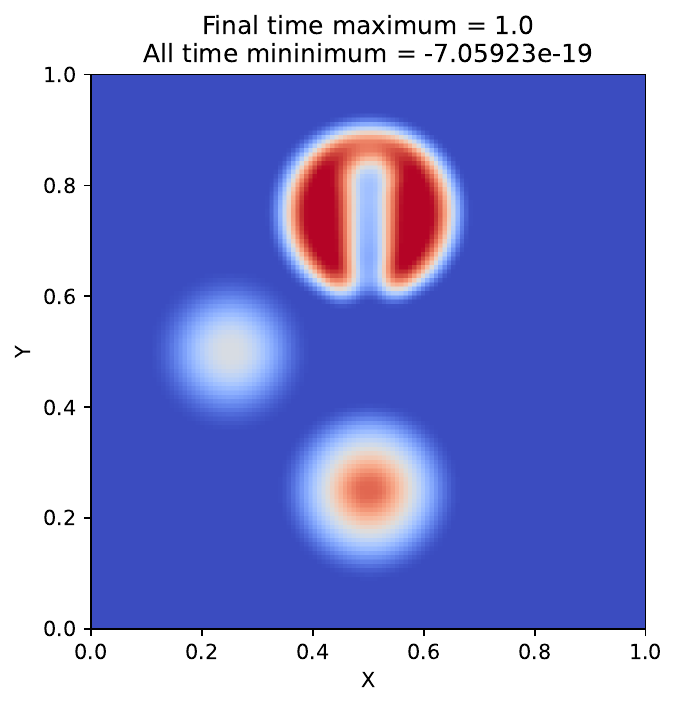}
\caption{SSP33:Woodfield$(R,4,0)$}
\end{subfigure} \\
 \begin{subfigure}[t]{0.19\textwidth}
\includegraphics[width=\textwidth,trim={5mm 7mm 5mm 0mm},clip]{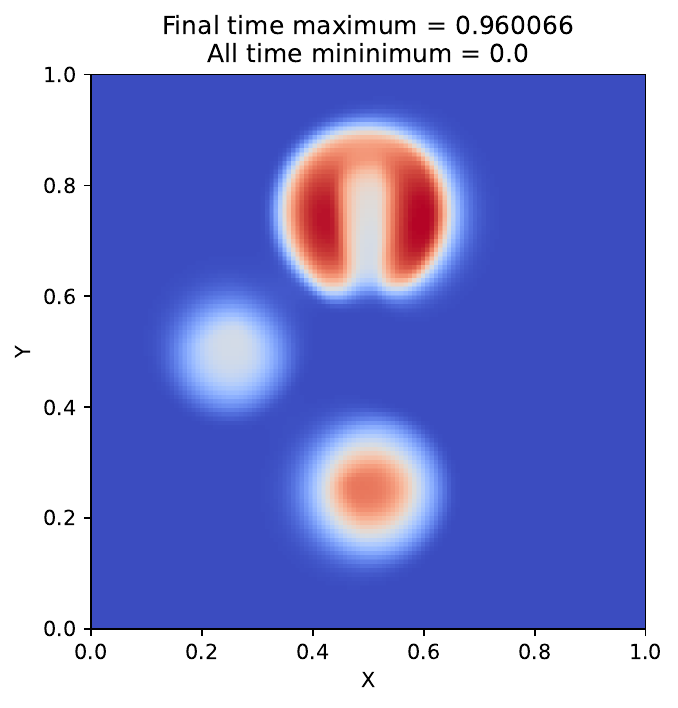}
\caption{SSP33:Ospre$_P(R)$}
\end{subfigure}	\begin{subfigure}[t]{0.19\textwidth}
\includegraphics[width=\textwidth,trim={5mm 7mm 5mm 0mm},clip]{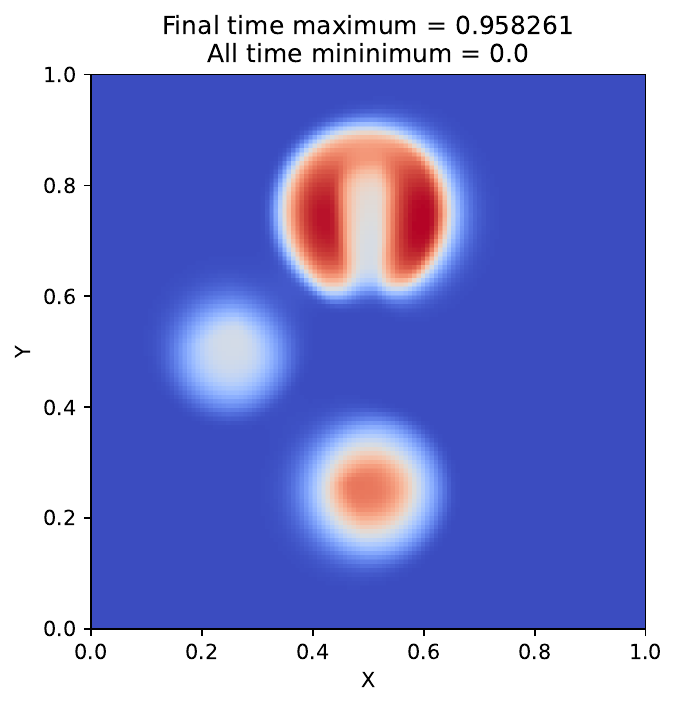}
\caption{SSP33:Ospre$(R)$}
\end{subfigure}	\begin{subfigure}[t]{0.19\textwidth}
\includegraphics[width=\textwidth,trim={5mm 7mm 5mm 0mm},clip]{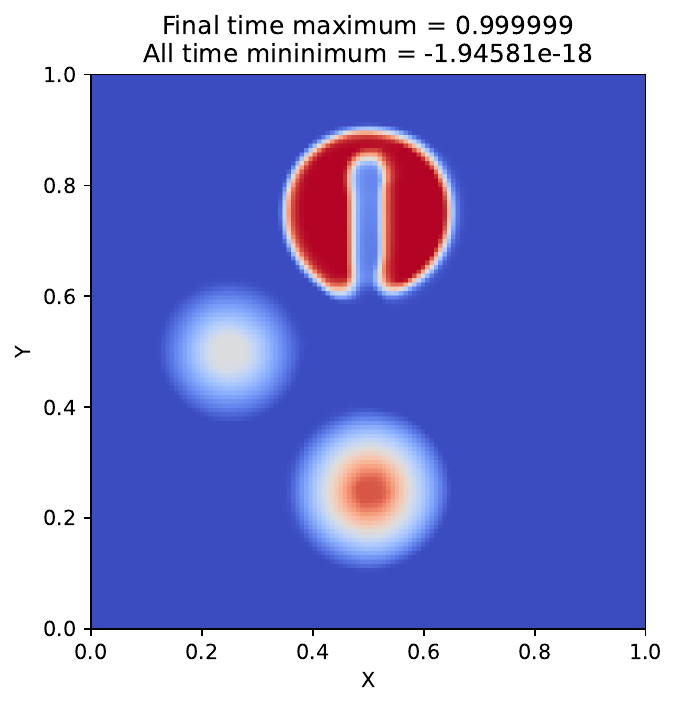}
\caption{SSP33:Superbee$(R)$}
\end{subfigure}	\begin{subfigure}[t]{0.19\textwidth}
\includegraphics[width=\textwidth,trim={5mm 7mm 5mm 0mm},clip]{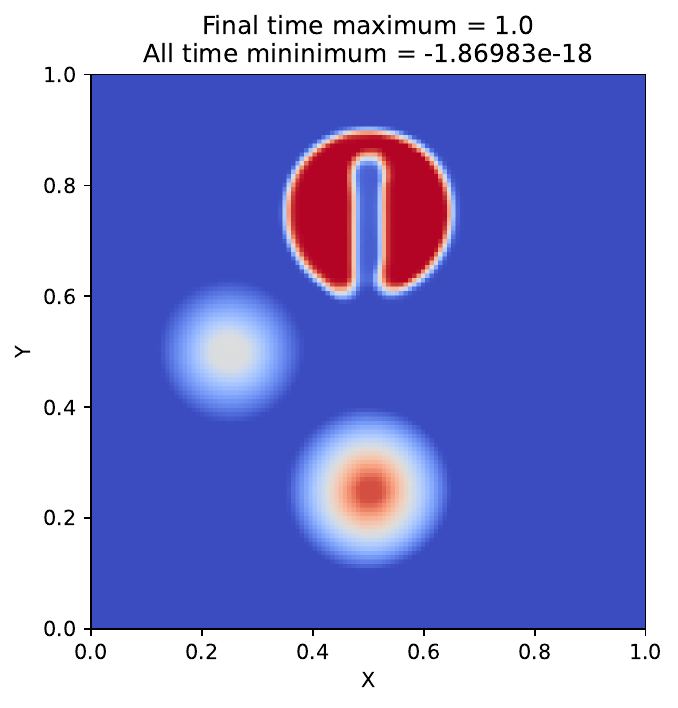}
\caption{SSP33:SuperbeeR$(R)$}
\end{subfigure}	
\begin{subfigure}[t]{0.19\textwidth}
\includegraphics[width=\textwidth,trim={5mm 7mm 5mm 0mm},clip]{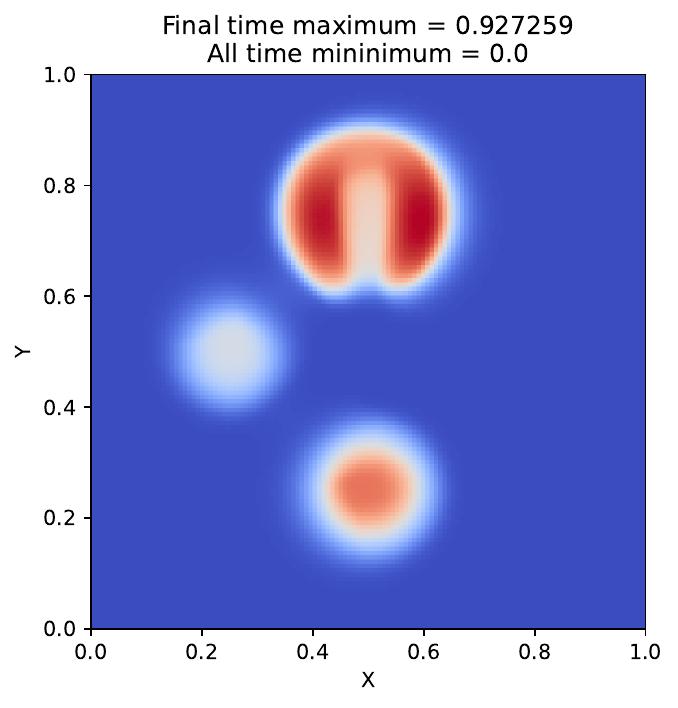}
 \caption{SSP33:VanAlbada$_P(R)$}
\end{subfigure} \\
\begin{subfigure}[t]{0.19\textwidth}
\includegraphics[width=\textwidth,trim={5mm 7mm 5mm 0mm},clip]{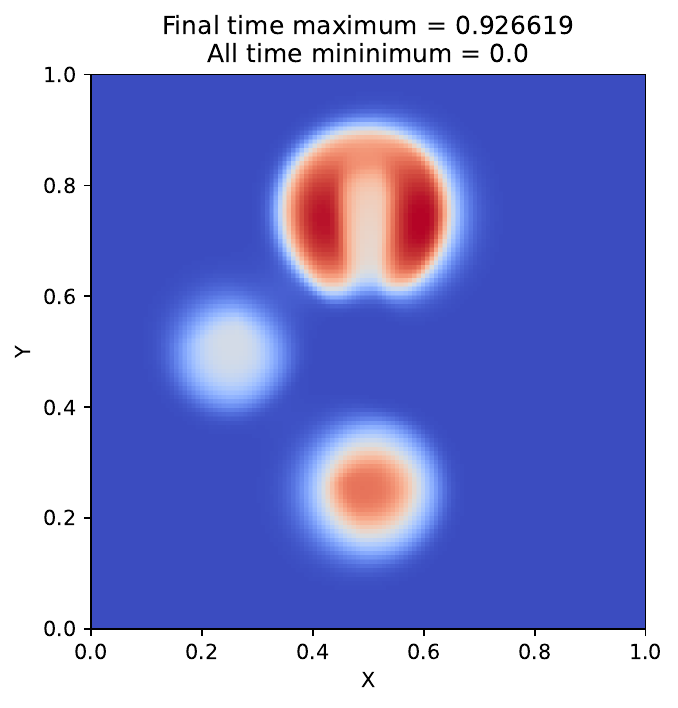}
\caption{SSP33:VanAlbada$(R)$}
\end{subfigure} 
\begin{subfigure}[t]{0.19\textwidth}
\includegraphics[width=\textwidth,trim={5mm 7mm 5mm 0mm},clip]{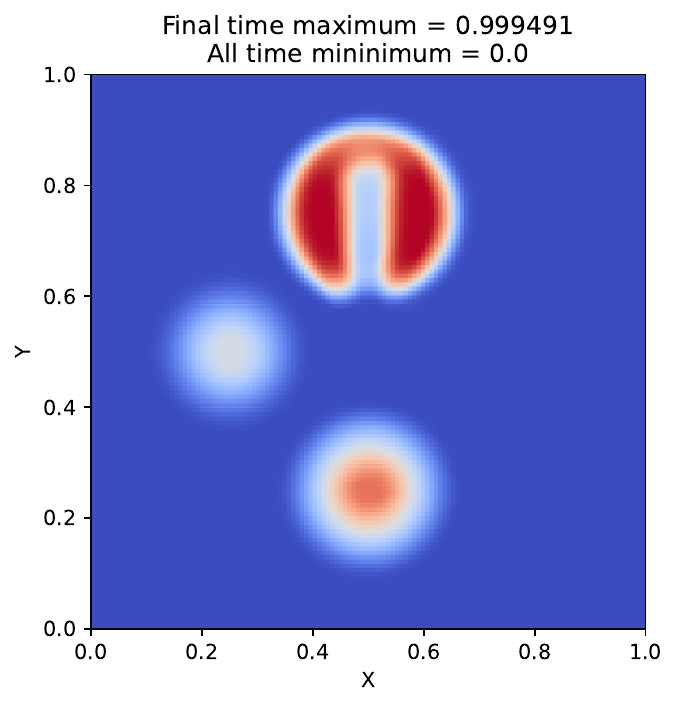}
\caption{SSP33:Differentiable$(R)$}\end{subfigure}
 \begin{subfigure}[t]{0.19\textwidth}
\includegraphics[width=\textwidth,trim={5mm 7mm 5mm 0mm},clip]{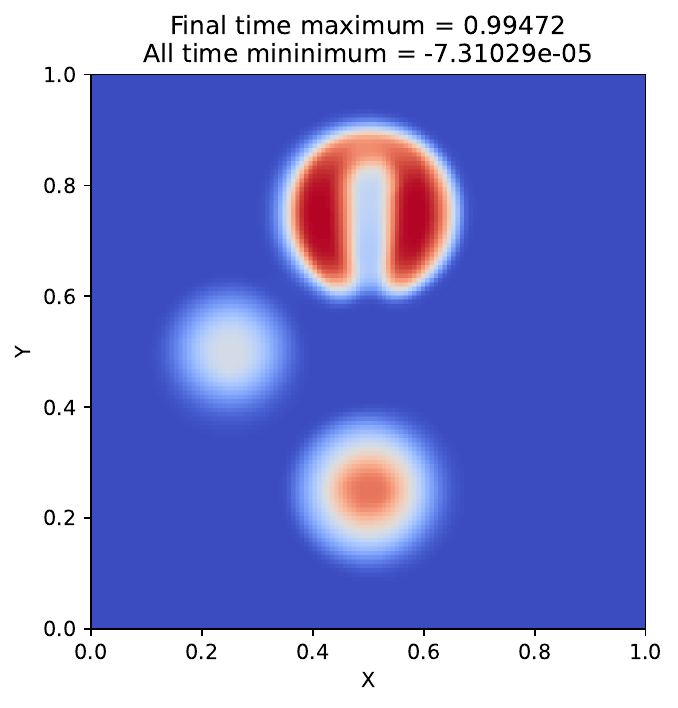}
\caption{SSP33:UTCDF$(R)$}
\end{subfigure} 
\begin{subfigure}[t]{0.19\textwidth}
\includegraphics[width=\textwidth,trim={5mm 7mm 5mm 0mm},clip]{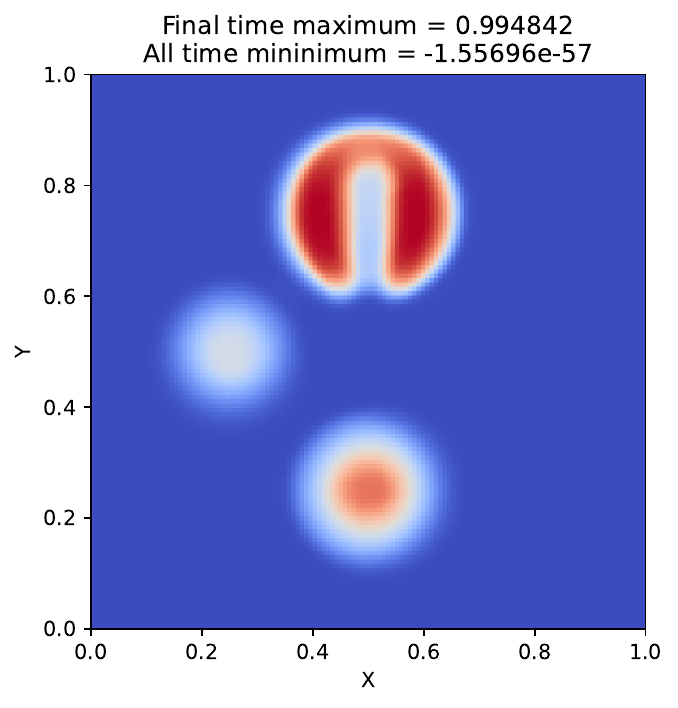}
\caption{SSP33:UTCDF$_{P}(R)$}
\end{subfigure} 
\begin{subfigure}[t]{0.19\textwidth}
\includegraphics[width=\textwidth,trim={5mm 7mm 5mm 0mm},clip]{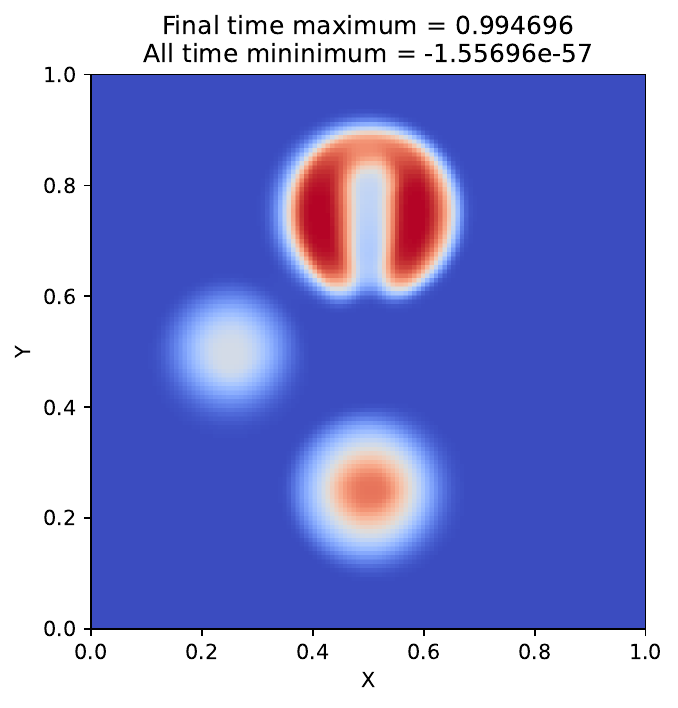}
\caption{SSP33:UTCDF$_{S}(R)$}
\end{subfigure} 
	\caption{ Final state of the solid body rotation of LeVeque initial conditions, any negative values below negative $-1\times 10^{-14}$ will be plotted as if $-0.5$ and will shift the entire colour range, so that the colour scheme highlights significant negatives at the final timestep should they appear. The Final time maximum $\max_{i,j}u^{n}$ and all time minimum values $\max_{i,j,k} u_{i,j}^{k}$are displayed on the top of each figure.  }\label{fig: sbr on zal schemes}
\end{figure}

\begin{figure}[H]
\centering
\includegraphics[scale=0.65]{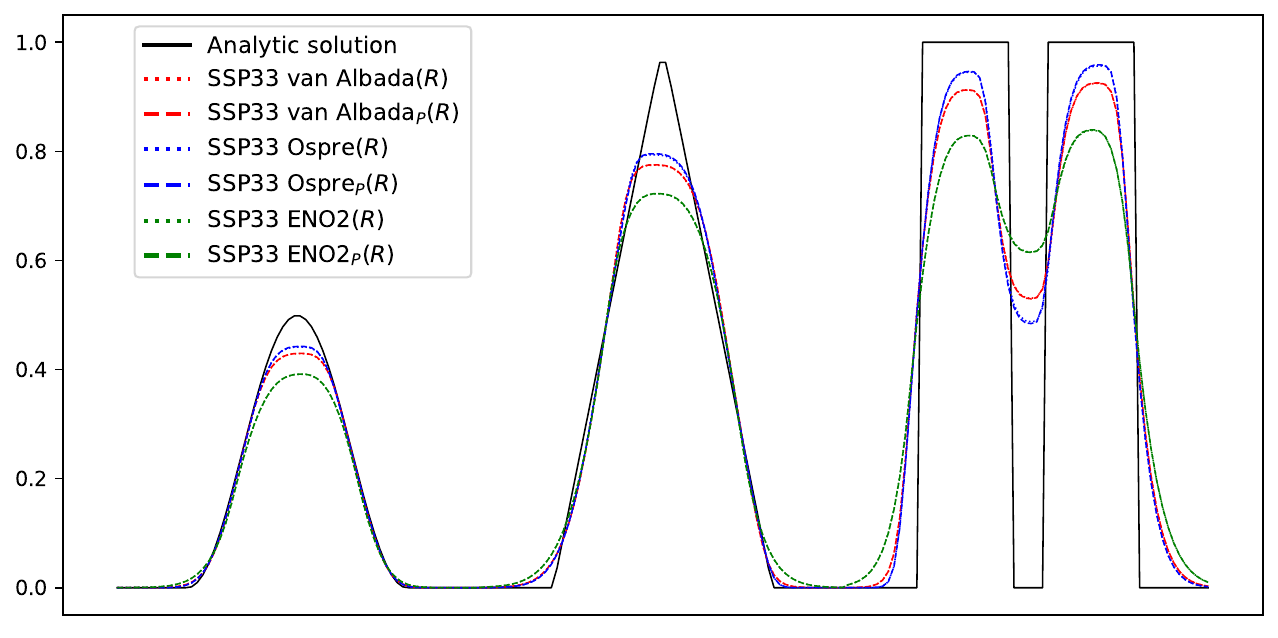}
\caption{Symmetric Spekreijse limiters and their Symmetric push into the Sweby region. We use the horizontal cross sections through the middle of the three LeVeque shapes; cosine, cone and the Zalesak slotted cylinder after the solid body rotation. Conclusion, all the pushed limiters are almost indistinguishable to the eye from the original limiters.}
\label{fig: Leveque slice ospre}
\end{figure}

\begin{figure}[H]
\centering
\includegraphics[scale=0.65]{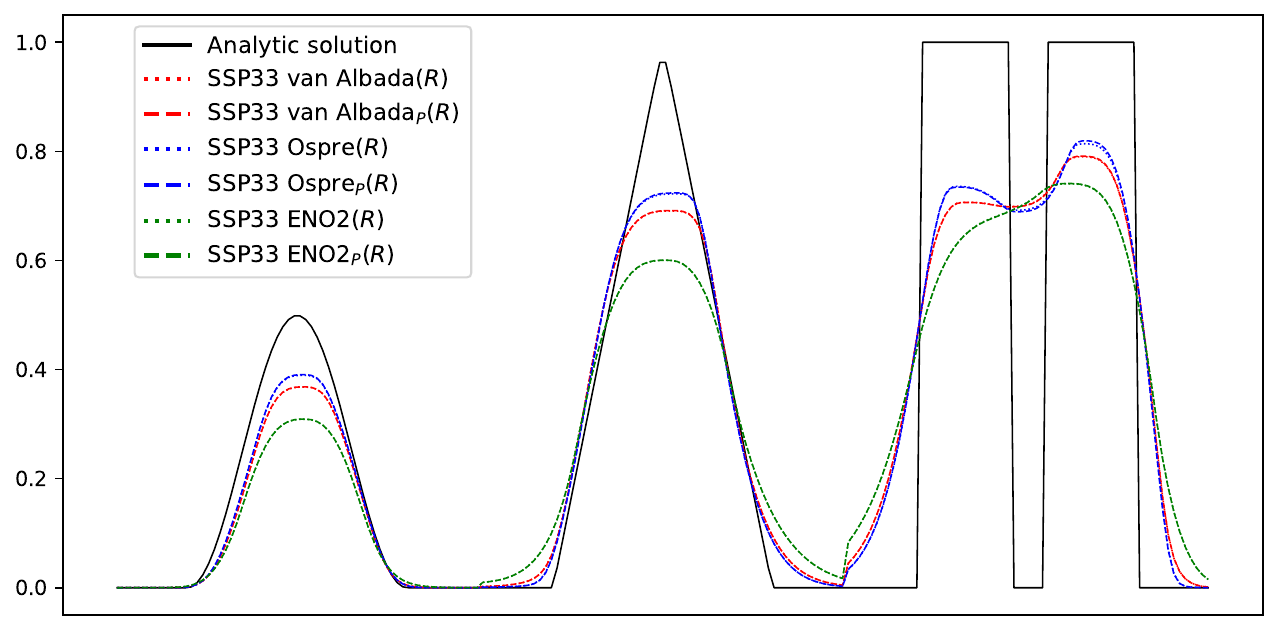}
\caption{Symmetric Spekreijse limiters and their Symmetric push into the Sweby region. We use the horizontal cross sections through the middle of the three LeVeque shapes; cosine, cone and the Zalesak slotted cylinder after the time reversing sinusoidal flow. Conclusion, all the pushed limiters are almost indistinguishable to the eye from the original limiters.}
\label{fig: Leveque slice ospre:def}
\end{figure}

\begin{figure}[H]
\centering
\includegraphics[scale=0.65]{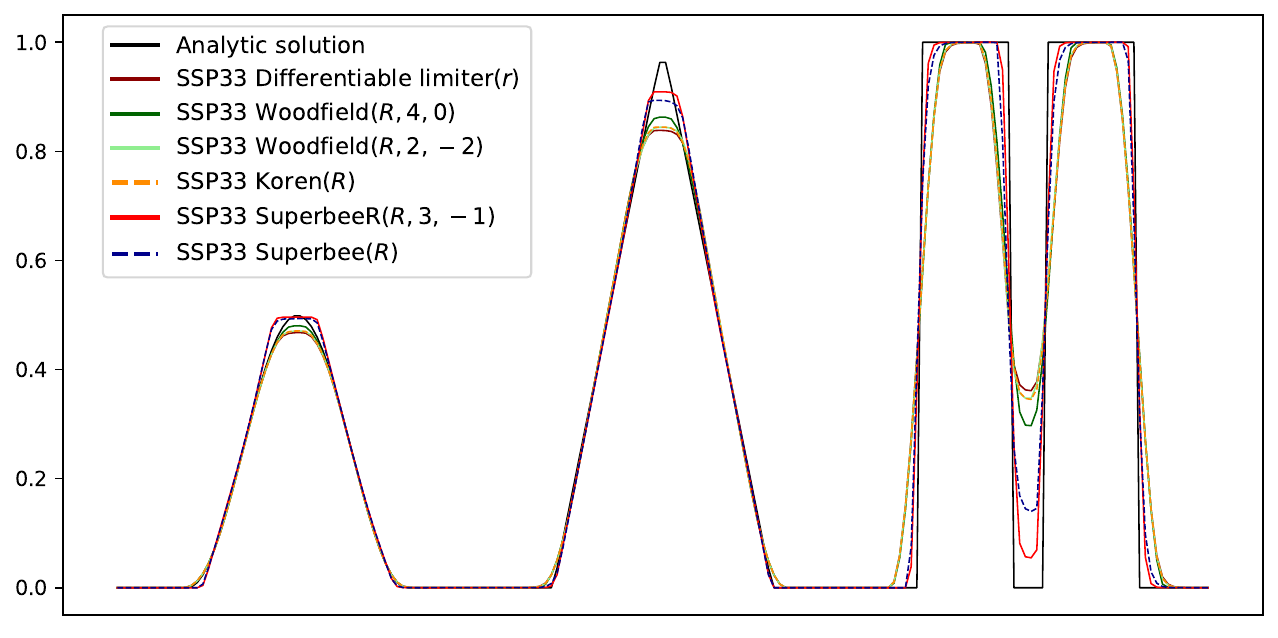}
\caption{Horizontal cross sections through the middle of the three LeVeque shapes; cosine, cone and the Zalesak slotted cylinder after the solid body rotation. The new differentiable limiter has results very similar to the Koren limiter. The differentiable limiter is more accurate than the Ospre limiter and is suitable for incompressible flow. The Woodfield(R,2,-2) limiter has similar accuracy to the Koren scheme, the Woodfield(R,4,0) limiter is more accurate than the Koren limiter, however it must be run at a reduced timestep. The SuperbeeR$(R,M,m)$ limiter is more compressive than the Superbee$(R)$ limiter.}
	\label{fig: Leveque slice koren}
\end{figure}

\begin{figure}[H]
\centering
\includegraphics[scale=0.65]{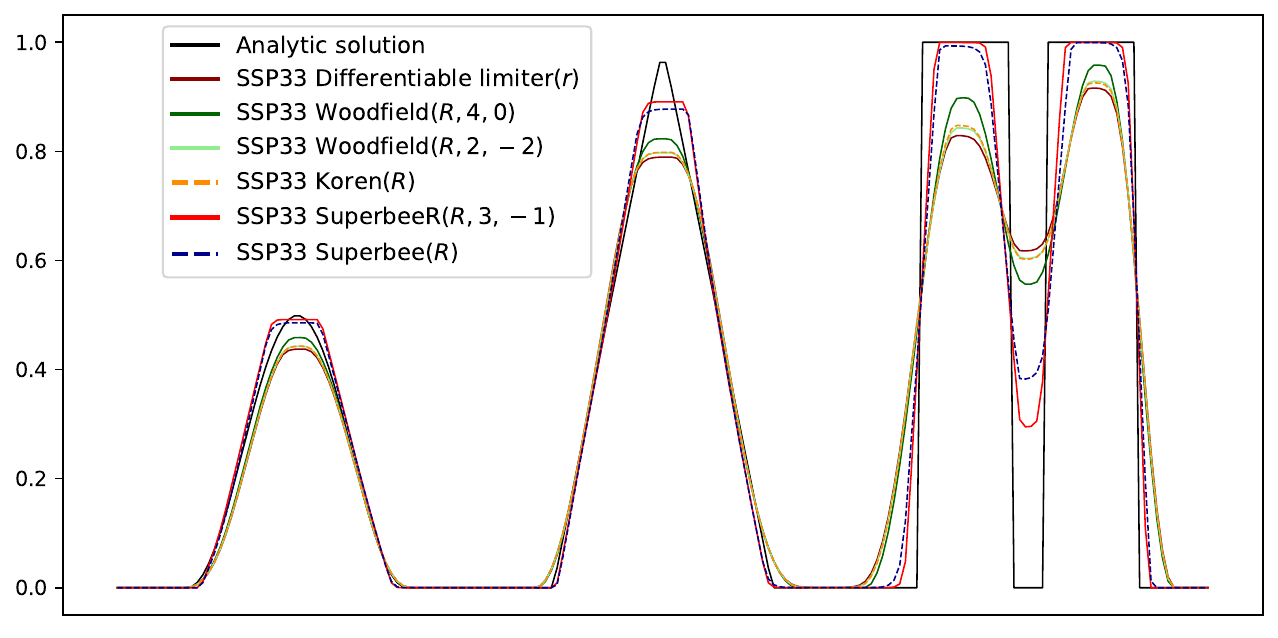}
\caption{Horizontal cross sections through the middle of the three LeVeque shapes; cosine, cone and the Zalesak slotted cylinder after the time reversing sinusoidal flow. The new differentiable limiter has results very similar to the Koren limiter. The differentiable limiter is more accurate than the Ospre limiter and is suitable for incompressible flow. The Woodfield(R,2,-2) limiter has similar accuracy to the Koren scheme, the Woodfield(R,4,0) limiter is more accurate than the Koren limiter, however it must be run at a reduced timestep. The SuperbeeR$(R,M,m)$ limiter is more compressive than the Superbee$(R)$ limiter.}
	\label{fig: Leveque slice koren:def}
\end{figure}

\Cref{fig: sbr on zal schemes} and \cref{fig: sinusoidal deformational reversing on zal schemes} contain the final tracer values of the  LeVeque initial conditions after the solid body rotation and the time reversing sinusoidal deformational flow respectively. \Cref{fig: Leveque slice ospre} and \cref{fig: Leveque slice ospre:def}, contain the horizontal cross section through the final tracer value of the middle of the three advected shapes, for the limiters (\cref{limiter: Eno2,limiter: Eno2P,limiter: Ospre,limiter: OspreP,limiter: Van Albada,limiter: Van Albaba P}) for the MVTS solid body rotation and the MVTV sinusoidal flow respectively. \Cref{fig: Leveque slice koren} and \cref{fig: Leveque slice koren:def}, contains the horizontal cross section through the final tracer value of the middle of the three advected shapes, for the other limiters including (\cref{limiter: differentiable new,limiter: superbeeR,limiter: Koren}), for the MVTS solid body rotation and the MVTV sinusoidal flow respectively.

The final tracer values of the  LeVeque initial conditions after the sinusoidal deformational time reversing flow plotted in \Cref{fig: sinusoidal deformational reversing on zal schemes}, contain negative values in the final time step for the schemes SSP33:ENO2$(R)$(fig c), SSP33:Ospre$(R)$(fig g) and SSP33:VanAlbada$(R)$(fig k). This observation is predicted (and proven in \cref{thm:necessary}) in \cref{sec:Applicability of the extended Spekreijse region} regarding the applicability of Spekreikse's limiter region to flux form incompressible advection. The RK4:Koren$(R)$ generated a $10^{-10}$ negative value over the integration period speculated as caused by the temporal integration as discussed in \cref{sec:RK4}. One observes sharper edges on the slotted cylinder when comparing SSP33:SuperbeeR$(R,M,m)$ (fig i) to SSP33:Superbee$(R)$ (fig h). One observes SSP33:Differentiable$(R)$ (fig l) is more visibly similar to the RK4:Koren$(R)$ scheme, than the other differentiable limiters SSP33:Ospre$(R)$ (fig g), SSP33:VanAlbada$(R)$ (fig k) which are seen to be successively more diffusive than the SSP33:Differentiable$(R)$ scheme, possibly due to the third order region being attained by this limiter as remarked on in \cref{eq:linear scheme pointwise truncation error}. 

 The final tracer values of the  LeVeque initial conditions after the solid body rotation of \Cref{fig: sbr on zal schemes} using the SSP33 timestepping scheme, contain no significant negative values in the final time step, this aligns with the theoretical predictions when combining by the MVTS property of the flow discussed in \cref{sec:Applicability of the extended Spekreijse region} and the SSP properties discussed in \cref{sec: temporal discretisation,sec:hidden maximum principle}. The RK4:Koren$(R)$ scheme generated a negative over the integration period as discussed in \cref{sec:RK4}. One observes sharper edges on the slotted cylinder when comparing SSP33:SuperbeeR$(R,M,m)$ (fig i) to SSP33:Superbee$(R)$(fig h). One observes the SSP33:Differentiable$(R)$ (fig l) limiter is visibly more similar to the RK4:Koren$(R)$ scheme, than the other differentiable limiters SSP33:Ospre$(R)$ (fig g), SSP33:VanAlbada$(R)$ (fig k) which are successively more diffusive.

\Cref{fig: Leveque slice ospre} and \cref{fig: Leveque slice ospre:def}, indicate that the push of the Ospre, VanAlbada and Eno2 limiters into the Sweby region had little effect on the accuracy of these schemes in the regimes of both sharp and smooth tracer values (for the timescales and flows in these particular experiments). The horisontal cross sections in \cref{fig: Leveque slice koren} and \cref{fig: Leveque slice koren:def} one observes that both the Superbee and the Superbee$(R,3,-1)$ have unphysical squaring of the cone and cosine lumps peak and more accurate representation of the discontinuous slotted cylinder profile. The Superbee$(R,M,m)$ limiter appears more compressive than Superbee. The horisontal cross sections in \cref{fig: Leveque slice koren} and \cref{fig: Leveque slice koren:def} indicate that the Differentiable limiter performs similar to the Koren$(R)$ limiter with a very small increase in diffusion due to the smoothing. The Woodfield$(R,2,-2)$ performed similar to the Koren$(R)$ limiter, indicating that in the negative gradient region in the $\theta=1$ region \cref{region:Woodfield} did not provide significant increase in accuracy under these test cases. The Woodfield$(R,4,0)$ performed more accurately than the Koren$(R)$ limiter, indicating additional accuracy in the $\theta=1$ region can be attained using the positive gradient region \cref{region:Woodfield}.


\subsubsection{Positivity}\label{sec:results:positivity}
\Cref{table: minimum values} contains the minimum values $\min_{\forall n,i,j} u_{i,j}^n$ attained for both the MVTV sinusoidal time reversing deformational flow \cref{test:sin deformation} as well as the MVTS solid body rotation \cref{test:solid body rotation} of the LeVeque initial conditions \cref{ic: Le veque}, for all limiters introduced in this paper. The main conclusion of \cref{table: minimum values} is that all limiters not strictly contained in the newly defined limiter regions \cref{region:Woodfield} have not remained positive for the sinusoidal MVTV flow. In particular, the Ospre$(R)$, VanAlbada$(R)$, and ENO2$(R)$, have positivity violations for deformational flow, demonstrating that the Spekreijse limiter region leads to schemes that are not bounded for the flux form advection equation. This is observed in sub-figures (b,g,k) of \cref{fig: sinusoidal deformational reversing on zal schemes}.

Also observed in \cref{table: minimum values,fig: sinusoidal deformational reversing on zal schemes}; limiters pushed into the Sweby region Ospre$_{P}(R)$, VanAlbada$_{P}(R)$, and ENO2$_{P}(R)$ remained strictly positive. The new limiters (Woodfield$(R,M,m)$, Differentiable$(R)$, SuperbeeR$(R,M,m)$) not contained in the Sweby region but contained in the newly defined limiter region, also remained strictly positive. RK4 is not positivity preserving. UTCDF$(R)$ is not positivity preserving, but UTCDF$_{S}(R)$, UTCDF$_{P}(R)$ are. In the second column of \cref{table: minimum values}, limiters in Spekreijse's limiter region remained positive for the mean value theorem satisfying flow, as predicted in \cref{sec:Applicability of the extended Spekreijse region}. All results in \cref{table: minimum values} agree with the theoretical contributions and predictions made in this paper up to machine precision. Column 3, is a different mean value theorem satisfying flow, for which significantly larger positivity violations are observed to emphasise the possibility of larger negatives.

\begin{table}[H]
\centering
\begin{tabular}{||c c c c c||}
\hline
Scheme&Limiter& MVTV-sin& MVTS-sbr & MVTV-sin 32 \\
& & $\min_{\forall n,i,j} u_{i,j}^n$& $\min_{\forall n,i,j} u_{i,j}^n$ &$\min_{\forall n,i,j} u_{i,j}^n$ \\ [0.5ex] 
\hline\hline

SSP33 & ENO2$(R)$ & $\b{-5.10824e-08}$ & 0.0 & $\b{-0.0128304}$   \\
SSP33 & Ospre$(R)$ & $\b{-5.21613e-11}$ & 0.0 & $\b{-0.0450631}$  \\
SSP33 & Vanalbada$(R)$ & $\b{-7.8611e-10}$ & 0.0 & $\b{-1.08423e-07}$   \\
\hline
SSP33 & ENO2$_P(R)$ & 0.0 & 0.0 & 0.0   \\

SSP33 & Ospre$_P(R)$ & 0.0 & 0.0 & 0.0   \\
SSP33 & Vanalbada$_P(R)$ & 0.0 & 0.0 & 0.0   \\
SSP33 & UTCDF$_P(R)$ & -1.07009e-56 & -1.55696e-57 & -7.39459e-34   \\
SSP33 & Koren$(R)$ & -8.73719e-19 & -4.61246e-19 & -3.62442e-19   \\
SSP33 & Superbee$(R)$ & -3.09042e-18 & -1.94581e-18 & -1.97676e-18   \\
\hline

SSP33 & Woodfield$(R,2,-1)$ & -8.91337e-19 & -5.86852e-19 & -4.15217e-19   \\
SSP33 & Woodfield$(R,4,0)$ & -1.0298e-18 & -7.05923e-19 & -4.22583e-19   \\
SSP33 & Differentiable$(r)$ & 0.0 & 0.0 & 0.0   \\

SSP33 & SuperbeeR$(R,3,-1)$ & -3.19545e-18 & -1.86983e-18 & -1.8609e-18  \\
SSP33 & UTCDF$_S(R)$ & -1.89143e-60 & -1.55696e-57 & -8.68458e-34   \\
\hline
RK4 & Koren$(R)$ & $\b{-9.80942e-11}$ & $\b{-2.55303e-10 }$& -3.72077e-19   \\
SSP33 & UTCDF$(R)$ & $\b{-5.70005e-05}$ & $\b{-7.31029e-05}$ & $\b{-0.000140264}$   \\
\hline
    \end{tabular}
	\caption{The minimum values attained in the full run of the experiments, bold indicates significant positivity violation. The Mean Value Theorem Violating sinusoidal flow (MVTV) is in the first and third} column, the Mean Value Theorem Satisfying (MVTS) rotational flow is in the second column. Only limiters not strictly contained in the new limiter regions \cref{region:Woodfield} are not positivity preserving for the sinusoidal flow. The RK4 scheme and the UTCDF limiter are not positivity preserving.  
	\label{table: minimum values}
\end{table}

\subsubsection{Accuracy}\label{sec:accuracy}
We compute the relative $l_2$ error norm $re_2(u,\Omega):= ||u - u_e||^2_{L^2(\Omega)}/||u_e||^2_{L^2(\Omega)}$ in subdomains of the Zalezak slotted cylinder $\Omega_{zal} = \lbrace (y,x) \in [5/8,7/8] \times [3/8,5/8] \rbrace $, the Cone $\Omega_{cone} = \lbrace (y,x) \in [3/8,5/8] \times [1/8,3/8]\rbrace$ and the Cosine lump $\Omega_{cos} = \lbrace (y,x) \in[3/8,5/8] \times [1/8,3/8] \rbrace$. \Cref{table: accuracy in different regions sbr} and \cref{table: accuracy in different regions sinusoidal} show the relative error norms at the final timestep in the regions $\Omega_{zal},\Omega_{cone},\Omega_{cos}$, for the solid body rotation and the sinusoidal flow respectively. The results indicate that the SuperbeeR$(R,3,-1)$ limiter has the smallest $L2$ error in the region of the slotted cylinder, also observed in the horizontal cross section \cref{fig: Leveque slice koren,fig: Leveque slice koren:def}. The Woodfield$(R,4,0)$ limiter has the smallest  $L2$ error in the region of the cone, and the cosine bump. The Woodfield$(R,4,0)$ scheme was more accurate than Koren limiter in all regions, which in turn was more accurate than the Woodfield$(R,2,-2)$ limiter in all regions. The Differentiable limiter was more accurate than the Ospre limiter which in turn was more accurate than the VanAlbada limiter observed in \cref{fig: Leveque slice ospre,fig: Leveque slice ospre:def,fig: Leveque slice koren,fig: Leveque slice koren:def}. Finally pushing the Ospre, VanAlbada and ENO2 limiters into the Sweby region, did not significantly decrease the relative error norm in any particular region, also observed in \cref{fig: Leveque slice ospre,fig: Leveque slice ospre:def}. This is interesting in light of the theoretical considerations in \cref{sec:accuracy of limited} since the ``pushed" schemes are reduced to the first order upwind scheme at extrema, in particular the ENO2$_{P}(R)$ no longer satisfies the second order at extrema condition \cite{hua1992possible} whilst ENO2$(R)$ does.

\begin{table}[H]
	\centering
	\begin{tabular}{||c c c c ||} 
	\hline
	Scheme & Slotted Cylinder & Cone & Cosine bell \\
	& $re_2(u,\Omega_{zal})$ & $re_2(u,\Omega_{cone})$ & $re_2(u,\Omega_{cos})$\\ [0.5ex] 
		\hline\hline
SSP33 Ospre$_P(R)$ & 0.351102 & 0.0673095 & 0.0658036 \\
SSP33 Ospre$(R)$ & 0.351987 & 0.0676014 & 0.0662881 \\
SSP33 van Albada$_P(R)$ & 0.371091 & 0.083057 & 0.076212 \\
SSP33 van Albada$(R)$ & 0.371456 & 0.0832153 & 0.0763878 \\
SSP33 ENO$2_{P}(R)$ & 0.422668 & 0.148316 & 0.106093 \\
SSP33 ENO$2(R)$ & 0.422894 & 0.148404 & 0.106193 \\
		\hline
SSP33 Woodfield$(R,4,0)$ & 0.274023 & 0.01206 & 0.0286216 \\
SSP33 Koren$(R)$ & 0.289004 & 0.0174684 & 0.0338112 \\
SSP33 Woodfield$(R,2,-2)$ & 0.290484 & 0.0181157 & 0.0342527 \\
SSP33 Differentiable limiter$(r)$ & 0.295057 & 0.0200713 & 0.0360145 \\
\hline
SSP33 SuperbeeR$(R,3,-1)$ & 0.168015 & 0.0496325 & 0.0341482 \\
SSP33 Superbee$(R)$ & 0.197235 & 0.0474037 & 0.0365515 \\
\hline
\end{tabular}
\caption{Relative $L_2$ error in different regimes of solid body rotation. }
	\label{table: accuracy in different regions sbr}
\end{table}

\begin{table}[h!]
	\centering
	\begin{tabular}{||c c c c ||} 
	\hline
	Scheme & Slotted Cylinder & Cone & Cosine bell \\
	& $re_2(u,\Omega_{zal})$ & $re_2(u,\Omega_{cone})$ & $re_2(u,\Omega_{cos})$\\ [0.5ex]  
		\hline\hline
SSP33 Ospre$_P(R)$ & 0.446784 & 0.177647 & 0.129118 \\
SSP33 Ospre$(R)$ & 0.447811 & 0.179705 & 0.13045 \\
SSP33 van Albada$_P(R)$ & 0.463153 & 0.219316 & 0.152391 \\
SSP33 van Albada$(R)$ & 0.463447 & 0.220061 & 0.152813 \\
SSP33 ENO$2_{P}(R)$ & 0.501686 & 0.34828 & 0.234589 \\
SSP33 ENO$2(R)$ & 0.501786 & 0.348622 & 0.234761 \\
\hline

SSP33 Woodfield$(R,4,0)$ & 0.373459 & 0.0762376 & 0.0640243 \\
SSP33 Koren$(R)$ & 0.390221 & 0.0931876 & 0.0768807 \\
SSP33 Woodfield$(R,2,-2)$ & 0.391025 & 0.0941322 & 0.0774227 \\
SSP33 Differentiable limiter$(r)$ & 0.396675 & 0.100286 & 0.0820969 \\
\hline
SSP33 SuperbeeR$(R,3,-1)$ & 0.279732 & 0.0779715 & 0.0718248 \\
SSP33 Superbee$(R)$ & 0.309864 & 0.0752206 & 0.0706076 \\
\hline
\end{tabular}
\caption{Relative $L_2$ error in different regimes for sinusoidal deformation. }
	\label{table: accuracy in different regions sinusoidal}
\end{table}

\section{Summary and Conclusions}

The main contribution of this paper is the derivation of two new limiter regions \cref{region:Woodfield}, larger than Sweby's admissible limiter region and smaller than the Spekreijse limiter region, sufficient for multidimensional incompressible flow to maintain a discrete local maximum principle (\cref{thm: divergence}). We show that the Spekreijse limiter region leads to flux form schemes that are not bounded for the advection equation unless a directional mean value theorem can be proven for each direction, this is demonstrated with a numerical example in \cref{table: minimum values} and by strictly necessary constraints in \cref{thm:necessary}.

In the new limiter regions, we have introduced; the first globally differentiable limiter entirely contained within the second order region capable of being monotone for flux form advection \cref{limiter: differentiable new}, a limiter \cref{limiter: Woodfield} that can be more accurate than the Koren limiter, and a limiter \cref{limiter: superbeeR} potentially more compressive than the Superbee limiter (more accurate in the regime of tracers with sharp gradients). 



The new limiter region provides a framework for proposing schemes that maintain a discrete local maximum principle when solving flux form equations when the flow is incompressible. We also have proven that this particular positive coefficient representation requires all symmetric limiters to be in the Sweby region \cref{sec:Details of proof}. We conjecture this to be a necessary condition for incompressible flow when symmetry is required on the limiter functions. We have also shown linear in-variance of the numerical scheme. Of practical merit we have numerically tested that pushing three known limiters in the Spekreijse region into the Sweby region had little negative consequence in terms of accuracy, we conjecture this to be related to the accuracy in the local neighbouhood of extrema considerations found in \cite{hua1992possible}.

Future work involves the study of the parameters $M,m$ in which the Woodfield$(R,M,m)$ limiter and the SuperbeeR$(R,M,m)$ limiter are most cost effective when the Courant number restriction is taken into account, and whether this changes when posed in the $\theta=0$ framework. The design of new limiters in the introduced region \cref{region:Woodfield} is open; accuracy, smoothness, and design criteria making use of less floating point operations would be of practical merit. We have only investigated $\theta = \lbrace 0,1\rbrace$, the generalisation by introducing a free parameter would be interesting and could lead to improved properties of limiting schemes. In the $\theta=0$ region, there are differentiable limiters in the second order region.
In the other region ($\theta=1$), there are no differentiable limiters in the second order region. One could conjecture that this asymmetry indicates that there may be less restrictive sufficient conditions in the $\theta=1$ framework than presented here using $\psi(R)<0$, for $R>0$. 

\section*{Acknowledgements}
During this work, the first author has been supported by an EPSRC studentship as part of the Centre for Doctoral Training in the Mathematics of Planet Earth (grant number EP/L016613/1).
Rupert Klein and Peter Sweby for comments leading to the improvement of this document. 
\bibliographystyle{abbrv}
\bibliography{BIBmain}

\appendix

\section{Details of proof for \cref{thm: divergence}}\label{sec:Details of proof}
\begin{proof} of \cref{thm: divergence}.
Expand both the method in \cref{sec:semi discrete scheme} 
\begin{align}
\begin{split}
\frac{\partial \bar{u}}{\partial t} 
& + c_{i+0.5,j}^{+}u^{R}_{i}
+ c_{i+0.5}^{-}u^{L}_{i+1}   
- c_{i-0.5}^{+}u^{R}_{i-1}
- c_{i-0.5}^{-}u^{L}_{i}   \\
& + c_{j+0.5}^{+}u^{U}_{j}
+ c_{j+0.5}^{-}u^{D}_{j+1} 
- c_{j-0.5}^{+}u^{U}_{j-1} 
- c_{j-0.5}^{-}u^{D}_{j} = 0, 
\end{split}\label{eq:method}
\end{align}
and a discrete form of the divergence-free condition
\begin{align}
\begin{split}
&c_{i+0.5}^{+}u+c_{i+0.5}^{-}u  -c_{i-0.5}^{+}u -c_{i-0.5}^{-}u \\
&+ c_{j+0.5}^{+}u+c_{j+0.5}^{-}u-c_{j-0.5}^{+}u -c_{j-0.5}^{-}u =0, 
\end{split}
\label{eq:divergence free}
\end{align}
in terms of their positive and negative components. Where here we adopt the notation where if one or other subscript is missing it is assumed to be at position $i$ or $j$ as appropriate. Taking away the divergence-free condition \eqref{eq:divergence free} from \eqref{eq:method} gives 
\begin{align}
\begin{split}
	\frac{\partial \bar{u}}{\partial t}& + c_{i+0.5,j}^{+}(u^{R}_{i}-u_{i})+c_{i+0.5}^{-}(u^{L}_{i+1} -u_{i}) -c_{i-0.5}^{+}(u^{R}_{i-1}-u_{i}) -c_{i-0.5}^{-}(u^{L}_{i} -u_{i})\\& + c_{j+0.5}^{+}(u^{U}_{j}-u_{i}) +c_{j+0.5}^{-}(u^{D}_{j+1}-u_{i}) -c_{j-0.5}^{+}(u^{U}_{j-1}-u_{i})  -c_{j-0.5}^{-}(u^{D}_{j}-u_{i}) = 0. \label{eq:removed div free}
\end{split}
\end{align}

Using $u_{i+1}-u_{i}= R_{i}(u_{i}-u_{i-1})$ and \cref{eq:UR}, one can derive the following identities
\begin{align}
    u_i^R - u_i &=  \frac{\theta}{2} \psi(R_i)(u_i - u_{i-1}) + \frac{1-\theta}{2} \psi(\frac{1}{R_{i}})(u_{i+1}-u_{i}),\\
 & = \big[ \frac{\theta}{2} \psi(R_i) + \frac{1-\theta}{2} R_{i}\psi(\frac{1}{R_{i}}) \big] (u_{i}-u_{i-1}).\label{eq:URmu}\\
& = \big[ \frac{\theta}{2} \frac{\psi(R_i)}{R_i} + \frac{1-\theta}{2}\psi(\frac{1}{R_{i}}) \big] (u_{i+1}-u_{i}).
\end{align}

using $u_{i-1}-u_{i-2} = 
\frac{1}{R_{i-1}}(u_i-u_{i-1})$ and \cref{eq:UR}, one can derive the following identities
\begin{align}
    u_{i-1}^R -u_i &= (u_{i-1} -u_i)+ \frac{\theta}{2} \psi(R_{i-1})(u_{i-1} - u_{i-2}) + \frac{1-\theta}{2} \psi(\frac{1}{R_{i-1}})(u_{i}-u_{i-1}),\\
    &= - \big[ 1 -  \frac{\theta}{2} \frac{\psi(R_{i-1})}{R_{i-1}} - \frac{1-\theta}{2} \psi(\frac{1}{R_{i-1}})\big] (u_{i}-u_{i-1}),\label{eq:URm1mu}\\
    &= - \big[ 1 -  \frac{\theta}{2} \frac{\psi(R_{i-1})}{R_{i-1}} \big] (u_{i}-u_{i-1}) +  \frac{1-\theta}{2} \psi(\frac{1}{R_{i-1}})\frac{1}{R_{i}}(u_{i+1}-u_{i}),
\end{align}

Using $-R_i(u_{i}-u_{i-1}) = (u_{i}-u_{i+1})$ and  and \cref{eq:UL}, one can derive the identities
\begin{align}
    u_i^L -u_i  &=  \frac{\theta}{2} \psi(\frac{1}{R_{i}})(u_i - u_{i+1}) - \frac{1-\theta}{2} \psi(R_{i})(u_{i}-u_{i-1}),\\
    &= \big[ -\frac{\theta}{2} \psi(\frac{1}{R_{i}})R_i  - \frac{1-\theta}{2} \psi(R_{i}) \big] (u_{i}-u_{i-1}),\\
&= \big[\frac{\theta}{2} \psi(\frac{1}{R_{i}}) + \frac{1-\theta}{2} \frac{\psi(R_{i})}{R_i} \big] (u_{i}-u_{i+1}),\label{eq:ULmu} 
\end{align}

Using $-R_{i+1}(u_{i+1}-u_{i}) =(u_{i+1}-u_{i+2})$ and \cref{eq:UL}, one can derive the identities
\begin{align}
    u_{i+1}^L -u_i  &=  (u_{i+1} - u_i) + \frac{\theta}{2} \psi(\frac{1}{R_{i+1}})(u_{i+1} - u_{i+2}) - \frac{1-\theta}{2} \psi(R_{i+1})(u_{i+1}-u_{i}),\\
    &=  (u_{i+1} - u_i) + -\frac{\theta}{2} \psi(\frac{1}{R_{i+1}})R_{i+1}(u_{i+1} - u_{i}) - \frac{1-\theta}{2} \psi(R_{i+1})(u_{i+1}-u_{i}),\\
    &= -\big[1  -\frac{\theta}{2} \psi(\frac{1}{R_{i+1}})R_{i+1} - \frac{1-\theta}{2} \psi(R_{i+1})\big] (u_{i}-u_{i+1}), \label{eq:ULp1mu}
\end{align}

As highlighted above, there are many achievable representations of the numerical scheme that can be attained using only the definition of $R$. In this work, we seek a representation with less stringent requirements on the limiter function for the representation to be of positive coefficient type.

The idea behind the representation is to group $u^{R}_{i}-u_i$, and $u^{R}_{i-1}-u_i$, together in terms of $(u_i-u_{i-1})$ only, using \cref{eq:URmu,eq:URm1mu}. And to also group $u^L_{i+1}-u_i$, $u^L_{i}-u_i$ together in terms of $u_{i+1}-u_i$, using \cref{eq:ULmu,eq:ULp1mu}.

This way we group like signed terms such as $c_{i+1/2}^{+}$, and $c_{i-1/2}^{+}$ together in the same positive coefficient. This, when substituted into \cref{eq:removed div free} results in the following representation 
	\begin{align}
	\begin{split}
		\frac{\partial u}{\partial t} &+ \bigg( c_{i+1/2}^{+}\big[ \frac{\theta \psi( R_{i} ) }{2}  + \frac{(1-\theta)}{2} R_{i}\psi( \frac{1}{R_{i}} ) \big] + c_{i-1/2}^{+}\big[ 1- \frac{\theta \psi(R_{i-1})}{2 R_{i-1}} - \frac{(1-\theta)}{2}\psi(\frac{1}{R_{i-1}}) \big] \bigg)[u_{i}-u_{i-1}]\\
		& + \bigg( -c^{-}_{i+1/2}\big[ 1-\frac{\theta}{2}\psi(\frac{1}{R_{i+1}})R_{i+1} - \frac{1-\theta}{2}\psi(R_{i+1}) \big] - c_{i-1/2}^{-}\big[ \frac{\theta}{2}\psi(\frac{1}{R_{i}}) + \frac{1-\theta}{2} \frac{\psi(R_{i})}{R_{i}}\big] \bigg)[u_{i}-u_{i+1}]\\
		& + \text{y-direction}. 
	\end{split}
	\end{align}
 This representation is fundamental to the construction of the limiter frameworks, introduced in this paper. But many different representations of the scheme are possible, and it may be the case that another representation can be proven positive coefficient, perhaps with less demands on the limiter function.

When $\theta = 1$ the representation reduces to 
\begin{align}
	\begin{split}
\frac{\partial u}{\partial t} 
&+ \bigg[ c_{i+1/2}^{+}\frac{\psi(R_{i}) }{2}+ c_{i-1/2}^{+}[1-\frac{\psi(R_{i-1})}{2R_{i-1}}] \bigg] (u_{i} - u_{i-1}) \\
&+\bigg[ -c^{-}_{i+1/2}[1- \frac{1}{2}R_{i+1}\psi(\frac{1}{R_{i+1}})] - c^{-}_{i-1/2} \frac{1}{2} \psi(\frac{1}{R_{i}}) \bigg] (u_{i} - u_{i+1}) + \text{y-direction} = 0,
\end{split}
\end{align}
where we read off the leading term coefficients
\begin{align}
	A_{i-1/2} &=\bigg(   \frac{1}{2} \bigg[ c_{i+0.5}^{+}\psi(R_{i}) \bigg] +c_{i-0.5}^{+} \bigg[1-\frac{1}{2}\psi(R_{i-1})/R_{i-1} \bigg] \bigg), \\
	B_{i+1/2} &= -\bigg( c_{i+0.5}^{-} \bigg[ 1 -\frac{1}{2}\psi(\frac{1}{R_{i+1}})R_{i+1} \bigg] + \frac{1}{2} \bigg[ c_{i-0.5}^{-}\psi(R_{i}^{-1}) \bigg] \bigg), \\
	C_{j-1/2} &= (i \mapsto j )\circ (A_{i-1/2}) , \quad D_{j+1/2} = (i \mapsto j )\circ (B_{i+1/2} ). 
\end{align}
We have introduced unorthodox notation $(i \mapsto j )\circ (f(i))$ to mean the same expression but with $i$ replaced with $j$.
Clearly 
\begin{align}
\psi(R) \geq 0  \quad \text{and} \quad \psi(S)/S \leq 2,
\end{align} 
are sufficient, for the scheme to be of positive coefficient type. 

It is also imposed as a necessary assumption on the limiter functions to be of this particular positive coefficient representation because of the arbitrary nature of the velocity field. Suppose that $c_{i+1/2}\geq 0$, then for the scheme to be of positive coefficient type at both $u_{i,j}$ and $u_{i+1,j}$ we must require that,
\begin{align}
	1/2c_{i+1/2}^+\psi(R_{i}) \geq 0, \quad 	c_{i+1/2}^+(1-1/2\psi(R_{i})R_{i}^{-1}) \geq 0.
\end{align}
 There could always exists a different positive coefficient representation of the scheme under some yet to be found transform or rearrangement.

We now attempt to find a sensible sufficient timestep restriction, the conditions
\begin{align}
	&0 \leq \psi(R)  \leq M_{\psi}, &\quad m_{\psi} \leq \psi(S)/S \leq 2, \\
	& m_{\psi} \leq \psi(1/T)T \leq 2, &\quad  0 \leq \psi(1/r) \leq M_{\psi}, \quad \forall R,S,T \in \mathbb{R},
\end{align} are sufficient for the following bounds
\begin{align}
\begin{split}
	A_{i-1/2} &\in [0, c_{i+0.5}^{+}M_{\psi}/2 + c_{i-0.5}^{+} (1-m_{\psi}/2) ] , \\
	B_{i+1/2} & \in [0 , -c_{i+0.5}^{-}[1 - m_{\psi}/2 ]-c_{i-0.5}^{-} M_{\psi}/2] ,\\
	C_{j-1/2} &\in [0, c_{j+0.5}^{+}M_{\psi}/2 + c_{j-0.5}^{+} (1- m_{\psi}/2) ] , \\
	D_{j+1/2} & \in [0 , -c_{j+0.5}^{-}[1 - m_{\psi}/2 ]-c_{j-0.5}^{-} M_{\psi}/2],
\end{split}
\end{align}
where we are yet to define the constants $ M_{\psi}\geq 0$, and $m_{\psi} \leq 2$.

The time step restriction is 
\begin{align}
	&	A_{i-1/2} + B_{i+1/2} + C_{i-1/2} + D_{i+1/2} \leq 1, 
\end{align}
which can be satisfied when
\begin{align}
\begin{split}
	&c_{i+0.5}^{+}M/2 + c_{i-0.5}^{+} (1-m/2) - c_{i+0.5}^{-}(1 -m/2)- c_{i-0.5}^{-} M/2 \\
	&+c_{j+0.5}^{+}M/2 + c_{j-0.5}^{+} (1-m/2) - c_{j+0.5}^{-}(1 -m/2)-c_{j-0.5}^{-} M/2  \leq 1 .
\end{split}
\end{align}
We lose some generality for a more convenient sufficient time step restriction. Define flow in and out Courant numbers by the following definitions
\begin{align}
	C_{i,j}^{in}&:= c_{i+0.5}^{+}-c_{i-0.5}^{-}+ c_{j+0.5}^{+} - c_{j-0.5}^{-} ,\\
	C_{i,j}^{out}&:= -c_{i+0.5}^{-} + c_{i-0.5}^{+} - c_{j+0.5}^{-} + c_{j-0.5}^{+}.
\end{align}
where these definitions are chosen based on how the flows effect the solution $du/dt$.
The time step restriction can be written as
\begin{align}
	C^{in}_{i,j} \frac{M_{\psi}}{2} + C_{i,j}^{out}(1- \frac{m_{\psi}}{2}) \leq 1,
\end{align}
using incompressibility $C_{i,j}^{in} = C_{i,j}^{out}$ we can write this as the following Courant number restriction
\begin{align}
	C \leq C_{FE} =  \frac{1}{(1+ \frac{M_{\psi} - m_{\psi}}{2}) }.
\end{align}
So far we have considered the $\theta = 1$ case, the general form is given by
	\begin{align}
	\begin{split}
		\frac{\partial u}{\partial t} &+ \bigg( c_{i+1/2}^{+}\big[ \frac{\theta \psi( R_{i} ) }{2}  + \frac{(1-\theta)}{2} R_{i}\psi( \frac{1}{R_{i}} ) \big] + c_{i-1/2}^{+}\big[ 1- \frac{\theta \psi(R_{i-1})}{2 R_{i-1}} - \frac{(1-\theta)}{2}\psi(\frac{1}{R_{i-1}}) \big] \bigg)[u_{i}-u_{i-1}]\\
		& + \bigg( -c^{-}_{i+1/2}\big[ 1-\frac{\theta}{2}\psi(\frac{1}{R_{i+1}})R_{i+1} - \frac{1-\theta}{2}\psi(R_{i+1}) \big] - c_{i-1/2}^{-}\big[ \frac{\theta}{2}\psi(\frac{1}{R_{i}}) + \frac{1-\theta}{2} \frac{\psi(R_{i})}{R_{i}}\big] \bigg)[u_{i}-u_{i+1}]\\
		& + \text{y-direction}. 
	\end{split}
	\end{align}
We now repeat the previous argument for $\theta =0$,
omitting some of the details. The below conditions
	\begin{align}
		\frac{1}{2} R_{i}\psi(\frac{1}{R_{i}}) & \geq 0, \\
		1-\frac{1}{2}\psi(\frac{1}{R_{i-1}})&\geq 0, \\
		1-\frac{1}{2}\psi(R_{i+1}) & \geq 0, \\
		\frac{\psi(R_{i})}{2R_{i}} & \geq 0, 
	\end{align}
 are sufficient for the positivity of the coefficients. This can be written more conveniently as
\begin{align}
	m_{\psi}\leq \psi(1/R) \leq 2, \quad 0 \leq S\psi(\frac{1}{S}) \leq M_{\psi}.
\end{align}
The timestep restriction for the discrete maximum principle for the $\theta =0$ scheme is
\begin{align}
	C\leq \frac{2}{2+M_{\psi}-m_{\psi}}.
\end{align}
\end{proof}

\section{Necessary conditions for sign preservation}
We detail some of the necessary conditions on the limiter function deduced by case-by-case considerations on possible velocity field and tracer values. This is done through means of contradiction-type arguments to attain the necessary conditions for sign preservation. 
\begin{proof}[Necessary]\label{proof:necessary}
Consider the forward Euler scheme given by \begin{align}
	\begin{split}
	u^{n+1}_{ij}=u_{ij}^n &- \bigg( c_{i+1/2}^{+}\big[ \frac{\theta \psi( R_{i} ) }{2}  + \frac{(1-\theta)}{2} R_{i}\psi( \frac{1}{R_{i}} ) \big] + c_{i-1/2}^{+}\big[ 1- \frac{\theta \psi(R_{i-1})}{2 R_{i-1}} - \frac{(1-\theta)}{2}\psi(\frac{1}{R_{i-1}}) \big] \bigg)[u_{i}-u_{i-1}]\\
		& - \bigg( -c^{-}_{i+1/2}\big[ 1-\frac{\theta}{2}\psi(\frac{1}{R_{i+1}})R_{i+1} - \frac{1-\theta}{2}\psi(R_{i+1}) \big] - c_{i-1/2}^{-}\big[ \frac{\theta}{2}\psi(\frac{1}{R_{i}}) + \frac{1-\theta}{2} \frac{\psi(R_{i})}{R_{i}}\big] \bigg)[u_{i}-u_{i+1}]\\
		&- \bigg( c_{j+1/2}^{+}\big[ \frac{\theta \psi( R_{j} ) }{2}  + \frac{(1-\theta)}{2} R_{j}\psi( \frac{1}{R_{j}} ) \big] + c_{j-1/2}^{+}\big[ 1- \frac{\theta \psi(R_{j-1})}{2 R_{j-1}} - \frac{(1-\theta)}{2}\psi(\frac{1}{R_{j-1}}) \big] \bigg)[u_{j}-u_{j-1}]\\
		& - \bigg( -c^{-}_{j+1/2}\big[ 1-\frac{\theta}{2}\psi(\frac{1}{R_{j+1}})R_{j+1} - \frac{1-\theta}{2}\psi(R_{j+1}) \big] - c_{j-1/2}^{-}\big[ \frac{\theta}{2}\psi(\frac{1}{R_{j}}) + \frac{1-\theta}{2} \frac{\psi(R_{j})}{R_{j}}\big] \bigg)[u_{j}-u_{j+1}].
	\end{split}
	\end{align}
Let $\theta=1$, $c_{i+1/2}>0$, $c_{i-1/2}=0$, and let $c_{j+1/2}$, $c_{j-1/2}$ be free, such that the divergence free assumption holds. Then let $u_{j+1}=0$, $u_{j-1}=0$, $u_{ij}^n=0$,  and then let
$u^n_{i-1}>0$ and let $u^n_{i+1}\geq 0$. The forward Euler scheme reduces to \begin{align}
u^{n+1}_{i,j} =  1/2c_{i+1/2}^{+}\psi(R_i^n) u^n_{i-1}\quad \text{where}\quad R_{i} = -\frac{u_{i+1}^n}{u_{i-1}^n} \in (-\infty,0]. 
\end{align}
If one supposes there exists $R<0$, such that $\psi(R)<0$, one can find $u_{i+1}$, $u_{i-1}$ such that $u^{n+1}_{ij}<0$ a change in sign. Therefore, $\psi(R)\geq 0$ is a necessary assumption when $R<0$ and $\theta=1$.

Now let $\theta=1$, $c_{i+1/2}=0$, $c_{i-1/2}>0$, and let $c_{j+1/2}$, $c_{j-1/2}$ be free, such that the divergence free assumption holds. Then let $u_{j+1}=0$, $u_{j-1}=0$, $u_{ij}^n=0$ and then let
$u^n_{i-1}>0$ and let $u^n_{i-1}\geq 0$. The forward Euler scheme reduces to \begin{align}
u^{n+1}_{i,j} = c_{i-1/2}^{+} \left(1- \frac{\psi(R_{i-1}^n)}{2R_{i-1}^n}\right) u^n_{i-1}\quad \text{where}\quad R_{i-1} = \frac{-u_{i-1}^n}{u_{i-1}^n - u_{i-2}^n} \in (-\infty,-1] \cup (0,\infty). 
\end{align}
 Therefore as can be proven through contradiction, $\psi(R)/R\leq 2$ is a necessary assumption for sign preservation for $\theta=1$ on the interval $R\in (-\infty,-1] \cup (0,\infty)$. 

Now consider $c_{i+1/2}<0$, $c_{i-1/2}=0$, and $u_{ij}^n=0$. Let $c_{j+1/2}$, $c_{j-1/2}$, be free such that the divergence free condition holds. Let $u_{j}=u_{j-1}=u_{j+1}$, then the scheme reduces to 
\begin{align}
    u^{n+1}_{ij} = (-c_{i-1/2}^{-})\left[ 1-\frac{1}{2}\psi(\frac{1}{R_{i+1}}) R_{i+1}\right]u_{i+1} \quad \text{where}\quad R_{i+1} = \frac{u_{i+2}-u_{i+1}}{u_{i+1}}\in [-1,\infty).
\end{align}
Therefore, using the intermediate variable $S_{i+1}=1/R_{i+1}\in (-\infty,-1)\cup (0,\infty)$, one can through contradiction show that 
\begin{align}
    \frac{\psi(S)}{S}\leq 2,
\end{align}
is a necessary assumption. Under similar considerations where instead $c_{i+1/2}=0$, $c_{i-1/2}<0$ one can show that $\psi(R)\geq 0$ for $R\in (-\infty,0)$, is a necessary assumption. There is likely more additional necessary assumptions required on the limiter function, that can be derived from additional case-by-case considerations. 
\end{proof}


\section{Linear invariance}\label{sec:Linear invariance}
In this section we prove linear invariance of the numerical scheme. 
\begin{definition}[Linear invariant scheme]
	The map $u_{i}\mapsto \alpha w_{i}+ \beta$ leaves the numerical method unchanged.
\end{definition}

\begin{theorem}[Linear invariance]
	The method described in \cref{sec:semi discrete scheme} is linear invariant for incompressible flow.
\end{theorem}

\begin{proof}
As in \cite{hundsdorfer1995positive} the transform $u_{i}\mapsto \alpha w_{i}+ \beta$ is investigated for linear invariance. However, we do the analysis after the use of a discrete divergence free condition,
\begin{align}
	\frac{\partial u}{\partial t}
	&+ \bigg(   \frac{1}{2} \bigg[ c_{i+0.5}^{+}\psi(R_{i}) \bigg] +c_{i-0.5}^{+} \bigg[1-\frac{1}{2}\psi(R_{i-1})/R_{i-1} \bigg] \bigg) (u_{i} - u_{i-1})\\
	& -\bigg( c_{i+0.5}^{-} \bigg[ 1 -\frac{1}{2}\psi(\frac{1}{R_{i+1}})R_{i+1} \bigg] + \frac{1}{2} \bigg[ c_{i-0.5}^{-}\psi(R_{i}^{-1}) \bigg] \bigg)( u_{i}-u_{i+1} )\\
	& +\text{same expression with } i\mapsto j = 0,
\end{align}
and observe
\begin{align}
	R_{i} =\frac{	 (u_{i+1} -  u_{i} )}{(u_{i} -u_{i-1}  )}  &\mapsto \frac{	(\alpha w_{i+1} + \beta )- (	\alpha w_{i} + \beta )}{(\alpha w_{i} + \beta )- (	\alpha w_{i-1} + \beta) } = \frac{	 (w_{i+1} -  w_{i} )}{(w_{i} -w_{i-1}  )} ,\\
	u_{i} - u_{i-1} &\mapsto \alpha (w_{i}-w_{i-1}), \\
	\frac{\partial u}{\partial t} \mapsto  \alpha \frac{\partial w}{\partial t},
\end{align}
gives exact linear invariance of the numerical method. Constants are preserved exactly for incompressible flow, and the equation is scaling invariant. This holds for both $\theta\in \lbrace 0,1 \rbrace$.
\end{proof}

\section{Time step restrictions}
The numerical scheme defined in \cref{sec:semi discrete scheme} has a time step restriction for a discrete maximum principle that depends on which limiter is used, whether a directional mean value theorem holds and what temporal discretisation is used. We collate \cref{table: cfl} of the maximum cell defined Courant number the Forward Euler scheme has a discrete maximum principle predicted by \cref{thm: divergence}, for some of the limiters introduced in this paper. This table can be used to check whether the scheme should be producing positive or discrete maximum principle satisfying results. The Courant number depends on the parameters $M$ and $m$, this effects the efficiency of the scheme, and can be made arbitrarily small for injudicious choice of limiter.
\begin{table}[h!] 
	\centering
	\begin{tabular}{||c c c c||  } 
		\hline
		Scheme&Limiter & \Cref{thm: divergence}-$C$&  Spekreijse-MVTS-$C$    \\ [0.5ex] 
		\hline\hline
		FE& $\psi(R)=0$ & 1 & 1 \\
		FE& $\psi(R)=aR+b$ & 0 & 0\\
		FE&vanAlbada$(R)$ & 0 & $2-\sqrt{2}$ \\
		FE&Ospre$(R)$ & $0$& $1/2$\\
		FE&ENO2$(R)$ & 0 & 1/2\\ 
		FE&vanAlbada$_P(R)$ & $\frac{4}{5+\sqrt{2}}$ & $\frac{4}{5+\sqrt{2}}$ \\
		FE&Ospre$_P(R)$ & $4/7$ & $4/7$\\
		FE&minmod$(R)$/ENO2$_P(R)$ & $2/3$ & $2/3$ \\
		FE&Koren$(R)$& $1/$2 & $1/2$\\ 
		FE&Woodfield(R,M,m)& $2/(2 +M - m)$ & $2/(2 +M - m)$\\ 
		FE&Differentiable(r)& $\frac{2}{4 + \sqrt{ \frac{5\sqrt{5} }{2} -\frac{11}{2} }}\approx0.4651$ & $\frac{2}{4 + \sqrt{ \frac{5\sqrt{5} }{2} -\frac{11}{2} }}\approx0.4651$\\ 
		RK4&Koren$(R)$& 0 & 0		\\
        FE & UTCDF$(R)$ & 0 & 0 \\
        FE & UTCDF$_S(R)$ & 1/2 & 1/2 \\
        FE & UTCDF$_P(R)$ & 1/2 & 1/2 \\
		\hline
	\end{tabular}
	\caption{The sufficient Courant number restrictions for temporal discretisations of the method \cref{sec:semi discrete scheme} to satisfy a local discrete maximum principle for different limiters.  Column one is the sufficient Courant number limit for incompressible flow and column two is the sufficient Courant number limit for directional mean value theorem satisfying flows. These time stepping criteria hold in more dimensions; however the definition of the Courant number changes.}
	\label{table: cfl}
\end{table}

	\section{Hidden maximum principles from internal strong stability}\label{sec:hidden maximum principle}
 
	When one requires that the substages satisfy a discrete maximum principle, and limits at every stage in the Shu Osher representation, one implies substage maximum principles. For example our SSP33 scheme satisfies the (hidden) maximum principles
	\begin{align}
		k^1 &\in[ \min \lbrace u_{i+1}^n,u_{i-1}^n,u^n,u_{j+1}^n,u_{j-1}^n \rbrace, \max \lbrace u_{i+1}^n,u_{i-1}^n,u^n,u_{j+1}^n,u_{j-1}^n \rbrace ], \\
		k^2 &\in[ \min \lbrace k^1_{i+1}, k^1_{i-1}, k^1, k^1_{j+1}, k^1_{j-1} \rbrace, \max \lbrace k^1_{i+1}, k^1_{i-1}, k^1, k^1_{j+1}, k^1_{j-1} \rbrace ], \\
		u^{n+1} & \in[ \min \lbrace k^2_{i+1}, k^2_{i-1}, k^2, k^2_{j+1}, k^2_{j-1} \rbrace, \max \lbrace  k^2_{i+1}, k^2_{i-1}, k^2, k^2_{j+1}, k^2_{j-1} \rbrace ],
	\end{align}
	instead of the original maximum principle 
	\begin{align}
		u^{n+1}_{i,j} &\in[ \min \lbrace u_{i+1}^n,u_{i-1}^n,u^n,u_{j+1}^n,u_{j-1}^n \rbrace, \max \lbrace u_{i+1}^n,u_{i-1}^n,u^n,u_{j+1}^n,u_{j-1}^n \rbrace ].
	\end{align}
	This implies a different maximum principle, 
	\begin{align}
		u^{n+1}_{i,j} & \in[ \min_{ \forall k,l \in \lbrace -3,-2,-1,0,1,2,3 \rbrace } u_{i+k,j+l} , \max_{ \forall k,l \in \lbrace -3,-2,-1,0,1,2,3 \rbrace } u_{i+k,j+l} ],
	\end{align}
	but enforces positivity and local boundedness of the substages, which could be essential in the context of a numerical method for a coupled system.

\section{Runge Kutta 4}\label{sec:RK4}
 The standard Runge Kutta 4 method well established for both its accuracy and efficiency is proposed for this particular spatial scheme by both \cite{koren1993robust} and \cite{hundsdorfer1995positive}. In \cite{koren1993robust} it is indicated monotonicity is supposed to be guaranteed by small timesteps, it has since been shown that the RK4 algorithm has no Shu Osher representation (without perturbation techniques and downwind/upwind biasing) and has no(zero) radius of monotonicity. Hundsdorfer et al. \cite{hundsdorfer1995positive} notes this but indicates the SSP literature is of little practical importance. In \cite{hundsdorfer1994method}, the same author performed similar experiments and quantifies the observation of $10^{-6}$ negative values \cite{hundsdorfer1994method}. From a practical perspective, the importance of using a SSP timestepping method depends on the range of acceptable monotonicity failure. These negatives would be considered not just significant but very large in our application where the positivity of density should not be negative. We consider roughly $10^{-14}$ the threshold of acceptable monotonicity failure, due to machine precision error in constructing a divergence free vector field, and other accumulation of machine precision errors at 64 bit arithmetic. At lower precision arithmetic, it may be the case that the SSP literature is of less practical importance. \Cref{table: minimum values} contains the minimum values attained by RK4 with the Koren$(R)$ limiter they are positivitity violating.

\section{Accuracy of un-limited scheme}\label{sec:accuracy of un-limited}
 Hundsdorfer et al. \cite{hundsdorfer1995positive} viewed the semi-discrete method \cref{sec:semi discrete scheme} as a finite difference method evolving point-wise values $u_{i,j}$ solving the singular point value equation. Hundsdorfer, and Trompert \cite{hundsdorfer1994method} showed that the third order upwind method solves the following modified equation using backward error analysis
\begin{align}
	\begin{split}
		u_t + (v_1 u)_x  &= - \frac{\Delta x^2}{24}\big[  u (v_{1})_{xxx} + 3 u_x (v_{1})_{xx} + 2 u_{xx} (v_{1})_x \big]   + (\text{y-direction}),\\
  &- \frac{\Delta x^3}{12}[u_{xxxx}v_1+u_{xxx}(v_1)_x] + O(\Delta x^4) + (\text{y-direction}).
	\end{split}\label{eq: finite difference error}
\end{align}
Where $v_1$ is the $x$ component of the velocity field. So the scheme is a formally third order finite difference method approximating the singular point equation in one and two dimensions for directionally constant flow.  
 
Zijlema and Wesseling \cite{zijlema1998higher} proposed the scheme as a finite volume method on cell averaged quantities approximating the integral form of the equation. When the method is interpreted as a finite volume scheme approximating the integral form of the equation, the formal truncation error analysis becomes different see \cite{nishikawa2021truncation,leonard1995order,nishikawa2021quick} for explanation. In which case in one dimension the finite volume interpretation of this scheme becomes truly third order with respect to the integral form of the equation (Aleksandar Donev has uploaded lecture notes proving this fact using the symbolic algebra package Mathematica \cite{Donev}). However, the scheme becomes second order in two dimensions or more because the one-dimensional flux integral is approximated by second order Gauss quadrature. 

The formal truncation analysis of the full method is tedious in the more general setting, and instead one designs limiters using simplifications. If it is assumed the flow is uniformly constant, and the equation is considered in the pointwise finite difference sense. Then the linear scheme associated with
\begin{align}
	\psi(R) = aR+b,  \label{eq:linear:unlimited}
\end{align}
has the semi-discrete truncation error 
\begin{align}
		\frac{u^{R}_{i}-u_{i-1}^{R}}{\Delta x} - (u_x)_{i} =  [a+b-1]\frac{\Delta x u_{xx} }{2!}  + [1-3b]\frac{\Delta x^2 u_{xxx} }{3!}  +(a+7b-1)\frac{\Delta x^3 u_{xxxx} }{4!}+... . \label{eq:linear scheme pointwise truncation error}
\end{align}
	So that $u\in C^2$ and $a+b=1$ is sufficient for second order, such that $\psi(R) = aR+ (1-a)$ passes through $(1,1)$. $u\in C^3$, $a=2/3, b=1/3$ is sufficient for third order. 

\section{Accuracy of limited scheme}\label{sec:accuracy of limited}
Spekreijse, showed that sufficient conditions on the limiter function $\psi(R)$ for second order accuracy are $\psi(1) = 1$, $\psi \in C^{2}$ near 1, \cite{S_1987}. Sufficient and necessary conditions on the limiter function $\psi(r)$ were shown by \cite{hua1992possible} to be $\psi(1) = 1$, and $\psi$ is Lipschitz continuous.

We also want to understand the effect limiting has on the truncation error particularly in the neighbourhood of smooth extrema where the gradient changes. It is often mistaken that TVD schemes are first order at noncritical extrema because of a slight oversight in the truncation analysis in Osher 1984 \cite{osher1984high},\cite{barth2003finite}. However, this is not quite the case as explained by Hua-mo \cite{hua1992possible} where formal truncation analysis of the $\theta = 0$ scheme is done by expanding the higher order correction in the neighbourhood of a critical point $x_{\alpha} = x_{i}+ \alpha \Delta x$. We summarise the results from \cite{hua1992possible} below.

\begin{theorem}[constant flow in one dimension near the critical points $\theta = 0$ \cite{hua1992possible}.]
	TVD schemes $\theta = 0$ may have second order accuracy at critical points if $\psi(r=3) + \psi(r=-1) = 2$. But cannot be uniformly second-order accurate in the neighbourhood of critical points. If $\psi(1) = 1$, then these TVD schemes have second-order accuracy in the region sufficiently far from the critical points of smooth solutions. 
	\begin{proof}
		\cite{hua1992possible}.
	\end{proof}
\end{theorem}

Similar conclusions for the $\theta = 1$ form is done in \cite{zijlema1998higher} however both theorem 3.1 and proof of theorem 3.1, have a few technical inaccuracies. Nevertheless, all the arguments are correct in spirit, second order at extrema is technically possible but impossible to achieve in the local neighbourhood of extrema.  We recommend following the analysis methods of \cite{hua1992possible} instead.

\begin{theorem}[\cite{hua1992possible}-for $\theta = 1$]
 Let $\tau_i$ denote the truncation error about point $x_{i}$, let $x_{\xi} = x_{i} + \xi \Delta x, \quad \xi \in [0,\infty)$, be a point of interest. 
 TVD schemes may have second order accuracy at critical points if  $3 \psi(R=1/3) - \psi(R=-1) = 2.$ But cannot be uniformly second-order accurate in the whole neighbourhood of critical points. If $\psi(1) = 1$, then TVD schemes have second-order accuracy in the region sufficiently far from the critical points of smooth solutions.
\end{theorem}
\begin{proof} See \cite{hua1992possible}, we include some sketch details of the proof in \cref{Near Extrema}.
\end{proof}
The reason this theoretical work is revisited by us here, is because some authors suggest using the second order at extrema condition in the design of limiters \cite{hua1992possible,zijlema1998higher} and others \cite{waterson1995unified} suggest that satisfying the theoretical second order at extrema condition does not have practical consequences. In our work we have suggested various modifications to limiters including pushing into the Sweby region, this breaks the $3\psi(R=1/3) - \psi(R=-1) = 2$ condition for the ENO2 limiter.

\section{On the order of total variation diminishing schemes near extrema}\label{Near Extrema}
In this section we present a sketch of the arguements of \cite{hua1992possible} but for the $\theta = 0$ case. We use the prime notation $u'$ to denote the spatial derivative $u_{x}$ for readability. The semi-discrete finite difference local truncation error $\tau$ of the spatial derivative $u_x$ at position $x_i$ is 
	\begin{align}
		\begin{split}
		\tau_i &:= \frac{u(x_i) - u(x_{i-1})}{\Delta x} + \frac{\psi(R_u(x_{i+1},x_{i},x_{i-1}))(u(x_{i}) - u(x_{i-1}))}{2 \Delta x} - \frac{\psi(R_u(x_{i},x_{i-1},x_{i-2}))(u(x_{i-1}) - u(x_{i-2}))}{2 \Delta x} - u'(x_i), \label{EQ:truncation error}
		\end{split}
	\end{align}
	where $R_u(a,b,c):= \frac{u(a) - u(b)}{u(b) -u(c)}$. As in \cite{hua1992possible} we wish to write this local truncation error at $x_{i}$ as a function of $u$ about an arbitrary point $x_{\xi} = x_{i}+\xi \Delta x$, to determine the order of accuracy within a distance of a critical point. 
	
First Taylor expand the first order upwind term in \cref{EQ:truncation error} about $x_i$ to cancel with $u_{x}(x_i)$, then re-expanding again about the arbitrary point $x_{\xi} = x_{i}+\xi \Delta x$ as in \cite{hua1992possible}. One can deduce the first order upwind scheme has the following expression of local truncation error
	\begin{align}
		\frac{u(x_i) - u(x_{i-1})}{\Delta x} -u'(x_i) &= - \Delta x \frac{u''(x_i)}{2!}+ \Delta x^2 \frac{u'''(x_i)}{3!} ...\\
		&=- \frac{\Delta x}{2!} [u''(x_{\xi}) - \xi \Delta x u'''(x_{\xi}) +\frac{\xi^2 \Delta x^2 u''''(x_{\xi})}{2!} ... ]\\
		&+  \frac{\Delta x^2}{3!}[u'''(x_{\xi}) - \xi \Delta x u''''(x_{\xi}) ...]\\
		& = - \frac{\Delta x}{2!} u''(x_{\xi})  +  (3\xi+1) \frac{\Delta x^2 u'''(x_{\xi})}{3!} + O(\Delta x^3).
	\end{align}

Now for the higher order flux limited correction in \cref{EQ:truncation error} use the following Taylor expansions about $x_{\xi}$
	\begin{align}
		u(x_{i+1}) - u(x_{i}) &=  \Delta x u'(x_{\xi}) -(\xi- 1/2) \Delta x^2 u''(x_{\xi}) -\frac{(\xi-1)^3 - \xi^3}{3!} \Delta x^3 u'''(x_{\xi}) -... \\
		u(x_{i}) - u(x_{i-1}) &= \Delta x u'(x_{\xi}) - (\xi+1/2) \Delta x^2 u''(x_{\xi}) - \frac{(\xi)^3 - (\xi+1)^3}{3!} \Delta x^3 u'''(x_{\xi}) -... \\
		u(x_{i-1}) - u(x_{i-2}) &= \Delta x u'(x_{\xi}) - (\xi+3/2) \Delta x^2 u''(x_{\xi}) - \frac{(\xi+1)^3 - (\xi+2)^3}{3!} \Delta x^3 u'''(x_{\xi}) -... 
	\end{align}
This allows us to write down an expression for the local truncation error in terms of $x_{\xi}$
	\begin{align}
		\begin{split}
		\tau_i &:=  - \frac{\Delta x}{2!} u''(x_{\xi})  +  (3\xi+1) \frac{\Delta x^2 u'''(x_{\xi})}{3!} \\
		& +\frac{1}{2}\psi(R_i)[ u'(x_{\xi}) -(\xi+ 1/2) \Delta x u''(x_{\xi}) -\frac{(\xi-1)^3 - \xi^3}{3!} \Delta x^2 u'''(x_{\xi})] \\
		& -\frac{1}{2}\psi(R_{i-1})[  u'(x_{\xi}) - (\xi+3/2) \Delta x u''(x_{\xi}) - \frac{(\xi+1)^3 - (\xi+2)^3}{3!} \Delta x^2 u'''(x_{\xi}) ] + O(\Delta x^3). 
		\end{split}
	\end{align}
	We collect terms of the appropriate order to define the truncation error as
	\begin{align}
		\begin{split}
		\tau_i &:= \frac{1}{2}[\psi(R_i)-\psi(R_{i-1})] u'(x_{\xi})  \\
		& + [ -1 - \psi(R_i)(\xi+1/2)+\psi(R_{i-1})(\xi+3/2)] \frac{\Delta x}{2!} u''(x_{\xi})\\
		&+\bigg[ (3\xi+1) -\frac{1}{2}\psi(R_i) [(\xi-1)^3 - \xi^3]+ \frac{1}{2}\psi(R_{i-1})[(\xi+1)^3 - (\xi+2)^3]\bigg]\frac{\Delta x^2 u'''(x_{\xi})}{3!} ,
		\end{split}
	\end{align}
here the gradient $R$ is also Taylor expanded about the arbitrary point $x_{\xi}$,
	\begin{align}
	R_{i} & = \frac{   u'(x_{\xi}) -(\xi- 1/2) \Delta x u''(x_{\xi}) -\frac{(\xi-1)^3 - \xi^3}{3!} \Delta x^2 u'''(x_{\xi}) -...}{u'(x_{\xi}) - (\xi+1/2) \Delta x u''(x_{\xi}) - \frac{(\xi)^3 - (\xi+1)^3}{3!} \Delta x^2 u'''(x_{\xi}) -... }\\
	R_{i-1} & = \frac{u'(x_{\xi}) - (\xi+1/2) \Delta x u''(x_{\xi}) - \frac{(\xi)^3 - (\xi+1)^3}{3!} \Delta x^2 u'''(x_{\xi}) -... }{ u'(x_{\xi}) - (\xi+3/2) \Delta x u''(x_{\xi}) - \frac{(\xi+1)^3 - (\xi+2)^3}{3!} \Delta x^2 u'''(x_{\xi}) -...}.
\end{align}

	This allows us to state the requirements of the scheme to be first second and third order. 
	\begin{definition}[First order requirement near $x_{\xi}$]
		First order accuracy requires
		\begin{align}
			\big[  \frac{\psi(R_{i})}{2}  -\frac{\psi(R_{i-1})}{2}  \big] u_{x}(x_{\xi})  = O(\Delta x).
		\end{align}
	\end{definition}
	
	\begin{definition}[Second order requirement near $x_{\xi}$]
		Second order accuracy requires
		\begin{align}
			\begin{split}
			&\big[  \frac{\psi(R_{i})}{2}  -\frac{\psi(R_{i-1})}{2}  \big] u_{x}(x_{\xi})  \\
			+& \big[ -1 - (\xi +\frac{1}{2} )\psi(R_{i})+ (\xi + \frac{3}{2})\psi(R_{i-1})\big] \frac{\Delta x u_{xx}(x_{\xi}) }{2!} = O(\Delta x^2).
			\end{split}
		\end{align}
	\end{definition}
	
	\begin{definition}[Third order order requirement near $x_{\xi}$]
		Third order accuracy requires
		\begin{align}
			\begin{split}
			&\big[  \frac{\psi(R_{i})}{2}  -\frac{\psi(R_{i-1})}{2}  \big] u_{x}(x_{\xi})  \\
			+& \big[ -1 - (\xi +\frac{1}{2} )\psi(R_{i})+ (\xi + \frac{3}{2})\psi(R_{i-1})\big] \frac{\Delta x u_{xx}(x_{\xi}) }{2!} \\
			+& \big[ (3\xi +1) - \frac{\psi(R_{i})}{2}[(\xi)^3 - (\xi+1)^3] + [(\xi+1)^3 - (\xi+2)^3]\frac{\psi(R_{i-1})}{2} \big] \frac{ \Delta x^2 (u_{xxx})(x_{\xi})}{3!} = O(\Delta x^3).
			\end{split}
		\end{align}
	\end{definition}

	There are now three cases to consider. We want to know the error at a critical point $x_{i}$, at a non critical point $x_{i}$, and the error at $x_{i}$ within the vicinity of a critical point $x_{\xi}$. Where we do not consider inflections as critical. One has to make an assumption on the size of $\Delta x$, in order to use the identity $1/(1-a) = 1 + a+ a^2 ...$, for $|a|\leq 1$, to approximate the expressions for $R_i,R_{i-1}$. We will also use the fact that $O(|f(x)|)=O(f(x))$, and the definition of Lipshitz continuity $\forall x,y, \exists K>0 \text{s.t.} |f(x)-f(y)|<K|x-y|$. One also has to make technical remarks about the smoothness of $u$ required for derivatives or more sharply, for the mean value theorem to apply (the inclusion of these technical remarks would be a important contribution to the literature).
	Case 1: We are at a critical point in position $i$, so that 
		\begin{align}
			\xi &= 0\\
			u_{x}(x_{\xi}) &=  u_{x}(x_{i}) = 0\\
			R_{i} &=-1 - \frac{2\Delta x}{3}  \frac{u_{xxx}(x_{i})}{u_{xx}(x_{i})} +O(\Delta x^2)\\
			R_{i-1} &=\frac{1}{3}+ \frac{2\Delta x}{3}  \frac{u_{xxx}(x_{i})}{u_{xx}(x_{i})}+O(\Delta x^2)
		\end{align}
		where $u_{xx}(x_i)\neq 0$.
		Case 2: We are at a non critical point i, so that
		\begin{align}
			\xi &= 0\\
			u_{x}(x_{\xi}) &=  u_{x}(x_{i})\neq 0\\
			R_{i} &= 1 + \Delta x u''(x_{i}) /u'(x_{i}) +O(\Delta x^3) \\
			R_{i-1} &=1 + \Delta x u''(x_{i}) /u'(x_{i}) + \Delta x^2 (4/3u_{xxx}x_{i} -3u_{xx}(x_{i})^2 )+O(\Delta x^3)
		\end{align}
		Case 3: We are at $x_i$ in the vicinity of a critical point at $x_{\xi}$
		\begin{align}
			\xi  \neq 0\\
			u_{x}(x_{\xi}) &= 0\\
			R_{i} &=\frac{(\xi - 1/2)}{(\xi+ 1/2)} + O(\Delta x) \\
			R_{i-1} &=\frac{(\xi + 1/2)}{(\xi + 3/2)} + O(\Delta x) 
		\end{align}
	
We satisfy first order accuracy requirements in case 1 and case 3 trivially because $u_x(x_{\xi}) = 0$.  In case 2 at a non-critical $x_{i}$ the first order requirement is 
\begin{align}
			\psi(1 + \Delta x u''(x_{i}) /u'(x_{i}) +O(\Delta x^3)) - \psi(1+ \Delta x u''(x_{i}) /u'(x_{i}) + O(\Delta x^2) ) = O(\Delta x). 
\end{align}
 A sufficient condition for this to be satisfied is Lipshitz continuity of the limiter function.

We first note that all inflection points will be second order because $u_x(x_{i}), u_{xx}(x_{i})= 0$. In case 1 at critical, non inflection points we have second order accuracy when 
\begin{align}
	-1 -1/2\psi(-1 - a \Delta x +O(\Delta x^2))+3/2 \psi(1/3+a \Delta x + O(x^2) ) = O(\Delta x), \label{eq:hurry}
\end{align}
where $a = \frac{2}{3}  \frac{u_{xxx}(x_{i})}{u_{xx}(x_{i})}$. If one first assumes 
\begin{align}
	1 -3/2\psi(1/3) + 1/2\psi(-1) = 0
\end{align}
and then adds this to \cref{eq:hurry}. One can deduce that Lipschitz continuity and $1 -3/2\psi(1/3) + 1/2\psi(-1) = 0$ are sufficient and necessary for second order, by the use of a contradiction arguement for the necessary condition. 

In case 2 at a non critical, non inflection point we have second order accuracy when both the conditions
		\begin{align}
			\psi(1 + a \Delta x +O(\Delta x^3)) - \psi(1+a \Delta x + O(\Delta x^2 ) = O(\Delta x^2) \\ 
			-1 -1/2\psi(1 + a \Delta x +O(\Delta x^3))+3/2 \psi(1+a \Delta x + b\Delta x^2 ) = O(\Delta x) 
		\end{align}
		are satisfied for $a= u''(x_{i}) /u'(x_{i})$. A sufficient and necessary condition for the first condition to hold is Lipschitz continuity. If this assumption is used on the second condition we can rewrite the second condition as 
		\begin{align}
			-1 + \psi(1+a \Delta x + O(\Delta x^2 ) = O(\Delta x) \label{eq:lipshitz above}
		\end{align}
		 If one assumes $\psi(1)=1$ and adds $1-\psi(1)=0$ to the above expression \cref{eq:lipshitz above} one can see that Lipschitz continuity and $\psi(1)=1$ are sufficient and necessary.
		We now consider second order accuracy in case 3, within a neighbourhood of a critical point which requires
		\begin{align}
			-\big[ 1 - ( \xi +\frac{1}{2} )\psi(R_{i})- (\xi + \frac{3}{2})\psi(R_{i-1})\big] \frac{u_{xx}(x_{\xi})}{2!} = O(\Delta x).
		\end{align}
		Leading to the condition 
		\begin{align}
			\big[ 1 + (\frac{1}{2} + \xi )\psi(\frac{\xi - 1/2}{\xi+ 1/2} + O(\Delta x))- (\xi + \frac{3}{2})\psi(\frac{\xi + 1/2}{\xi + 3/2} + O(\Delta x) )\big] = O(\Delta x),\label{eq: lip xi}
		\end{align}
		assuming differentiability of the limiter then a sufficient condition for the above to hold is the following condition
		\begin{align}
			1 + ( \xi +\frac{1}{2} )\psi(\frac{\xi - 1/2}{\xi+ 1/2})- (\xi + \frac{3}{2})\psi(\frac{\xi + 1/2}{\xi + 3/2})  = 0 \label{eq:A37}
		\end{align}
	If the above property \cref{eq:A37} is assumed and it is subtracted from \cref{eq: lip xi}, the differentiability is relaxed to the Lipschitz continuity and can be proven both necessary and sufficient.
	
	The second order in the neighbourhood of an extrema condition \cref{eq:A37} is satisfied for any second order linear scheme $\psi(R) = aR+b$ when $1-a-b=0$. But the only solution that satisfies $\psi(0)=0$ is given by $\psi(R) = aR$ which is not within monotonicity constraints for large gradients \cite{hua1992possible}. This \cref{eq:A37} condition reduces to the second order $\psi(1) = 1$ non-extrema condition in the limit $\xi \rightarrow \pm \infty$ sufficiently away from extrema. The \cref{eq:A37} condition reduces to the previously derived second order at extrema condition $2 = 3\psi(1/3) - \psi(-1)$ in the limit $\xi \rightarrow 0$. Concluding the known accuracy results in the neighbourhood of extrema.

 \section{Relationship to other limiter regions.}
\label{sec:Relationship to other limiters}

This appendix reviews some other limiter regions, the limiter regions reviewed are different in definition, construction and aims than those presented in this paper. Spekreijse in 1987 \cite{S_1987} introduced generalised extended limiter regions using the $\theta = 0$, framework allowing for non zero values of $\psi(R)$ for $R<0$. This same type of free parameter generalisation was later applied to the Sweby region in the $\theta=1$ framework in the work by Cahouet-Coquel in 1989
\cite{jacques1989uniformly} see also \cite{hua1992possible}, in which non zero values are allowed for $\psi(r)$ for $r<0$. It is worth noting that in both $\theta \in \lbrace 0,1\rbrace$ frameworks, a mean value theorem is always imposed, these extended limiter regions have not been proven positive coefficient for flux form advection for general velocity fields. This is due to their prevalent use in flux splitting or flux vector splitting frameworks.
In 2010, Dubois \cite{dubois2010nonlinear} introduced two limiter regions in the context of one-dimensional flows, based on the monotonicity and convexity of the subcell reconstruction. Hyman in 1983 \cite{hyman1983accurate} and Huynh in 1993 \cite{huynh1993accurate} extended the monotonicity preserving constraints in which the MP-Parabolic construction or a Hermite cubic interpolant is monotone. When written in terms of a limiter function, this is found to be $\psi \in [-\min(3,3R), \min(3,3R)]$, $0\leq \psi(R)\leq 3\min(1,R)$ respectively. Čada and Torrilhon, in \cite{vcada2009compact} design limiter functions based of the regions defined by both Hyman and Huynh. Typically, limiters deriving from Hyman's work have local monotone properties from the reconstruction, but not formal strict maximum principles (similar to WENO), unless they are further restricted to the Spekreijse limiter region. The above limiter regions are different to those in this paper, in terms of construction, aim and definition, in particular, the extent to which the above limiter regions are proven to have a discrete maximum principle for flux form incompressible advection has not been extablished.

\subsection{Error handling}
Here we introduce another table, using the same test case used for \cref{table: minimum values}, but with simpler error handling (without conditionals handling case by case considerations in division).

\begin{table}[H]
    \centering
    \begin{tabular}{|| c c c c c ||}
\hline
		Scheme&Limiter& MVTV-sin& MVTS-sbr & MVTV-sin 32 \\
		& & $\min_{\forall n,i,j} u_{i,j}^n$& $\min_{\forall n,i,j} u_{i,j}^n$ &$\min_{\forall n,i,j} u_{i,j}^n$ \\ [0.5ex] 
		\hline\hline

SSP33 & ENO2$(R)$ & $\b{-5.10824e-08}$ & -1.63505e-16 & $\b{-0.0128308}$   \\
SSP33 & Ospre$(R)$ & $\b{-5.21618e-11}$ & -4.08984e-16 & $\b{-0.0450631}$   \\
SSP33 & Vanalbada$(R)$ & $\b{-7.8611e-10}$ & -3.69424e-16 & $\b{-1.08423e-07}$   \\
\hline
SSP33 & ENO2$_P(R)$ & -2.52212e-16 & -1.76984e-16 & -7.11199e-17   \\
SSP33 & Ospre$_P(R)$ & -7.20499e-16 & -5.87407e-16 & -8.58698e-17   \\
SSP33 & Vanalbada$_P(R)$ & -3.53086e-16 & -4.21421e-16 & -7.9581e-17   \\
SSP33 & UTCDF$_P(R)$ & -1.15881e-15 & -7.18115e-16 & -2.96113e-16   \\
SSP33 & Koren$(R)$ & -7.50797e-16 & -9.29711e-16 & -4.70821e-15   \\
SSP33 & Superbee$(R)$ & -7.03817e-16 & -7.35427e-16 & -1.18831e-14   \\
\hline
SSP33 & SuperbeeR$(R,3,-1)$ & -1.16908e-15 & -8.03216e-16 & -9.65484e-16   \\

SSP33 & Woodfield$(R,2,-1)$ & -7.55775e-16 & -1.15396e-15 & -8.68959e-16   \\
SSP33 & Woodfield$(R,4,0)$ & -8.93315e-16 & -8.78768e-16 & -4.14971e-15   \\
SSP33 & Differentiable$(r)$ & -5.86937e-16 & -2.30743e-16 & -1.05367e-14   \\
SSP33 & UTCDF$_S(R)$ & -1.26445e-15 & -7.46022e-16 & -2.3749e-16   \\
\hline
RK4 & Koren$(R)$ & $\b{-9.80942e-11}$ & $\b{-2.55303e-10}$ & -4.6105e-15   \\
SSP33 & UTCDF$(R)$ & $\b{-5.70005e-05}$ & $\b{-7.31029e-05}$ & $\b{-0.000140264}$  \\
\hline\hline
    \end{tabular}
    \caption{Same as \cref{table: minimum values}, but using simple error handling for division, division is defined as $DV(a,b):= (a/(b+2.2e-16))$. We observe the same results up to machine precision.}
    \label{tab:my_label}
\end{table}

\setcitestyle{square}


\end{document}